\documentclass[trsc,nonblindrev]{informs3noheader} 
\OneAndAHalfSpacedXI 

\usepackage[authoryear]{natbib}

\usepackage{tikz}
\usepackage{fix-cm}
\usepackage[edges]{forest}
\usepackage{multicol}
\usepackage{enumitem}
\usepackage{csquotes}
\usepackage{float}
\usepackage{needspace}
\usepackage{orcidlink}
\usepackage[pagewise]{lineno}
\usepackage{booktabs}
\usepackage{makecell}
\usepackage{eurosym}
\usepackage{mathrsfs}
\usepackage{xcolor}
\providecolor{scolor}{gray}{0}
\usepackage{etoolbox}
\usepackage{xpatch}
\usepackage{xspace}
\usepackage{bbding}
\usepackage{bm}
\usepackage{algorithm}
\usepackage{algorithmic}
\usepackage{subcaption}

\providecommand{\sep}{;\space}

\hypersetup{citecolor=blue, linkcolor=blue, urlcolor=blue}

\def\tsc#1{\csdef{#1}{\textsc{\lowercase{#1}}\xspace}}
\tsc{WGM}
\tsc{QE}

\makeatletter
\def\changeBibColor#1{%
  \in@{#1}{}%
  \ifin@\color{blue}\else\normalcolor\fi
}

\xpatchcmd{\@bibitem}
  {\item}
  {\changeBibColor{#1}\item}
  {}
  {\PackageWarning{mypatch}{Failed to patch \string\@bibitem}}

\xpatchcmd{\@lbibitem}
  {\item}
  {\changeBibColor{#2}\item}
  {}
  {\PackageWarning{mypatch}{Failed to patch \string\@lbibitem}}
\makeatother

\begin{document}
\let\WriteBookmarks\relax
\def\floatpagepagefraction{1}
\def\textpagefraction{.001}

\RUNTITLE{Reinforcement Learning and Operational Research}  


\TITLE{Reinforcement Learning and Operational Research: A Review of Methodological Integration and Research Directions}

\ARTICLEAUTHORS{%
\AUTHOR{Yahan Lu$^{a}$, Dongyang Xia$^{a}$, Nursen Aydin$^{b,*}$, Shadi {Sharif Azadeh}$^{a}$}
\AFF{$^a$ Department of Transport \& Planning, Delft University of Technology, The Netherlands \\ 
$^b$ Warwick Business School, University of Warwick, UK\\
$*$: Corresponding author\\
}
}
\ABSTRACT{
The growing demand for real-time, data-driven decision-making in complex and dynamic systems is placing increasing pressure on Operational Research (OR) methodologies. Sequential decision-making under uncertainty has long been studied in OR through Markov decision processes, dynamic programming, and their approximations, and reinforcement learning (RL) provides simulation-based methods for learning decision policies in such problems. Recent research shows an increasing interest in integrating RL with OR to address dynamic decision-making problems and enhance heuristic and exact methods for combinatorial optimization. The overarching goal across these efforts is to improve solution quality, computational efficiency, and robustness. Given the diversity of integration approaches and application settings, there is a clear need for a systematic and technically detailed review of how RL empowers OR methods. To address this gap, this paper presents a structured review of two methodological roles of RL in OR. The first concerns sequential decision problems in dynamic environments, where states evolve and decisions are repeated over time. The second concerns static combinatorial optimization problems, where inputs are fixed before the solution process begins and RL either constructs solutions or supports heuristic and exact algorithms. Digital twins provide an optional deployment and evaluation layer for training, testing, and validating these methods. We critically synthesize recent advances across these roles and uses, highlighting their advantages, implementation requirements, limitations, and challenges. Finally, we consolidate guidance on when the reviewed approaches may be considered and outline directions for future research on the methodological and practical integration of RL and OR. 
}


\maketitle

\section{Introduction}
\label{Sec_introduction}


Operational Research (OR) provides a foundational framework for supporting decision-making at the strategic, tactical, and operational levels.
Advances in information systems, together with increasing urbanization, mobility, digitalization, and the rise of customer-centric services, have created a pressing need for anticipatory, real-time decision-making. While classical OR methods address these needs to a significant extent, processing continuous information streams and optimizing operations under highly dynamic and uncertain conditions have become increasingly challenging. 

OR has long studied sequential decisions under uncertainty. Markov decision processes and stochastic dynamic programming provide a foundation for optimizing decisions over time \citep{Puterman1994}. Approximate dynamic programming addresses the computational difficulties of large decision problems through simulation and value function approximation \citep{powell2007approximate,Powell2011}.  Neuro-dynamic programming, which combines dynamic programming with neural network approximations, laid much of the groundwork for modern RL \citep{BertsekasTsitsiklis1996}. Linear programming approaches approximate value functions in large Markov decision processes and can provide bounds on the optimal value or the approximation error, with applications such as revenue management \citep{deFariasVanRoy2003,Adelman2007}. Stochastic programming represents uncertainty and recourse decisions \citep{Sen1999}, with statistical methods available to assess solution quality \citep{Bayraksan2006}. Adaptive multistage sampling estimates optimal values in Markov decision processes and is a precursor of the Monte Carlo tree search methods now widely used in reinforcement learning (RL) \citep{Chang2005}.

Within this broader context, RL provides a framework for learning decision policies through interaction with an environment. Many RL algorithms can be viewed as simulation-based approximate methods for solving Markov decision processes when the transition model is unknown or the state space is too large for exact dynamic programming \citep{BertsekasTsitsiklis1996,Powell2011}. RL algorithms have attracted growing attention in recent years and have demonstrated strong performance in high-dimensional Markov decision processes. Beyond dynamic sequential decision-making, RL has also shown potential to enhance OR methodologies for static combinatorial optimization problems (COPs), where the problem inputs, including any description of uncertainty, are specified before the solution process begins. In this setting, RL either constructs solutions or guides heuristic and exact algorithms, such as branch-and-cut and cutting-plane generation. Digital twins can also support RL by providing models of physical systems for training, testing, and validating policies.

In recent years, several review papers have summarized advances in RL. These include general surveys from a computer-science perspective \citep{Kaelbling1996}, reviews of multi-agent RL \citep{Busoniu2010, Hu2024}, and domain-specific overviews in areas such as continuous control \citep{Recht2019}, building control systems \citep{Wang2020control}, industrial process control \citep{Nian2020}, and robotics \citep{Singh2022}. Other studies examine RL from a safety perspective in applications such as autonomous driving and power systems \citep{Gu2024, Sureview2025}, while several reviews focus on deep RL methods more broadly \citep{Mousavi2016, Arulkumaran2017, Wang2020, Ladosz2022} or in specialized fields such as fluid mechanics \citep{Garnier2021} and healthcare \citep{Yu2023}. Additionally, dedicated reviews on model-free RL \citep{Shakya2023} and model-based RL \citep{Luo2024} further enrich this landscape. \citet{Mazyavkina2021} review RL for combinatorial optimization problems (COPs), while \citet{Yan2022review} and \citet{Farazi2021} provide broad coverage of applications in logistics and supply chain management and transportation, respectively. \citet{Boute2022} provide a roadmap for deep RL in inventory control, and \citet{Hildebrandt2023} discuss opportunities for RL in stochastic dynamic vehicle routing. Despite this growing body of work, these surveys focus primarily on RL algorithms, often within a single application area or algorithmic family, and only partly address how RL is integrated with OR methodologies to improve decision-making and optimization performance. The interaction between RL and OR, especially the use of RL to enhance, support, or extend traditional OR algorithms, remains insufficiently explored.

Interest in integrating machine learning with OR algorithms has grown, as reflected in several reviews \citep{Talbi2016,Song2019,Bengio2021,Karimi2022}. For instance, \citet{Bengio2021} provide a methodological review of machine learning (ML) for combinatorial optimization, focusing on how learning components can be embedded within OR algorithms. Their taxonomy contrasts supervised/unsupervised learning with RL and describes integration patterns ranging from end-to-end learning to algorithm configuration and solver guidance. \citet{Karimi2022} examine the use of ML to design components of metaheuristics, with RL treated as a relatively small subset compared with supervised and unsupervised learning. \citet{Scavuzzo2026} review how machine learning can augment branch-and-bound for mixed-integer linear programming, including RL approaches to branching decisions. The unpublished tutorial by \citet{Wu2025}, retained solely as a background reference, provides an introduction on using deep reinforcement learning (DRL), a subset of RL, to solve COPs, with an emphasis on the foundational concepts and key background knowledge underlying the end-to-end learning paradigm.

Despite these contributions, the literature still lacks a dedicated and systematic overview of the integration of RL and OR, particularly regarding the approaches to applying RL to dynamic sequential decision problems and static COPs, the use of RL to guide heuristic and exact optimization methods, and its integration with digital twins (DTs) for OR applications, together with a clear discussion of the associated research opportunities. Specifically, two major gaps remain. First, existing review papers do not specifically systematize how RL enhances heuristics, exact methods, and commercial solvers nor clarify the mechanisms through which RL improves algorithmic performance. Second, there is no unified taxonomy that explicitly separates RL for dynamic sequential decision-making problems from RL for static COPs, while also distinguishing both methodological roles from the deployment and evaluation layer enabled by digital twins. Such a taxonomy is needed to map RL to concrete integration points in OR workflows, including algorithm selection, initialization, operator learning, parameter tuning, bound tightening, policy validation, and pre-deployment testing.

This paper addresses these gaps by providing a comprehensive review of how RL empowers OR. The work most closely related to ours is given by \citet{Bengio2021}, which surveys ML for combinatorial optimization problems but treats RL as one element within a broader landscape. In contrast, our study provides a focused and in-depth analysis of how RL enhances OR models and workflows beyond classical combinatorial optimization. Before outlining our contributions, we clarify the scope of OR problems considered in this review. In dynamic sequential decision problems, the system state evolves between decisions, new information is revealed over time, and the objective is to find a policy that specifies decisions as the system evolves. In static COPs, the problem inputs, including any description of uncertainty such as scenarios or ambiguity sets, are specified before the solution process begins, and the objective is to construct a feasible solution such as a route, schedule, assignment, or facility set. Although RL-based methods for static COPs may construct a solution step by step or guide an iterative search, these steps are part of the algorithm that solves a given problem instance rather than decisions in an evolving system. Some studies solve a static problem repeatedly within a dynamic setting, for example in a rolling-horizon scheme. We classify such studies by what RL learns. Studies in which RL learns the operational policy are reviewed in Section \ref{Sec_MDP}, and studies in which RL improves the algorithm that solves each static problem are reviewed in Section \ref{Sec_algorithm}. This distinction underlies the two methodological roles in our taxonomy. To ensure analytical depth, this review focuses on transportation and logistics and on closely related operations problems, such as inventory control, production scheduling, and maintenance. The taxonomy is developed and examined mainly in this context. A small number of selected studies from other domains are included to illustrate RL methods and their integration with OR. Some organizing dimensions may also be relevant to healthcare and energy systems, but the limited examples reviewed here do not establish their applicability to these domains.

Our contributions are five-fold. First, we introduce a taxonomy that distinguishes two methodological roles of RL in OR from the use of digital twins for deployment and evaluation. The first role concerns repeated decisions in dynamic environments, where states evolve over time. The second concerns static COPs, where the optimization instance is specified before the solution process begins and RL constructs solutions or guides heuristic and exact algorithms. For this second role, we also summarize the functions and integration points of RL within OR algorithms. Digital twins form an optional layer that can support both roles through policy training, testing, validation, and operational control. Second, we synthesize the main categories of RL algorithms for solving sequential decision-making problems with different structural characteristics. Third, we systematize how RL can guide both heuristic and exact optimization methods for static COPs, and discuss when RL improves solution quality and convergence relative to traditional OR baselines. Fourth, we examine how digital twins support RL training and acceleration, and how RL supports operational control with DTs. We compare these mechanisms and discuss model accuracy, resource requirements, and validation in physical systems. Fifth, we consolidate guidance on when to consider different RL algorithm families or hybrid RL-OR methods, and we highlight key opportunities, challenges, and future research directions in this interdisciplinary area. 

The remainder of the paper is organized as follows. Section \ref{Sec_knowledge} provides background on RL. Section \ref{Sec_reviewMethod} outlines the review methodology. Section \ref{Sec_MDP} focuses on the use of RL for dynamic sequential decision-making problems. Section \ref{Sec_algorithm} reviews state-of-the-art methods combining RL with OR algorithms for static COPs. Section \ref{Sec_digitalTwins} examines digital twins for the deployment and evaluation of RL in OR. Finally, Section \ref{Sec_conclusion} concludes the paper.

\section{RL framework}
\label{Sec_knowledge}

RL is a subfield of machine learning focused on learning decision-making policies through repeated interaction and feedback. It has been widely applied in areas such as automatic control, robotics, transportation logistics, and game playing. For a comprehensive introduction, we refer to \citet{sutton2018}. In this section, we provide a concise overview of key RL concepts (Section \ref{sec_keyconcept}) and the main categories of RL algorithms (Section \ref{Sec_RLalgorithms}).

\subsection{Key concepts of RL}
\label{sec_keyconcept}

RL models decision-making as an interaction between an agent and an environment. The key elements include a \textit{state}, an \textit{agent}, a \textit{policy}, a \textit{reward signal}, a \textit{value function}, and optionally a \textit{model} of the environment. As indicated in Figure \ref{fig_RL}, at each decision-making stage $t$, the agent observes a state $S_t$, selects an action $A_t$ according to a policy $\pi$, receives a reward $R_{t+1}$, and the environment transitions to a new state $S_{t+1}$. The objective is to learn a policy that maximizes long-term reward.

Specifically, the \textit{policy} determines which action the agent takes in each state. The \textit{value function} estimates the expected long-term return from a state or state-action pair under a given policy. A \textit{model}, when present, predicts state transitions and rewards. The agent is not explicitly told which actions to take; rather, it must learn through trial and error which actions yield high long-term returns. Consequently, actions influence not only immediate reward but also the sequence of future states and rewards.

\par\addvspace{\intextsep}
\noindent\begin{minipage}{\linewidth}
    \captionsetup{type=figure,font={small,sf},labelfont={bf,color=scolor},labelsep=colon,skip=6pt}
    \linespread{1}\selectfont
    \sffamily\small
    \centering

\tikzset{every picture/.style={line width=0.75pt}} 

\begin{tikzpicture}[x=0.75pt,y=0.75pt,yscale=-1,xscale=1]

\draw   (268,24) -- (338,24) -- (338,64) -- (268,64) -- cycle ;
\draw   (255,153) -- (363,153) -- (363,207) -- (255,207) -- cycle ;
\draw [line width=1.5]    (338.5,42.5) -- (413.5,42.5) ;
\draw [line width=1.5]    (413.5,42.5) -- (413,183.5) ;
\draw [line width=1.5]    (413,183.5) -- (366.5,183.5) ;
\draw [shift={(362.5,183.5)}, rotate = 360] [fill={rgb, 255:red, 0; green, 0; blue, 0 }  ][line width=0.08]  [draw opacity=0] (11.61,-5.58) -- (0,0) -- (11.61,5.58) -- cycle    ;
\draw [line width=1.5]    (181,182.5) -- (223.8,182.6) ;
\draw [line width=1.5]    (181.5,42.5) -- (181,182.5) ;
\draw [line width=1.5]    (264.5,42.98) -- (181.5,42.5) ;
\draw [shift={(268.5,43)}, rotate = 180.33] [fill={rgb, 255:red, 0; green, 0; blue, 0 }  ][line width=0.08]  [draw opacity=0] (11.61,-5.58) -- (0,0) -- (11.61,5.58) -- cycle    ;
\draw    (255.5,170) -- (227.5,170) ;
\draw [shift={(225.5,170)}, rotate = 360] [color={rgb, 255:red, 0; green, 0; blue, 0 }  ][line width=0.75]    (10.93,-3.29) .. controls (6.95,-1.4) and (3.31,-0.3) .. (0,0) .. controls (3.31,0.3) and (6.95,1.4) .. (10.93,3.29)   ;
\draw    (199.8,169.8) -- (225.5,170) ;
\draw    (199.8,169.8) -- (200.2,55.8) ;
\draw    (200.2,55.8) -- (265.4,56.58) ;
\draw [shift={(267.4,56.6)}, rotate = 180.68] [color={rgb, 255:red, 0; green, 0; blue, 0 }  ][line width=0.75]    (10.93,-3.29) .. controls (6.95,-1.4) and (3.31,-0.3) .. (0,0) .. controls (3.31,0.3) and (6.95,1.4) .. (10.93,3.29)   ;
\draw [line width=1.5]    (227.8,182.6) -- (255,182.6) ;
\draw [shift={(223.8,182.6)}, rotate = 0] [fill={rgb, 255:red, 0; green, 0; blue, 0 }  ][line width=0.08]  [draw opacity=0] (11.61,-5.58) -- (0,0) -- (11.61,5.58) -- cycle    ;
\draw [color={rgb, 255:red, 128; green, 128; blue, 128 }  ,draw opacity=1 ] [dash pattern={on 4.5pt off 4.5pt}]  (225,148.2) -- (224.6,203) ;

\draw (282,37) node [anchor=north west][inner sep=0.75pt]   [align=left] {{\small Agent}};
\draw (267.4,172.2) node [anchor=north west][inner sep=0.75pt]   [align=left] {{\small Environment}};
\draw (419,106.5) node [anchor=north west][inner sep=0.75pt]   [align=left] {action \textit{A{\footnotesize t}}};
\draw (125,109) node [anchor=north west][inner sep=0.75pt]   [align=left] {state \textit{S{\footnotesize t}}};
\draw (203.4,63.8) node [anchor=north west][inner sep=0.75pt]   [align=left] {reward \textit{R{\footnotesize t}}};
\draw (225,151.4) node [anchor=north west][inner sep=0.75pt]   [align=left] {\textit{{\small R}{\scriptsize t+}}{\scriptsize 1}};
\draw (225,187.6) node [anchor=north west][inner sep=0.75pt]   [align=left] {\textit{{\small S}{\scriptsize t+}}{\scriptsize 1}};

\end{tikzpicture}
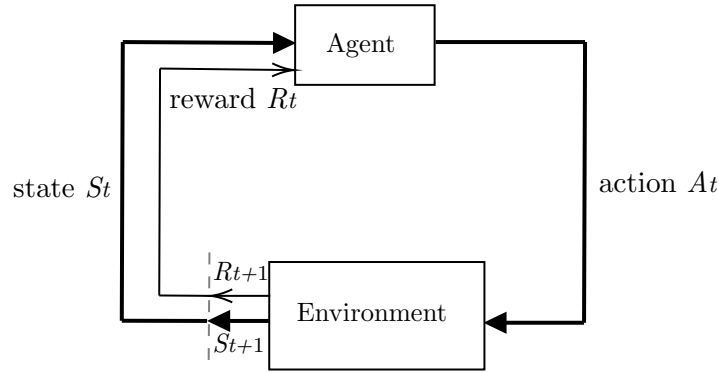
    \caption{Illustration of the general RL setting.}
    \label{fig_RL}
\end{minipage}
\par\addvspace{\intextsep}

RL offers a flexible framework for learning policies that optimize long-term performance in stochastic and dynamic environments without requiring an explicit analytical model, and it can be combined with optimization or control methods. Two general challenges are the design of informative reward functions, especially when rewards are sparse, and the balance between exploring new actions and exploiting actions already known to be effective.

\subsection{RL algorithms}
\label{Sec_RLalgorithms}

RL algorithms can be broadly categorized into \textit{model-free} and \textit{model-based} approaches, depending on whether the agent has access to (or learns) an explicit model of the environment.

\subsubsection{Model-free RL algorithms}
Model-free algorithms learn the value of actions directly through extensive trial-and-error interactions with the environment. Depending on their learning objective, they can be grouped into three categories:

\textbf{(i) Value-based RL algorithms.} These algorithms estimate the value of state-action pairs and select actions by choosing those with the highest estimated values. The policy is derived indirectly by acting greedily with respect to the learned value function. Representative algorithms include Q-learning \citep{Watkins1992}, Deep Q-Network (DQN) \citep{Mnih2015}, Double DQN \citep{Hasselt2016}, and Dueling DQN \citep{Wang2016}.

\textbf{(ii) Policy-Based RL algorithms.} These algorithms learn a policy function directly by estimating the probability of taking an action in a given state. Instead of evaluating individual actions, they focus on maximizing the expected return of the entire policy. This is typically achieved through gradient-based updates, increasing the likelihood of effective actions. Examples include Trust Region Policy Optimization (TRPO) \citep{Schulman2017TRPO} and Proximal Policy Optimization (PPO), used, for example, for inventory control by \citet{Dehaybe2024}.

\textbf{(iii) Actor-Critic RL algorithms.} These algorithms combine the advantages of value-based and policy-based approaches. The actor proposes actions according to a policy network, while the critic evaluates them using a value function. Both components are trained simultaneously to improve decision-making efficiency and stability. Representative algorithms include Advantage Actor-Critic (A2C) \citep{Vijay1999}, Deep Deterministic Policy Gradient (DDPG) \citep{Lillicrap2015}, Asynchronous Advantage Actor-Critic (A3C) \citep{Mnih2016}, Soft Actor-Critic (SAC) \citep{Haarnoja2018}, and Twin Delayed DDPG (TD3) \citep{Fujimoto2018}.

\subsubsection{Model-based RL algorithms}

Model-based RL uses an explicit model of state transitions and rewards to evaluate candidate actions, generate imagined experience, or improve a policy. The model may be specified from prior knowledge or learned from data. The defining feature is its use in planning or policy or value updates. Training an otherwise model-free agent in a simulator is not sufficient for this classification. With a known model, the main computational task is planning or policy optimization. With a learned model, model estimation and policy improvement are coupled, and prediction errors must be controlled.

Representative algorithms include Dyna-Q \citep{Sutton1991}, Stochastic Ensemble Value Expansion (STEVE) \citep{Buckman2019}, Model-Based Policy Optimization (MBPO) \citep{Janner2019}, and Robust Adversarial Model-Based Offline Reinforcement Learning (RAMBO-RL) \citep{Rigter2022}. By simulating interactions within the learned model, these algorithms reduce reliance on costly real-world exploration and often accelerate learning.

\section{Taxonomy and review methodology}
\label{Sec_reviewMethod}

In this section, we first present the taxonomy of our review and then outline the search strategies and inclusion criteria.

\subsection{Taxonomy}
\label{sec_Objective}

The purpose of this review is to examine how RL contributes to OR in transportation and logistics and closely related operations problems. As illustrated in Figure \ref{fig_Taxonomy}, we distinguish two methodological roles of RL in OR. These concern dynamic sequential decision problems and solution methods for static COPs. Digital twins provide an optional deployment and evaluation layer that can support both roles.

\textbf{(i) RL for solving dynamic sequential decision problems.} These problems involve dynamic systems in which states evolve between decisions and critical information (e.g., demand, disturbances, opponent behavior) is revealed over time. The objective is to find a policy for repeated decisions that optimizes cumulative system performance, measured by criteria such as discounted reward, finite-horizon return, or long-run average reward. Such problems are typically formulated as Markov Decision Processes and solved using RL algorithms.

\textbf{(ii) RL as an optimization tool for static combinatorial optimization problems.} In static COPs, the problem inputs, including any description of uncertainty such as scenarios or ambiguity sets, are specified before the solution process begins, and the objective is to construct a discrete solution (e.g., a route, schedule, assignment, or facility set) satisfying the relevant constraints. Time is not dynamic in this setting: even when RL constructs a solution step by step or guides an iterative search process, the sequence is part of the algorithm that solves a given problem instance rather than the evolution of an operational system. RL can be employed as an end-to-end solver, integrated with heuristic methods, or embedded into exact algorithms to enhance the performance of OR solution methods.

\textbf{Digital twins for RL deployment and evaluation.} DTs combine digital models with operating data to support policy training, testing, validation, and operational control. Their role concerns how RL methods are developed and used in relation to a physical system. We therefore distinguish this optional layer from the two methodological roles above.

\par\addvspace{\intextsep}
\noindent\begin{minipage}{\linewidth}
    \captionsetup{type=figure,font={small,sf},labelfont={bf,color=scolor},labelsep=colon,skip=6pt}
    \linespread{1}\selectfont
    \sffamily\small
    \centering
    {\color{black}\sffamily\linespread{1.0}\selectfont\setlength{\parindent}{0pt}%
    \newcommand{\rlband}[2]{\begin{minipage}[t][2.85cm][t]{4.5cm}\raggedright
        \hyphenpenalty=10000\exhyphenpenalty=10000
        {\fontsize{9}{11}\selectfont\textbf{#1}}\par\vspace{4pt}#2\end{minipage}}%
    \newcommand{\rlitem}[1]{\hangindent=0.85em\hangafter=1\makebox[0.85em][l]{\textbullet}#1\par\vspace{1.5pt}}%
    \newcommand{\rltag}[1]{\setlength{\fboxsep}{2pt}\fcolorbox{black}{white}{\fontsize{7.5}{9}\selectfont\textbf{#1}}}%
    \begin{tikzpicture}[x=1cm,y=1cm,font=\fontsize{9}{11}\selectfont,
        band/.style={draw=black!40,fill=white,line width=0.5pt,rounded corners=3pt,
                     align=left,inner sep=5pt,font=\fontsize{8.5}{10}\selectfont},
        role/.style={draw=black,fill=white,line width=0.9pt,rounded corners=3pt,
                     align=center,text width=4.5cm,minimum height=2.42cm,inner sep=6pt},
        method/.style={draw=black,fill=black!7,line width=0.6pt,rounded corners=3pt,
                       align=center,text width=3.3cm,minimum height=1.6cm,inner sep=4pt}]

        \node[band,anchor=south] (dynband) at (-5.5,4.75)
            {\rlband{Dynamic sequential decision problems}
             {\rlitem{States evolve between decisions}
              \rlitem{Information is revealed over time}
              \rlitem{RL learns a policy for dynamic decisions}}};
        \node[band,anchor=south] (statband) at (-0.45,4.75)
            {\rlband{Static combinatorial optimization problems}
             {\rlitem{\textcolor{black}{Inputs are specified before solving begins}}
              \rlitem{Solution steps are algorithmic, not system evolution}
              \rlitem{RL solves the COP alone or integrates with OR algorithms}}};

        \node[role] (dynamic) at (-5.5,3.34)
            {\textbf{RL for solving dynamic\\sequential decision problems}\\[7pt]
             RL acts as the policy that\\maps observed states\\to decisions};
        \node[role] (static) at (-0.45,3.34)
            {\textbf{RL for solving\\static COPs}\\[7pt]
             RL builds solutions directly\\or guides the steps of\\OR algorithms};

        \node[align=center,text width=3.6cm,font=\fontsize{8.5}{10}\selectfont]
            at (4.35,8.35) {\textbf{RL for solving static COPs}};
        \node[method] (solver) at (4.35,7.15)
            {Independent RL solver\\[3pt]\rltag{RL REPLACES}};
        \node[method] (heuristic) at (4.35,5.07)
            {RL with heuristic\\algorithms\\[3pt]\rltag{RL INSIDE}};
        \node[method] (exact) at (4.35,2.99)
            {RL with exact\\algorithms\\[3pt]\rltag{RL INSIDE}};

        \draw[line width=0.6pt] (static.east) -- (2.3,3.34);
        \draw[line width=0.6pt] (2.3,2.99) -- (2.3,7.15);
        \draw[->,line width=0.6pt] (2.3,7.15) -- (solver.west);
        \draw[->,line width=0.6pt] (2.3,5.07) -- (heuristic.west);
        \draw[->,line width=0.6pt] (2.3,2.99) -- (exact.west);

        \node[draw=black,dashed,line width=0.7pt,rounded corners=3pt,
              align=center,text width=13.7cm,minimum height=1.5cm,inner sep=6pt]
              (dt) at (-0.91,1.0)
            {\textbf{Digital twins for RL deployment and evaluation}\\[4pt]
             Optional layer supporting both methodological roles\\[3pt]
             Training, testing, policy validation, and operational control};

        \draw[->,line width=0.7pt] (dt.north -| dynamic.south) -- (dynamic.south);
        \draw[->,line width=0.7pt] (dt.north -| static.south) -- (static.south);

        \node[anchor=north,align=center,text width=13.7cm,
              font=\fontsize{8}{10}\selectfont] at (-0.91,0.0)
            {RL = Reinforcement Learning. OR = Operational Research.\\
             COPs = Combinatorial Optimization Problems.};
    \end{tikzpicture}
    }
    \caption{Taxonomy of RL in OR.}
    \label{fig_Taxonomy}
\end{minipage}
\par\addvspace{\intextsep}

We include studies that satisfy the following conditions:
(a) the study addresses an optimization or control problem relevant to OR, and
(b) the study explicitly involves an RL component. For DT-related work, we include studies where the digital twin serves as an operational decision-support system, enabling optimization, control, or performance evaluation typical of OR applications. A formal OR model is not required, but the study must define operational objectives and constraints and evaluate them through the DT.

Application coverage focuses on transportation and logistics and closely related operations problems. A small number of selected studies from other domains are included to illustrate RL methods and DT integration, and they are marked in the classification tables. These examples do not provide a comprehensive review of those domains.

\subsection{Review methodology}
\label{sec_Search}

To identify application studies for the main review of RL and OR, we conducted literature searches across major scientific databases, including Scopus, Web of Science Core Collection, IEEE Xplore, ScienceDirect, SpringerLink, and INFORMS PubsOnline. In addition, we performed a venue-restricted search of leading OR journals (e.g., INFORMS journals, European Journal of Operational Research (EJOR), and Transportation Research Parts A–F) as well as top-tier AI conferences (NeurIPS, ICML, ICLR, AAAI) to capture online first and early access papers not yet indexed in databases.

The keyword strategy below applies to the main review of studies that use RL for OR problems or combine RL with OR algorithms. Foundational OR literature and public software libraries and benchmark resources were identified through separate targeted searches using method names, publication titles, and resource names. Publication details were checked against publisher records, and resource information was checked against official documentation and repositories. These supporting references provide methodological context or evaluation guidance and are not included in the application study counts or classification tables. Requiring RL terminology can miss related work described only in terms of dynamic programming or approximate dynamic programming. The separate searches provide selected background references to explain these connections, rather than systematic coverage of the broader OR literature on sequential decisions under uncertainty.

We adopted a layered keyword search strategy that combines three mandatory blocks and an optional refinement block:
\[
\textbf{Query} = (\mathcal{A}\ \text{RL})~\wedge~(\mathcal{B}\ \text{applications})~\wedge~(\mathcal{C}\ \text{OR methods}).
\]

\noindent\textbf{$\mathcal{A}$ (RL) includes RL algorithms.} Terms are grouped by families as follows:
\begin{description}[leftmargin=1.2em, itemsep=0.25ex, topsep=0.4ex, font=\itshape]
  \item[$\bullet$ Generic \& model-based RL family.] 
  ``reinforcement learning'' $\vee$ ``deep reinforcement learning'' $\vee$ 
  ``model-free reinforcement learning'' $\vee$ ``model-based reinforcement learning'' $\vee$
  ``model based reinforcement learning'' $\vee$ ``MBRL'' $\vee$
  ``model-based value expansion''.
  
  \item[$\bullet$ Value-based RL family.] 
  ``value-based reinforcement learning'' $\vee$ ``Q-learning'' $\vee$
  ``deep Q-Network'' $\vee$ ``double deep Q-network'' $\vee$ ``dueling deep Q-network''.
  
  \item[$\bullet$ Policy-based RL family.]  
  ``policy-based reinforcement learning'' $\vee$ 
  ``trust region policy optimization'' $\vee$ ``proximal policy optimization'' $\vee$
  ``deep deterministic policy gradient'' $\vee$ 
  ``twin delayed deep deterministic policy gradient''.
  
  \item[$\bullet$ Actor-critic RL variants.] 
  ``actor-critic'' $\vee$ ``advantage actor-critic'' $\vee$
  ``asynchronous advantage actor-critic'' $\vee$ ``soft actor-critic''.
\end{description}

\noindent\textbf{$\mathcal{B}$ (applications) consists of research domains.} Domains are grouped by themes as follows:
\begin{description}[leftmargin=1.2em, itemsep=0.25ex, topsep=0.4ex, font=\itshape]
  \item[$\bullet$ Transport \& logistics.] 
  ``transport'' $\vee$ ``logistics''$\vee$ ``vehicle routing'' $\vee$
  ``scheduling'' $\vee$ ``timetabling'' $\vee$ ``line planning'' $\vee$ ``network design''.
  \item[$\bullet$ Operations management.] 
  ``inventory control'' $\vee$ ``assortment'' $\vee$ ``pricing'' $\vee$
  ``revenue management'' $\vee$ ``facility location'' $\vee$ ``resource allocation''.
\end{description}

The search focuses on transportation and logistics and related operational decisions. Healthcare and energy systems are not treated as separate search domains.

\noindent\textbf{$\mathcal{C}$ (OR methods).} OR methods are grouped by categories as follows:
\begin{description}[leftmargin=1.2em, itemsep=0.25ex, topsep=0.4ex, font=\itshape]
  \item[$\bullet$ Mathematical programming \& decomposition.] 
  ``operational research'' $\vee$ ``mixed-integer programming'' $\vee$
  ``integer programming'' $\vee$ ``linear programming'' $\vee$
  ``combinatorial optimization'' $\vee$ ``column generation'' $\vee$
  ``Benders decomposition'' $\vee$ ``Lagrangian relaxation'' $\vee$ ``decomposition''.
  \item[$\bullet$ Branching \& cutting families. ] 
  ``branch-and-price'' $\vee$ ``branch-and-cut'' $\vee$ ``cutting plane'' $\vee$ ``branch and bound''.
  \item[$\bullet$ Constraint-based family.] 
  ``constrained (or relaxed) decision diagrams'' $\vee$ ``constraint programming''.
  \item[$\bullet$ Solvers and generic.] 
  ``CPLEX'' $\vee$ ``GUROBI'' $\vee$ ``Markov decision process'' $\vee$ ``MDP'' $\vee$
  ``digital twin'' $\vee$ ``digital replica''.
\end{description} 

To improve coverage of model-based RL within transportation and logistics, we supplemented the search using the three keyword blocks with targeted combinations of (``model-based reinforcement learning'' OR ``model based reinforcement learning'' OR MBRL OR ``world model'' OR ``Dyna-Q'') and (transportation OR logistics OR routing OR delivery OR scheduling OR inventory OR warehouse OR freight OR traffic OR rail OR port OR ``ride-hailing''). Block $\mathcal{C}$ was not mandatory in this supplement, because relevant applications may not use OR terminology in their titles, abstracts, or keywords. Publisher and conference records and author or institutional repositories were checked to verify publication status and whether an explicit model supports planning or policy or value updates. The same publication status criteria apply to the supplementary search.
 
We followed PRISMA-style stages for systematic review:
Identification, Screening, Eligibility, and Inclusion. The inclusion criteria are as follows: (i) the paper includes an explicit RL component, (ii) there is a clear OR connection through models, algorithms, or applications, (iii) the application study was published between 2016 and 2026, and (iv) the study has been formally published in a journal or in computing conference proceedings, including articles published online first or as early access. We cite the published version when a preprint and a published article describe the same study. Manuscripts available only on arXiv, other unpublished working papers, and workshop contributions without formal proceedings are excluded from the reviewed evidence and classification tables. Foundational background references may predate the period covered by the application studies. The sole unpublished manuscript retained for background is the closely related tutorial by \citet{Wu2025}, which is cited in the Introduction to position our contribution and is not counted among the reviewed application studies. Public software repositories and benchmark websites are also cited as resources for the evaluation guidance in Section~\ref{sec_benchmarking}. These resources are not counted as reviewed application studies.

\section{RL: A dynamic sequential decision-making tool}
\label{Sec_MDP}

In this section, we review studies that apply RL algorithms to solve sequential decision-making problems in dynamic environments. In these problems, the system state evolves over time, future information is not fully known before decisions are made, and decisions must be optimized sequentially to improve long-term performance. These problems are typically formulated as Markov Decision Processes \citep[MDPs, e.g.,][]{powell2007approximate}, where the aim is to learn a policy that selects actions over time to maximize cumulative rewards based on state observations and feedback from the environment.

RL can be compared with established OR methods for dynamic sequential decision-making problems in terms of the model information they require, the approximations they use, and how their solutions are evaluated. Standard dynamic programming uses specified transition and reward models \citep{Puterman1994}. Model-free RL learns values or policies from sampled experience, while model-based RL uses known or learned models. Approximate dynamic programming (ADP) also uses simulation, value function approximation, and policy search \citep{powell2007approximate,Powell2011,Bertsekas2012}. These ideas overlap with those used in the value-based, policy-based, and actor-critic methods reviewed below, and the two literatures often describe them under different names. This overlap is visible in the reviewed studies. \citet{Dehaybe2024} compare deep reinforcement learning with dynamic programming-based heuristics, \citet{Yan2023} compare their RL approach with an ADP baseline, and \citet{Verleijsdonk2024} use approximate policy iteration. The comparison therefore concerns the information and approximations used, as well as the costs of training and computing decisions. Multistage stochastic programming represents uncertainty through scenarios, models recourse decisions explicitly, and includes operational constraints directly in the optimization model \citep{Sen1999}.  For multistage stochastic linear programs, stochastic dual dynamic programming (SDDP) approximates expected future costs with piecewise linear functions constructed from Benders cuts \citep{Pereira1991}. SDDP and value-based RL both approximate future performance, but SDDP constructs its approximations from the optimization model and its dual solutions, while RL can fit value functions or policies from sampled experience. Under its standard assumptions, SDDP constructs a lower bound on the optimal cost that supports convergence assessment, whereas the RL methods reviewed here generally do not provide a comparable bound. Standard SDDP relies on convex recourse and commonly assumes stagewise-independent uncertainty, although extensions handle some forms of dependence. RL methods can also be applied to problems with discrete decisions, nonlinear dynamics, or simulator-based models. Comparisons should therefore use consistent objectives, constraints, and uncertainty assumptions.

Information relaxation bounds the optimal value by relaxing restrictions on access to future information and imposing suitable penalties \citep{Brown2010}. These bounds could help assess the quality of feasible RL policies under the same problem model, complementing policy learning with an evaluation of solution quality. Robust MDPs optimize policies against the least favorable transition model in a specified set \citep{Iyengar2005}. Distributionally robust MDPs instead use the least favorable admissible distribution of uncertain model parameters \citep{XuMannor2012}. These criteria address model uncertainty and differ from evaluation under a single assumed model. Robust optimization also appears in the reviewed studies. \citet{Li2025} combine robust optimization with adversarial RL for bus service scheduling, and \citet{Li2024COR} and \citet{An2026} use RL-guided heuristics to solve distributionally robust optimization models. Comparisons with RL should separate the uncertainty criterion from the choice of learning or optimization method. For comprehensive treatments of these methods, we refer to \citet{BirgeLouveaux2011} and \citet{Shapiro2021} for stochastic programming and to \citet{Powell2019,Powell2022} for a unified framework of sequential decisions under uncertainty that encompasses ADP and RL.

We classify the reviewed studies based on the RL algorithms they adopt and annotate each work using a set of \textit{orthogonal labels} that capture key modeling and implementation aspects. We first analyze the literature across different RL algorithm families and discuss the problem characteristics suited to each, and then outline key challenges and future research directions. Specifically, Section~\ref{sec_value} reviews studies that solve MDPs using value-based methods. Section~\ref{sec_policy} focuses on policy-based RL algorithms. Section~\ref{sec_actor} surveys actor-critic variants. Section~\ref{sec_model} discusses model-based RL approaches. Section~\ref{sec:open} outlines cross-cutting challenges and future research directions. Given the uneven volume of literature, our coverage follows a consistent convention: for Sections~\ref{sec_value} and ~\ref{sec_actor}, where the number of studies is large, we synthesize cross-study patterns and cite representative works; for Section~\ref{sec_policy}, we briefly discuss each study, while Section~\ref{sec_model} groups the model-based literature by application and explains how each method uses its model.

\subsection{Solving MDPs with value-based RL algorithms}
\label{sec_value}

Given the large number of studies in the category of applying value-based RL algorithms to MDPs, we synthesize cross-study patterns rather than reviewing related studies individually. Table \ref{tab_valueBaedMDP} summarizes state-of-the-art studies, reporting domain/task, problem signature, the employed RL algorithm, training mode, and baselines. Representative applications include railway scheduling~\citep[][]{Semrov2016}, crowdsourced urban delivery~\citep[e.g.,][]{Ahamed2021}, same-day delivery with vehicles and drones~\citep[e.g.,][]{Chen2022, Chen2023}, inventory control~\citep[e.g.,][]{Oroojlooyjadid2022, Mao2025}, reusable resource allocation and pricing~{\hypersetup{citecolor=blue}\citep[][]{WangSongLiuZhao2026}}, adaptive bus dispatching under adversarial demand perturbations~\citep[][]{Li2025}, and strategic transportation planning~\citep[][]{Yu2025}. Most of these studies feature discrete actions, high observability, and mostly episodic horizons.

Overall, value-based methods are most frequently used in resource allocation and scheduling, urban delivery and crowdsourcing, inventory and credit control, and digital marketing-settings with predominantly discrete operational decisions (e.g., replenishment, dispatch/retention, promotion selection).  Most studies assume full observability, with only a few relying on partial or binary observations. Time horizons are typically finite or episodic (daily or period-based), with occasional infinite-horizon discounted models in inventory control. The underlying dynamics commonly incorporate randomness or non-stationarity, such as demand variation, stochastic travel times, perishability, user heterogeneity, and shifting reward structures. These are addressed through period-by-period learning or restart strategies. Training is conducted both online and offline, the latter often using historical logs or simulation, consistent with the off-policy nature of value-based methods. Baseline comparisons range from rule-based and heuristic policies to behavioral policies and human benchmarks, with fewer comparisons against simulation-only strategies. 

In conclusion, sequential decision-making problems with discrete actions, high observability, and episodic or discrete-time formulations are well suited to value-based RL algorithms. Transportation planning and inventory control can also benefit from offline learning pipelines, a setting where value-based methods are particularly effective.

\begin{table}[h]
\centering
\caption{Classification of papers studying value-based RL algorithms for sequential decision-making.}
\label{tab_valueBaedMDP}
\resizebox{\textwidth}{!}{%
\begin{tabular}{@{} r| l l l l l @{}}
\toprule
Publications & Domain / Task & Problem signature & RL method & Training & Baselines \\
\midrule
\citet{Semrov2016} &
\begin{tabular}[c]{@{}l@{}} Railway / \\ Train rescheduling\end{tabular} &
\begin{tabular}[c]{@{}l@{}} A: discrete; O: full; Horizon: \\ episodic; Dyn: disturbances \end{tabular}&
Q-learning &
Online &
\begin{tabular}[c]{@{}l@{}} FIFO;\\ Random walk \end{tabular}\\
\citet{Ahamed2021} &
\begin{tabular}[c]{@{}l@{}} Crowdsourced \\ urban delivery \end{tabular} &
\begin{tabular}[c]{@{}l@{}}  A: discrete; O: full; \\ Horizon: episodic; \textcolor{black}{Dyn: deterministic}\end{tabular} &
DQN &
Offline &
Heuristics\\
\citet{Chen2022} &
\begin{tabular}[c]{@{}l@{}} Same-day delivery \\ with vehicles and drones\end{tabular} &
\begin{tabular}[c]{@{}l@{}}  A: discrete; O: full; \\ Horizon: daily; Dyn: stochastic \end{tabular}  &
DQN &
Offline&
PFA   \\
\citet{Oroojlooyjadid2022} &
\begin{tabular}[c]{@{}l@{}} Beer game / \\ Multi-echelon inventory \end{tabular} &
\begin{tabular}[c]{@{}l@{}}A: discrete; O: partial/local;\\ Horizon: finite; Dyn: lead-time lags\end{tabular} &
DQN&
Online &
\begin{tabular}[c]{@{}l@{}} Base-stock rules; \\ Behavioral policies \end{tabular} \\
\citet{Guo2022} &
\begin{tabular}[c]{@{}l@{}} Synchromodal logistics / \\ Shipment matching \end{tabular} &
\begin{tabular}[c]{@{}l@{}} A: discrete; O: full; \\ Horizon: finite episodic; \\ Dyn: dynamic and stochastic \\ travel times \end{tabular} &
Q-learning &
Offline &
 MA \\
\citet{Wang2023} &
\begin{tabular}[c]{@{}l@{}}Digital marketing /\\ Sequential personalized \\ promotions\textsuperscript{*} \end{tabular}&
\begin{tabular}[c]{@{}l@{}} A: discrete; O: full;  \\ Horizon: 30-day episodic; \\ Dyn: stochastic users\end{tabular} &
DDDQN &
Online &
\begin{tabular}[c]{@{}l@{}} Mass policies; \\ Myopic; Heuristics; \\ DQN variants \end{tabular}\\
\citet{Chen2023} & Same-day delivery &
\begin{tabular}[c]{@{}l@{}} A: discrete; O: full; \\ Horizon: daily; Dyn: stochastic\end{tabular}&
DQN &
Offline  &
\begin{tabular}[c]{@{}l@{}} Bucket policy; \\ Reserved-vehicle policy\end{tabular} \\
\citet{Alfonso2024} &
\begin{tabular}[c]{@{}l@{}}Credit limit \\ adjustment\textsuperscript{*}\end{tabular} &
\begin{tabular}[c]{@{}l@{}} A: binary; O: full; \\ Horizon: episodic over customers; \\ Dyn: stochastic\end{tabular} &
\begin{tabular}[c]{@{}l@{}} Double \\ Q-learning\end{tabular} &
Offline  &
\begin{tabular}[c]{@{}l@{}} Alternative \\ strategies\end{tabular} \\
{\hypersetup{citecolor=blue}\citet{WangSongLiuZhao2026}} &
\begin{tabular}[c]{@{}l@{}} Reusable resource \\ allocation \& pricing \end{tabular} &
\begin{tabular}[c]{@{}l@{}}  A: discrete; O: full; \\ Horizon: finite; \\ \textcolor{black}{Dyn: stochastic rental times} \end{tabular} &
 \begin{tabular}[c]{@{}l@{}} Value-based \\ episodic RL \end{tabular} &
Online &
\begin{tabular}[c]{@{}l@{}}  Full-information greedy \\ oracle; $\varepsilon$-greedy variants \end{tabular}\\
\citet{Yu2025} &
\begin{tabular}[c]{@{}l@{}} Transportation / \\ Strategic planning \end{tabular} &
\begin{tabular}[c]{@{}l@{}}  A: continuous; \\ O: partial; \\ \textcolor{black}{Horizon: finite; Dyn: deterministic} \end{tabular}&
DQN  &
Simulator&
\begin{tabular}[c]{@{}l@{}} Human reactive; \\ Random; \\  Value-only policy \end{tabular} \\
\textcolor{black}{\citet{Li2025}} &
\textcolor{black}{\begin{tabular}[c]{@{}l@{}} Bus services / \\ Adaptive dispatching \end{tabular}} &
\textcolor{black}{\begin{tabular}[c]{@{}l@{}} A: discrete headways; O: full; \\ Horizon: episodic; Attacker--defender; \\ Dyn: adversarial demand perturbations \end{tabular}} &
\textcolor{black}{\begin{tabular}[c]{@{}l@{}} Adversarial \\ DDDQN \end{tabular}} &
\textcolor{black}{Simulator} &
\textcolor{black}{\begin{tabular}[c]{@{}l@{}} Defender-only \\ DDDQN variants \end{tabular}} \\
\citet{Mao2025} &
\begin{tabular}[c]{@{}l@{}} Nonstationary MDPs / \\ Inventory control \end{tabular}  &
\begin{tabular}[c]{@{}l@{}} A: discrete; O: full; \\ Horizon: episodic; \\ \textcolor{black}{Dyn: non-stationary} \end{tabular}  &
RestartQ-UCB &
Simulator & \begin{tabular}[c]{@{}l@{}} Q-Learning UCB; \\ $\varepsilon$-greedy \end{tabular}\\
\bottomrule
\end{tabular}%
}
\\[2pt]
\begin{minipage}{\textwidth}
\footnotesize\raggedright
\textit{Notes:} [a] A = action type; O = observation; Dyn = dynamics. FIFO = First-In-First-Out, DQN = Deep Q-Network, PFA = Policy
Function Approximation, MA = Myopic approach, DDDQN = Double Dueling DQN, RestartQ-UCB = Restarted Q-learning with Upper Confidence Bounds.
\par
[b] * marks problems or application domains outside transportation, logistics, and closely related operations problems; the corresponding studies are included for methodological illustration.
\end{minipage}
\end{table}

\subsection{Solving MDPs with policy-based RL algorithms}
\label{sec_policy}

Table \ref{tab_policeBaedMDP} summarizes advancements in solving MDPs using policy-based RL algorithms. These methods have been applied across diverse domains including inventory lot-sizing \citep{Dehaybe2024}, container port truck dispatching \citep{Jin2024}, urban rail transit rescheduling \citep{Ying2024}, equipment maintenance \citep{Verleijsdonk2024}, lateral transshipment in health supply chains \citep{Vanvuchelen2024}, rack storage assignment and replenishment in robotic warehouses \citep{Teck2025}, and dynamic job-shop scheduling \citep{Yu2026}. As can be observed from Table \ref{tab_policeBaedMDP}, most applications involve real-time decision-making problems. Action types are predominantly discrete, with a few examples involving continuous actions such as inventory control. Observability is generally full. Most problems feature finite or episodic horizons, while some inventory and maintenance settings operate under infinite-horizon discounting. All studied environments involve randomness or non-stationarity arising from factors such as demand variations, order or job arrivals, schedule disruptions, or system degradation. Training strategies vary accordingly. Online learning is used for scheduling and inventory applications, simulators are used for transshipment, warehouse, and job-shop applications, and simulation-based rollouts are employed in maintenance problems. The baseline algorithms also vary widely, ranging from dynamic programming-based heuristics and handcrafted dispatching, storage, and transshipment rules to DRL-based hyper-heuristics, metaheuristics, and distributed PPO.

\begin{table}[h]
\centering
\caption{Classification of papers studying policy-based RL algorithms for sequential decision-making.}
\label{tab_policeBaedMDP}
\resizebox{\textwidth}{!}{%
\begin{tabular}{@{} r| l l l l l @{}}
\toprule
Publications & Domain / Task & Problem signature & RL method & Training & Baselines \\
\midrule
\citet{Dehaybe2024} &
Inventory lot-sizing  &
\begin{tabular}[c]{@{}l@{}}A: continuous; O: full; Horizon:\\ rolling discounted; Dyn: non-stationary\end{tabular} &
PPO &
Online &
\begin{tabular}[c]{@{}l@{}} DP-based heuristics \end{tabular} \\
\citet{Jin2024} & \begin{tabular}[c]{@{}l@{}} Container port \\ truck dispatching \end{tabular} & 
\begin{tabular}[c]{@{}l@{}} A: discrete; O: full; \\ Horizon: finite;\ Dyn: stochastic \end{tabular} & 
PPO & Online & \begin{tabular}[c]{@{}l@{}}  Handcrafted rule; \\ DRL-HH \end{tabular} \\
\citet{Ying2024} & \begin{tabular}[c]{@{}l@{}}Urban rail / \\ Disruption \\ rescheduling \end{tabular} & 
\begin{tabular}[c]{@{}l@{}} A: discrete; O: full; \\ Horizon: finite; Dyn: non-stationary \end{tabular} & 
\begin{tabular}[c]{@{}l@{}} Multi-agent \\ PPO \end{tabular} & Online & \begin{tabular}[c]{@{}l@{}} Metaheuristics; \\ distributed PPO \end{tabular}\\
\citet{Verleijsdonk2024} &
DTMPA &
\begin{tabular}[c]{@{}l@{}}  A: discrete per engineer; O: full; Horizon: \\  infinite; Dyn: stochastic degradation\end{tabular} &
\begin{tabular}[c]{@{}l@{}} API with \\ DCL \end{tabular} &
\begin{tabular}[c]{@{}l@{}} In-simulation \\ rollouts \end{tabular}&
\begin{tabular}[c]{@{}l@{}} Ranking \\ heuristics\end{tabular} \\
\textcolor{black}{\citet{Vanvuchelen2024}} &
\textcolor{black}{\begin{tabular}[c]{@{}l@{}} Health supply chain / \\ Lateral transshipment \end{tabular}} &
\textcolor{black}{\begin{tabular}[c]{@{}l@{}} A: continuous, mapped to feasible \\ quantities; O: full; Horizon: infinite; \\ Dyn: seasonal demand and lead times \end{tabular}} &
\textcolor{black}{PPO} &
\textcolor{black}{Simulator} &
\textcolor{black}{\begin{tabular}[c]{@{}l@{}} Transshipment \\ heuristics \end{tabular}} \\
\textcolor{black}{\citet{Teck2025}} &
\textcolor{black}{\begin{tabular}[c]{@{}l@{}} Robotic mobile fulfillment / \\ Rack storage assignment \\ and replenishment \end{tabular}} &
\textcolor{black}{\begin{tabular}[c]{@{}l@{}} A: discrete with masking; O: full; \\ Horizon: episodic; \\ Dyn: stochastic order requests \end{tabular}} &
\textcolor{black}{PPO} &
\textcolor{black}{Simulator} &
\textcolor{black}{\begin{tabular}[c]{@{}l@{}} Storage and replenishment \\ decision rules; Full-information \\ MILP (small instances) \end{tabular}} \\
\textcolor{black}{\citet{Yu2026}} &
\textcolor{black}{\begin{tabular}[c]{@{}l@{}} Manufacturing / \\ Dynamic job-shop \\ scheduling \end{tabular}} &
\textcolor{black}{\begin{tabular}[c]{@{}l@{}} A: discrete (PDRs); O: full; \\ Horizon: episodic; Dyn: new job \\ arrivals, uncertain processing times \end{tabular}} &
\textcolor{black}{PPO} &
\textcolor{black}{Simulator} &
\textcolor{black}{\begin{tabular}[c]{@{}l@{}} PDRs; \\ DRL methods \end{tabular}} \\
\bottomrule
\end{tabular}%
}
\\[2pt]
\begin{minipage}{\textwidth}
\footnotesize\raggedright
\textit{Notes:} [a] A = action type; O = observation; Dyn = dynamics. PPO = Proximal Policy Optimization, DRL-HH = Deep Reinforcement Learning-based
Hyper-Heuristic, DTMPA = Dynamic Traveling Multi-Maintainer with Alerts, API = Approximate Policy Iteration, DCL = Deep Controlled Learning\textcolor{black}{, MILP = Mixed-Integer Linear Programming, PDR = Priority Dispatching Rule, DRL = Deep Reinforcement Learning}.
\end{minipage}
\end{table}

We briefly discuss representative studies that span discrete versus continuous action spaces and single-agent versus multi-agent control. \citet{Teck2025} use PPO to decide, as order requests arrive, where an inventory rack is stored after picking and whether it is first sent for replenishment in a robotic mobile fulfillment system. The action space is discrete, with infeasible storage zones masked, and the policy is trained in simulation and compared with commonly used storage and replenishment decision rules. \citet{Dehaybe2024} consider continuous ordering decisions under non-stationary demand in a rolling-discounted setting, using PPO trained online in simulation. The problem features continuous (or mixed) actions with full observability, comparisons are made against dynamic programming-based heuristics, illustrating policy-based RL’s suitability for continuous-control inventory problems. \citet{Ying2024} adopt a multi-agent PPO framework for real-time rescheduling under stochastic disruptions. Their setting involves discrete operational decisions with full observability, online training in a simulator, and baselines including metaheuristics and distributed PPO variants, highlighting policy-based RL in time-critical operational control.

In summary, policy-based RL is often preferred in scenarios characterized by high randomness, online or rolling decision-making, and continuous or mixed action spaces. These methods directly optimize the policy objective, demonstrate robustness to sparse or delayed rewards, and can be continuously improved through simulation. They are particularly suitable for problems with continuous or hybrid action spaces, explicit policy constraints or action masking requirements, and multi-agent interactions that involve cooperation or competition.

\subsection{Solving MDPs with actor-critic RL algorithms}
\label{sec_actor}

Given the breadth and heterogeneity of studies employing actor–critic RL algorithms to solve MDPs, we focus on distilling recurring patterns rather than reviewing individual studies one by one. Table~\ref{tab_actorMDP} summarizes recent work in this area. As shown, these methods are most often applied in settings with continuous or mixed action spaces, and in environments characterized by stochastic dynamics, multi-agent interactions, or safety-critical constraints. Application domains include railway scheduling and timetabling, including energy-efficient timetable rescheduling under station disturbances, ride-hailing dispatch, autonomous mobility-on-demand systems, large-scale network control, safety-constrained continuous control, manufacturing resource and worker allocation, and volunteer assignment in disaster response.

A comprehensive review of these studies reveals three key trends. First, algorithm choice is closely tied to the action space: DDPG and PPO-Lagrangian are applied to continuous or mixed-action problems, while Multi-Agent Actor-Critic (MAA2C), multi-agent DDPG, and multi-agent SAC are adopted in decentralized decision-making settings. Second, training methods reflect system accessibility. Online learning is prevalent in operational problems that require real-time adaptation, simulators are used when real-world interaction is costly or risky, and imagined rollouts are leveraged to enhance data efficiency. Third, benchmark baselines remain diverse, ranging from heuristic rules and greedy strategies to metaheuristics, exact methods for small instances (e.g., dynamic programming), MPC, and various PPO variants. 

\begin{table}[h]
\centering
\caption{Classification of papers studying actor-critic RL algorithms for sequential decision-making.}
\label{tab_actorMDP}
\resizebox{\textwidth}{!}{%
\begin{tabular}{@{} r| l l l l l @{}}
\toprule
Publications & Domain / Task & Problem signature & RL method & Training & Baselines \\
\midrule
\citet{Ying2020} & \begin{tabular}[c]{@{}l@{}} Metro rescheduling / \\ Headway control \end{tabular} &
\begin{tabular}[c]{@{}l@{}} A: continuous; O: full; \\ Horizon: finite; Dyn: non-stationary\end{tabular} &
DDPG &
Online & DE; PSO; GA \\
\citet{Liu2020} &
\begin{tabular}[c]{@{}l@{}} Dynamic selective \\ maintenance \end{tabular} &
\begin{tabular}[c]{@{}l@{}}A: large discrete; O: full;\\ Horizon: finite; Dyn: multi-state\end{tabular} &
Actor-critic &
Online &
DP \\
\citet{Lee2021} & \begin{tabular}[c]{@{}l@{}} Disaster response / \\ DSPA\textsuperscript{*} \end{tabular} &
\begin{tabular}[c]{@{}l@{}} A: discrete; O: partial; \\ Horizon: finite; Multi-agent; \\ \textcolor{black}{Dyn: stochastic arrivals} \end{tabular} &
\begin{tabular}[c]{@{}l@{}}Multi-agent \\ actor-critic\end{tabular} &
Online &
\begin{tabular}[c]{@{}l@{}}FCFS; Naive policy; \\ Oracle upper bound \end{tabular} \\
\citet{Zhu2021} & Ride-sourcing & 
\begin{tabular}[c]{@{}l@{}} A: continuous policy over zones;  \\ O: aggregated; Horizon: finite; \\ \textcolor{black}{Dyn: stochastic order matching} \end{tabular} &
Actor-critic & Simulator &\begin{tabular}[c]{@{}l@{}}  Multi-agent \\ RL variants\end{tabular} \\
\citet{Litrain2022} &
\begin{tabular}[c]{@{}l@{}} Railway / \\ Train timetabling \end{tabular} &
\begin{tabular}[c]{@{}l@{}} A: discrete; O: local (actors), \\ global (critic); Horizon: episodic; \\ \textcolor{black}{Dyn: deterministic} \end{tabular}&
MAA2C &
Simulator &
GA; PSO \\
\citet{Ying2022} & \begin{tabular}[c]{@{}l@{}} Urban rail / \\ Coordinated operations \end{tabular}  &
\begin{tabular}[c]{@{}l@{}} A: continuous embedding; O: full; \\ Horizon: infinite, discounted; \\ \textcolor{black}{Dyn: time-varying demand} \end{tabular} &
Multi-agent DDPG & Online & \begin{tabular}[c]{@{}l@{}} GA; \\ Distributed DDPG \end{tabular} \\
\citet{Enders2023} &
\begin{tabular}[c]{@{}l@{}} AMoD / \\ Request assignment \\ \& rejection \end{tabular} &
\begin{tabular}[c]{@{}l@{}} A: discrete; O: full; \\ Horizon: finite episodic; \\ Dyn: stochastic arrivals \end{tabular} &
Multi-agent SAC &
Simulator &
\begin{tabular}[c]{@{}l@{}} Greedy; \\ MPC \end{tabular} \\
\citet{Ma2024} &
\begin{tabular}[c]{@{}l@{}} Large-scale network \\ control /  Traffic, power\textsuperscript{*}, \\ pandemic\textsuperscript{*}, CACC \end{tabular} &
\begin{tabular}[c]{@{}l@{}} A: mainly continuous;\\
O: partial/local; Horizon: finite episodic;\\
Dyn: stochastic exogenous inputs and \\ non-stationary interactions \end{tabular} &
\textcolor{black}{\begin{tabular}[c]{@{}l@{}} Model-based \\ decentralized PPO \end{tabular}}  &
Online &
PPO variants \\
\citet{Jayant2022} &
\begin{tabular}[c]{@{}l@{}} Robot navigation / \\ Safety Gym \\ (PointGoal, CarGoal)\textsuperscript{*} \end{tabular} &
\begin{tabular}[c]{@{}l@{}} A: continuous; O: modified observations; \\ M: learned dynamics ensemble; \\ Horizon: infinite, discounted; \\ Episodes: $T=750$; model rollouts: $H=80$; \\ Dyn: unknown; expected-cost constraint \end{tabular} &
\begin{tabular}[c]{@{}l@{}} MBPPO-Lagrangian \\ (model-based \\ actor--critic) \end{tabular} &
\begin{tabular}[c]{@{}l@{}} Online (on-policy); \\ simulator interaction \\ + model rollouts \end{tabular} &
\begin{tabular}[c]{@{}l@{}} Unconstrained PPO; \\ PPO-Lagrangian; CPO; \\ safe-LOOP (MBRL) \end{tabular} \\
\citet{Panda2024} & \begin{tabular}[c]{@{}l@{}} Manufacturing / \\ Dynamic resource matching \end{tabular}&
\begin{tabular}[c]{@{}l@{}} A: continuous quantities; O: full; \\ Horizon: finite; Dyn: stochastic demand  \end{tabular}&
DKDDPG &
Online &
\begin{tabular}[c]{@{}l@{}} Exact (small); DKQL; \\ DQN; DPG; DDPG \end{tabular}\\
{\hypersetup{citecolor=blue}\citet{LiSeru2025}} & \textcolor{black}{\begin{tabular}[c]{@{}l@{}} Seru production / \\ Dynamic worker allocation \end{tabular}}&
\textcolor{black}{\begin{tabular}[c]{@{}l@{}}A: discrete; O: full; \\ Horizon: finite; Dyn: stochastic arrivals\end{tabular}} &
\textcolor{black}{AC-Ptr} &
\textcolor{black}{Simulator}  &
\textcolor{black}{\begin{tabular}[c]{@{}l@{}} Standard AC; PPO; \\ Approximation algorithm\end{tabular}}  \\
\textcolor{black}{\citet{Meng2025}} &
\textcolor{black}{\begin{tabular}[c]{@{}l@{}} Disaster response / \\ Volunteer assignment \end{tabular}} &
\textcolor{black}{\begin{tabular}[c]{@{}l@{}} A: discrete assignments; O: full; \\ Horizon: finite (24 periods); \\ Dyn: stochastic task arrivals \end{tabular}} &
\textcolor{black}{\begin{tabular}[c]{@{}l@{}} nAPI \\ (actor-critic) \end{tabular}} &
\textcolor{black}{Simulator} &
\textcolor{black}{\begin{tabular}[c]{@{}l@{}} c$\mu$, SPT, SRTF rules; \\ PPO; A2C; GA; PSO; ACO \end{tabular}} \\
\textcolor{black}{\citet{Wang2025}} &
\textcolor{black}{\begin{tabular}[c]{@{}l@{}} Urban rail / Energy-efficient \\ timetable rescheduling \end{tabular}} &
\textcolor{black}{\begin{tabular}[c]{@{}l@{}} A: continuous (coasting speed, \\ dwell time); O: local; Multi-agent; \\ Horizon: episodic; Dyn: station \\ disturbances \end{tabular}} &
\textcolor{black}{DMARL} &
\textcolor{black}{Simulator} &
\textcolor{black}{\begin{tabular}[c]{@{}l@{}} IGA; DDPG; DRLA; \\ Supervised learning; \\ No action \end{tabular}} \\
\bottomrule
\end{tabular}%
}
\\[2pt]
\begin{minipage}{\textwidth}
\footnotesize\raggedright
\textit{Notes:} [a] A = action type; O = observation; Dyn = dynamics. DDPG = Deep Deterministic Policy Gradient, DE = Differential Evolution, PSO = Particle Swarm Optimization, GA = Genetic Algorithm, DP = Dynamic Programming, DSPA = Decentralized Selective Patient Admission, FCFS = First-Come First-Serve, MAA2C = Multi-Agent
Actor-Critic, SAC = Soft Actor-Critic, MPC = Model Predictive Control, DKDDPG = DDPG with Domain-Knowledge Q Penalty, DKQL = Domain Knowledge-informed Q-learning, DQN = Deep Q-Network, DPG = Deterministic Policy Gradient, AC = Actor-Critic, PPO = Proximal Policy Optimization, AC-Ptr = Actor-Critic with Pointer Networks, AMoD = Autonomous Mobility-on-Demand, CACC = Cooperative Adaptive Cruise Control, CPO = Constrained Policy Optimization, MBPPO = Model-Based Proximal Policy Optimization, MBRL = Model-Based Reinforcement Learning, nAPI = Neural Approximate Policy Iteration, SPT = Shortest Processing Time, SRTF = Shortest Remaining Time First, A2C = Advantage Actor-Critic, ACO = Ant Colony Optimization, DMARL = Distributed Multi-Agent Reinforcement Learning, IGA = Iterative Genetic Algorithm, DRLA = Deep Reinforcement Learning Algorithm.
\par
[b] * marks problems or application domains outside transportation, logistics, and closely related operations problems; the corresponding studies are included for methodological illustration.
\end{minipage}
\end{table}

In summary, actor-critic RL algorithms are particularly effective for exploring high-dimensional action spaces, adapting to non-stationary dynamics through continual policy updates, and enabling coordination in multi-agent settings via decentralized actors with a shared or centralized critic. Safe or constrained variants can incorporate penalties or action masks during learning, offering practical advantages in transportation, mobility and process control applications. 

\subsection{Solving MDPs with model-based RL algorithms}
\label{sec_model}

Model-based RL uses a known or learned model of state transitions and rewards for forward planning, imagined experience, or the improvement of policies or value estimates. An explicit model must participate in these computations. Merely training a model-free algorithm in a simulator does not meet this criterion. Table~\ref{tab_modelBasedMDP} covers fifteen formally published studies within transportation and logistics, including freight operations, fleet dispatch, warehouse control, retail replenishment, traffic control, and autonomous driving. The table uses the same columns as Tables~\ref{tab_policeBaedMDP} and~\ref{tab_actorMDP}, and reports the environment model under the problem signature. 

For operational planning with explicit dynamics, \citet{Pan2021} use a multiclass multilane cell transmission model to evaluate freeway control policies during cross-entropy policy search. \citet{Alcaraz2022} address driving, break, and rest decisions along long distance freight routes: a stochastic sample model supports model predictive control, Monte Carlo tree search, and rollout evaluation under uncertain travel and service times. \citet{Yan2023} integrate empirical state transitions, sample average approximation, and SARSA($\Delta$) updates to coordinate charging and order dispatch for an e-hailing electric vehicle fleet. This is a dynamic operational control problem, even though a mathematical programming solver is used within the decision process. For energy efficient routing, \citet{Leech2025} combine a traffic simulation model with tree search and policy and value learning to evaluate route choices.

Learned models also support warehouse and replenishment decisions. \citet{Andersen2019} learn predictive warehouse dynamics and generate short imagined trajectories for policy training in automated storage and retrieval. \citet{Demizu2023} combine offline learning of new product demand with online replenishment planning across retail stores. These studies extend the applications beyond vehicle control to warehouse and retail inventory operations.

For traffic and rail control, \citet{Wang2023TrafficInference} learn a traffic inference model that predicts subsequent states and rewards and generates experience for signal control. \citet{Devailly2022} combine graph representations with latent dynamics as in MuZero and tree search to train coordinated traffic signal policies. The learned policies can also be deployed without online search. \citet{Huang2021} use learned intersection dynamics to produce imagined trajectories for meta-RL, supporting adaptation across intersections. \citet{Zhu2025SLFMLight} introduce SLFMLight, which combines DDPG with a stochastic link flow model for traffic signal control. Predictions of future link densities provide rewards and states for updating the actor and critic networks. \citet{Zhu2025VCTS} learn neural dynamics for virtually coupled train sets and use iterative random shooting to generate trajectories that guide policy search.

Four studies in Table~\ref{tab_modelBasedMDP} address autonomous driving at the vehicle level. We include them as transportation-related sequential decision problems that illustrate how learned models support planning or policy learning when experimentation with physical vehicles is costly or risky. \citet{Wu2022} learn an uncertainty-aware ensemble dynamics model and adaptively truncate synthetic rollouts used for actor-critic learning. \citet{Diehl2022} learn stochastic dynamics from logged data and combine model predictive planning with a learned action prior and value estimates in UMBRELLA. \citet{Gao2025} develop a predictive individual world model that combines latent dynamics with trajectory prediction for surrounding agents. The model generates imagined trajectories for actor-critic training. Think2Drive similarly trains an actor-critic planner using imagined latent trajectories for urban driving \citep{Li2025Think2Drive}.

The number of model-based RL applications identified here depends partly on the search and classification criteria and does not establish that such applications are rare in OR. Requiring explicit OR keywords can miss relevant transportation and logistics studies. Limiting inclusion to formally published journal articles and papers in computing conference proceedings also excludes emerging work available only as preprints. The remaining technical obstacles provide a plausible explanation for differences in the use of model-based RL across applications: stochastic demand and interacting agents complicate model estimation. Model errors can accumulate during planning. Repeated planning adds computational cost to decisions that must be made quickly. These issues motivate short rollouts, uncertainty handling, and compact model representations in the reviewed studies \citep{Janner2019, Wu2022, Diehl2022}. Nevertheless, this search within transportation and logistics cannot establish how widely model-based RL is used across OR. We therefore interpret Table~\ref{tab_modelBasedMDP} as evidence of applications in several areas within the examined transportation and logistics settings, rather than as proof of a general absence of research.

\begin{table}[!htbp]
\centering
\caption{Classification of papers studying model-based RL algorithms for sequential decision-making in transportation and logistics.}
\label{tab_modelBasedMDP}
\linespread{1}\selectfont
\renewcommand{\arraystretch}{0.9}
\resizebox{\textwidth}{!}{%
\begin{tabular}{@{}r|lllll@{}}
\toprule
Publications & Domain / Task & Problem signature & RL method & Training & Baselines \\
\midrule
\citet{Andersen2019} & \begin{tabular}[c]{@{}l@{}} Warehouse storage \\ and retrieval \end{tabular} & \begin{tabular}[c]{@{}l@{}} A: controller dependent \\ O: warehouse observations \\ M: learned dynamics \\ Horizon: episodic \end{tabular} & \begin{tabular}[c]{@{}l@{}} DVAE-2 \\ with PPO \end{tabular} & \begin{tabular}[c]{@{}l@{}} Model learning and \\ imagined policy updates \end{tabular} & \begin{tabular}[c]{@{}l@{}} DQN \\ PPO \end{tabular} \\ \addlinespace[1pt]
\citet{Pan2021} & \begin{tabular}[c]{@{}l@{}} Integrated freeway \\ traffic control \end{tabular} & \begin{tabular}[c]{@{}l@{}} A: traffic control variables \\ O: traffic state \\ M: known multiclass \\ multilane CTM \end{tabular} & \begin{tabular}[c]{@{}l@{}} Cross-entropy \\ policy search \end{tabular} & \begin{tabular}[c]{@{}l@{}} Policy evaluation \\ using the known model \end{tabular} & No optimal control \\ \addlinespace[1pt]
\citet{Alcaraz2022} & \begin{tabular}[c]{@{}l@{}} Long distance \\ freight operations \end{tabular} & \begin{tabular}[c]{@{}l@{}} A: discrete \\ O: route and driver state \\ M: stochastic route model \\ Horizon: finite \end{tabular} & \begin{tabular}[c]{@{}l@{}} MPC with MCTS \\ and rollout \end{tabular} & \begin{tabular}[c]{@{}l@{}} Online planning \\ without prior training \end{tabular} & \begin{tabular}[c]{@{}l@{}} Base policy \\ DQN \end{tabular} \\ \addlinespace[1pt]
\citet{Yan2023} & \begin{tabular}[c]{@{}l@{}} EV fleet charging \\ and order dispatch \end{tabular} & \begin{tabular}[c]{@{}l@{}} A: charging and dispatch \\ O: fleet and order state \\ M: empirical transitions \\ Horizon: finite \end{tabular} & \begin{tabular}[c]{@{}l@{}} SARSA($\Delta$) \\ with SAA \end{tabular} & \begin{tabular}[c]{@{}l@{}} Online value updates \\ and SAA planning \end{tabular} & \begin{tabular}[c]{@{}l@{}} Heuristic \\ ADP \\ Model-free RL \end{tabular} \\ \addlinespace[1pt]
\citet{Wu2022} & \begin{tabular}[c]{@{}l@{}} Autonomous \\ vehicle control \end{tabular} & \begin{tabular}[c]{@{}l@{}} A: continuous \\ O: camera observations \\ M: learned ensemble \\ Horizon: episodic \end{tabular} & UA-MBRL & \begin{tabular}[c]{@{}l@{}} Truncated imagined \\ rollouts and \\ actor-critic updates \end{tabular} & \begin{tabular}[c]{@{}l@{}} Vanilla-SAC \\ MAB-SAC \end{tabular} \\ \addlinespace[1pt]
\citet{Diehl2022} & Automated driving & \begin{tabular}[c]{@{}l@{}} A: continuous \\ O: driving history \\ M: stochastic ensemble \\ Horizon: finite MPC \end{tabular} & UMBRELLA & \begin{tabular}[c]{@{}l@{}} Offline model learning \\ and MPC planning \end{tabular} & \begin{tabular}[c]{@{}l@{}} BC, BC-VAE \\ MBOP, NMP, MPUR \\ Human, Naive \end{tabular} \\ \addlinespace[1pt]
\citet{Wang2023TrafficInference} & \begin{tabular}[c]{@{}l@{}} Traffic signal \\ control \end{tabular} & \begin{tabular}[c]{@{}l@{}} A: discrete \\ O: traffic state \\ M: next-state and \\ reward predictors \end{tabular} & \begin{tabular}[c]{@{}l@{}} Model-based \\ Q-learning \end{tabular} & \begin{tabular}[c]{@{}l@{}} Model learning and \\ generated experience \\ for Q-value updates \end{tabular} & \begin{tabular}[c]{@{}l@{}} Fixedtime, SOTL \\ DQN, A3C \\ Ape-X \end{tabular} \\ \addlinespace[1pt]
\citet{Demizu2023} & \begin{tabular}[c]{@{}l@{}} Retail inventory \\ replenishment \end{tabular} & \begin{tabular}[c]{@{}l@{}} A: discrete order quantities \\ O: inventory and features \\ M: learned demand \\ Horizon: 200 days \end{tabular} & \begin{tabular}[c]{@{}l@{}} BNN/MAML \\ and random \\ shooting \end{tabular} & \begin{tabular}[c]{@{}l@{}} Offline model learning \\ and online planning \end{tabular} & \begin{tabular}[c]{@{}l@{}} Moving average \\ TRPO \\ Perfect information \end{tabular} \\ \addlinespace[1pt]
\citet{Devailly2022} & \begin{tabular}[c]{@{}l@{}} Traffic signal \\ control \end{tabular} & \begin{tabular}[c]{@{}l@{}} A: discrete \\ O: graph traffic state \\ M: latent dynamics \\ Horizon: episodic \end{tabular} & MuJAM & \begin{tabular}[c]{@{}l@{}} Latent dynamics \\ and tree search \\ for policy training \end{tabular} & \begin{tabular}[c]{@{}l@{}} Fixed Time \\ Greedy \\ IG-RL \end{tabular} \\ \addlinespace[1pt]
\citet{Huang2021} & \begin{tabular}[c]{@{}l@{}} Traffic signal \\ control \end{tabular} & \begin{tabular}[c]{@{}l@{}} A: discrete phases \\ O: intersection state \\ M: learned dynamics \end{tabular} & ModelLight & \begin{tabular}[c]{@{}l@{}} Imagined rollouts \\ and meta-RL \end{tabular} & \begin{tabular}[c]{@{}l@{}} MetaLight \\ FRAP++ \\ MAML \end{tabular} \\ \addlinespace[1pt]
\citet{Gao2025} & \begin{tabular}[c]{@{}l@{}} Urban autonomous \\ driving \end{tabular} & \begin{tabular}[c]{@{}l@{}} A: discrete \\ O: partial \\ M: individual world model \\ Horizon: episodic \end{tabular} & PIWM & \begin{tabular}[c]{@{}l@{}} Actor-critic training \\ in imagined latent space \end{tabular} & \begin{tabular}[c]{@{}l@{}} Random, DQN \\ PPO, SAC \\ DreamerV3 \end{tabular} \\ \addlinespace[1pt]
\citet{Li2025Think2Drive} & \begin{tabular}[c]{@{}l@{}} Urban autonomous \\ driving \end{tabular} & \begin{tabular}[c]{@{}l@{}} A: discrete controls \\ O: privileged BEV state \\ M: latent world model \\ Horizon: episodic \end{tabular} & Think2Drive & \begin{tabular}[c]{@{}l@{}} Actor-critic training \\ in imagined latent space \end{tabular} & \begin{tabular}[c]{@{}l@{}} Roach (PPO) \\ TCP \end{tabular} \\ \addlinespace[1pt]
\citet{Leech2025} & \begin{tabular}[c]{@{}l@{}} Energy efficient \\ vehicle routing \end{tabular} & \begin{tabular}[c]{@{}l@{}} A: discrete route choices \\ O: traffic images \\ M: traffic simulation \\ Horizon: episodic \end{tabular} & \begin{tabular}[c]{@{}l@{}} Tree search and \\ policy and value \\ networks \end{tabular} & \begin{tabular}[c]{@{}l@{}} Model rollouts \\ for tree search and \\ network training \end{tabular} & \begin{tabular}[c]{@{}l@{}} Shortest route \\ DQN \end{tabular} \\ \addlinespace[1pt]
\citet{Zhu2025VCTS} & \begin{tabular}[c]{@{}l@{}} Virtually coupled \\ train control \end{tabular} & \begin{tabular}[c]{@{}l@{}} O: positions and velocities \\ M: learned neural \\ train dynamics \end{tabular} & \begin{tabular}[c]{@{}l@{}} IRS-guided \\ policy search \end{tabular} & \begin{tabular}[c]{@{}l@{}} Model learning and \\ trajectory optimization \\ for policy training \end{tabular} & \begin{tabular}[c]{@{}l@{}} Classical control \\ methods \end{tabular} \\ \addlinespace[1pt]
\citet{Zhu2025SLFMLight} & \begin{tabular}[c]{@{}l@{}}Traffic signal\\control\end{tabular} & \begin{tabular}[c]{@{}l@{}}A: continuous green ratios\\O: partial\\M: known SLFM\\Horizon: episodic\end{tabular} & \begin{tabular}[c]{@{}l@{}}SLFMLight\\(DDPG)\end{tabular} & \begin{tabular}[c]{@{}l@{}}Traffic flow predictions\\for actor and critic\\updates\end{tabular} & \begin{tabular}[c]{@{}l@{}}FixedTime, MaxPressure\\DDPGLight\\CoLight, MPLight\end{tabular} \\ \addlinespace[1pt]
\bottomrule
\end{tabular}%
}
\\[2pt]
\begin{minipage}{\textwidth}
\footnotesize\raggedright
\textit{Notes:} [a] A = action type, O = observation, M = environment model,
CTM = Cell Transmission Model, MCTS = Monte Carlo Tree Search,
BNN = Bayesian Neural Network, IRS = Iterative Random Shooting, BEV = Bird's-Eye View,
SLFM = Stochastic Link Flow Model, ADP = Approximate Dynamic Programming, MPC = Model Predictive Control, SAA = Sample Average Approximation.
Observation descriptions are used where full or partial observability is not established. Other method and baseline names follow the original studies.
\end{minipage}
\end{table}

\subsection{Discussion \& future research directions}
\label{sec:open}

We now consolidate the main insights across RL algorithm families, uncover their fundamental limitations, and point to promising directions for future research.

\textbf{Insights across RL algorithm categories.} 
The analyses above show that value-based methods are most often applied in environments  with discrete action spaces, full observability, and finite or segmented decision horizons. In contrast, settings characterized by strong stochasticity, rolling online control, continuous or hybrid action spaces, or multi-agent interactions tend to favor policy-based or actor-critic approaches. When interactions are costly or risky, the environment is non-stationary, safety constraints must be observed, or historical logs are available, model-based RL algorithms become a more appropriate choice. Regarding the performance-cost trade-offs, model-based and offline methods typically achieve higher sample efficiency, but this comes with substantial training and engineering overhead due to model learning, simulator construction, and uncertainty quantification. These approaches also face scalability challenges, such as high-dimensional encoder-decoder architectures, long-horizon credit assignment, and sparse or delayed reward signals.

Overall, model-free and model-based RL represent two distinct paradigms that differ fundamentally in how they interact with and learn from the environment. In model-free RL frameworks, the agent directly maps observations to actions based on past experiences. It learns to act by trial and error, updating its behavior using reward prediction errors. This makes model-free methods effective in highly complex or unpredictable environments where modeling transitions is infeasible or unnecessary and where real-time decision-making can be supported through abundant interactions or reliable simulators. Model-based RL, in contrast, constructs an internal model of state transitions and rewards, enabling planning through simulated trajectories. These methods excel in structured or deterministic environments, where long-term predictions are more reliable. Model-based RL algorithms are often favored when sample efficiency is critical, or prior knowledge about the system’s structure is available. They allow agents to learn faster with fewer interactions, albeit at the cost of higher training and modeling complexity.

Ultimately, the choice between model-free and model-based RL algorithms depends on the nature of the environment and practical considerations such as data availability, computational resources, and operational requirements such as safety, interpretability, or planning depth.

\textbf{Limitations and future research directions.} Model-free RL algorithms have made remarkable progress in recent years. However, their low sample efficiency presents a significant limitation, that is, they typically require extensive real-world interactions with the environment to collect sufficient training data and learn effective policies. This constraint hinders their applicability to real-world domains where interactions are costly, time-consuming, or risky, confining their practical use largely to settings with accessible and reliable simulators. Moreover, model-free methods lack the ability to simulate hypothetical future trajectories, relying solely on real interactions rather than imagined rollouts. This limitation not only reduces the efficiency of exploration but also restricts the agent’s ability to anticipate long-term consequences. For challenges specific to value-based model-free RL algorithms, we refer to \citet{Park2024}.

As for model-based RL algorithms, despite their potential to improve sample efficiency and enable planning through simulated rollouts, they also face important limitations. First, they rely on accurate modeling of environment dynamics. In systems with complex, discontinuous, or stochastic behavior, learning such models can be difficult and may outweigh potential benefits, making these algorithms less effective. Second, model-based RL introduces additional sources of estimation error: not only from value function approximation but also from the learned dynamics model. When both are imperfect, compounded errors can significantly degrade policy quality. Third, these algorithms are highly sensitive to model accuracy. An inaccurate model can mislead the agent during planning, resulting in poor generalization or unsafe decisions, especially in long-horizon tasks where prediction errors accumulate. Therefore, reliable model learning is essential, yet often challenging in high-dimensional or partially observable environments.

Looking forward, an important future research direction is the development of hybrid frameworks that enable flexible coordination between model-free and model-based RL algorithms within a single decision process, thereby fully combining their advantages. Specifically, one direction worth exploring is the development of centralized decision-making architectures that can dynamically coordinate the use of model-free or model-based RL algorithms based on the characteristics of the problem, such as decision-making stages, uncertainty, and dynamics. Such centralized coordination frameworks could allow agents to exploit model-based planning for rapid learning, and then switch to model-free strategies for robustness and policy generalization, or alternate between the two to improve both computational efficiency and solution quality.

The resources and evaluation criteria in Section~\ref{sec_benchmarking} also support comparisons in dynamic sequential decision problems. Such comparisons should use the same demand and travel time scenarios and report operational costs or rewards, service levels, and constraint violations.

\section{RL: Enhancer of solution methods for COPs}
\label{Sec_algorithm}

In this section, we focus on the research stream in which RL serves as a tool to enhance OR solution methods for static COPs. Here, COPs refer to problems where the optimization instance, including any uncertainty model, is specified before the solution process begins, and the goal is to construct a solution (e.g., a route or schedule) for that instance, rather than to learn a policy for repeated decisions in an evolving operational system. A growing body of work treats RL as an end-to-end, independent solution method for solving COPs. Although many RL-based COP methods generate solutions through a sequence of construction decisions or guide iterative search decisions, these decisions belong to the algorithmic solution process for fixed problem inputs; they should not be confused with temporal decisions in an evolving operational system. This line of research leverages the exploration capabilities of RL to learn a mapping from problem inputs to solutions and to directly search the solution space. In addition, RL can be integrated with heuristic or exact algorithms to tailor solution methodologies for COPs, typically by guiding parts of the search process, such as selecting promising neighborhoods or cutting planes.

Section~\ref{sec_RLsolver} reviews studies in which RL serves as an end-to-end and independent solution method, Section~\ref{sec_RLheuristic} examines RL-enhanced heuristic algorithms, and Section~\ref{sec_RLexact} discusses RL-enhanced exact algorithms. In each subsection, we analyze the relevant literature, summarize the associated challenges and limitations, and propose potential directions for future research. Section~\ref{sec_benchmarking} provides common evaluation criteria and resources for reproducible experiments.

\subsection{RL as an end-to-end and independent solution method}
\label{sec_RLsolver}

When RL is used as an end-to-end solution method for static COPs, it learns a direct mapping from fixed problem inputs to feasible solutions, which are typically decoded using greedy, sampling-based, or beam-search strategies. We first review classic COPs such as the Traveling Salesman Problem (TSP), Vehicle Routing Problem (VRP), and Minimum Vertex Cover Problem (MVCP) in Section~\ref{sec_RLsolverCCOP}, as early research on end-to-end RL methods primarily focused on these problems. We then examine applications in more domain-specific combinatorial optimization contexts in Section~\ref{sec_RLsolverOCOP}. Finally, we discuss key challenges and potential future research directions in Section~\ref{sec_RLsolverFuture}.

In this research context, terms such as \emph{Pointer Networks} and \emph{Transformers} refer to the neural network \emph{architectures} used to parameterize policies or value functions. These choices are orthogonal to the underlying RL algorithms (value-based, policy-based, actor-critic) and to whether the method is model-free or model-based. These architectures are frequently adopted in this line of research, as they naturally align with the structural properties of classic COPs with routing decisions and support the construction of high-quality solutions in an end-to-end manner. Specifically, pointer networks~\citep[see,][]{Vinyals2015} are sequence-to-sequence models whose decoder \enquote{points} to input items, making them particularly suitable for constructing permutations (e.g., routing or scheduling). 
Transformers~\citep[see,][]{Vaswani2018, Kool2018} use self-attention mechanisms to encode variable-size sequences or sets and support autoregressive decoding, offering improved scalability and performance on large-scale problem instances.

\subsubsection{Literature analysis of classic COPs}
\label{sec_RLsolverCCOP}

Focusing on RL as an end-to-end, independent solution method for solving classic COPs, Table~\ref{tab:RLOptimization} classifies representative studies along four dimensions: (i) \emph{RL paradigm and architecture} (e.g., as policy gradient, value-based, actor-critic; Pointer Networks, attention/Transformers); (ii) \emph{Solved COPs}; (iii) \emph{Baselines}, i.e., benchmark algorithms used for comparison (e.g., Christofides, OR-Tools, Concorde); and (iv) \emph{Problem scale}. As shown in Table~\ref{tab:RLOptimization}, graph-structured problems, especially TSP and VRP, are the most frequently studied, and classic heuristics are commonly used as baselines. Across these studies, three recurring patterns emerge:
(i) routing problems often use policy gradient methods with pointer or Transformer decoders, while general graph problems prefer value-based methods with graph neural network (GNN) encoders;
(ii) recent studies shift from pointer networks to Transformer architectures as problem scale increases;
(iii) baseline selection varies across domains, making it difficult to compare results across studies due to the lack of standardized evaluation protocols.

\begin{table}[h]
\centering
\caption{Classification of papers studying RL as an end-to-end and independent solution method for solving classic COPs.}
\label{tab:RLOptimization}
\resizebox{\textwidth}{!}{%
\begin{tabular}{@{} r|l l l c @{}}
\toprule
Publications & RL paradigm / architecture & Solved COPs & Baselines & Problem scale \\
\midrule
\citet{Nazari2018}  & Policy gradient; Pointer-like seq2seq 
                   & VRP 
                   & Classic heuristics and OR-Tools 
                   & 10--100 nodes \\
\citet{Kool2018}  & Policy gradient; Transformer
                   & TSP, VRP
& \begin{tabular}[c]{@{}l@{}} OR-Tools, classical heuristics,\\  Concorde (exact), and GUROBI\end{tabular}
& 20--100 nodes
  \\
\citet{Dai2017}     & Value-based; GNN  
                   & TSP, MVCP\textsuperscript{*}, MAXCUT\textsuperscript{*} 
                   & Approximation and heuristics 
                   & 50--100 nodes \\
\citet{Barrett2020} & Value-based; GNN + vertex flipping 
                   & MAXCUT\textsuperscript{*} 
                   & \begin{tabular}[c]{@{}l@{}} An RL-based heuristic \\ and a greedy algorithm \end{tabular}
                   & 20--200 vertices \\
\citet{Li2021}    & Policy gradient; Pointer Network
                   & MOTSP 
                   & NSGA-II, MOEA/D, MOGLS 
                   & 20--500 nodes \\                
\citet{Jin2023}     & Policy gradient; Transformer pointer 
                   & TSP 
                   & Existing deep-learning algorithms 
                   & 20--500 nodes \\
\bottomrule
\end{tabular}%
}
\\[2pt]
\begin{minipage}{\textwidth}
\footnotesize\raggedright
\textit{Notes:} [a] TSP = Traveling Salesman Problem, VRP = Vehicle Routing Problem, GNN = Graph neural networks, MVCP = Minimum Vertex Cover Problem, MAXCUT = Max-Cut Problem, MOTSP =  Multiobjective Traveling Salesman Problem.
\par
[b] * marks problems or application domains outside transportation, logistics, and closely related operations problems; the corresponding studies are included for methodological illustration.
\end{minipage}
\end{table}

We now discuss representative formally published studies. For instance, \citet{Nazari2018} adapted pointer-like sequence models to solve VRPs with both static and dynamic elements, while \citet{Dai2017} combined RL with graph embeddings to address a range of graph-structured COPs (e.g., Minimum Vertex Cover Problem, Max-Cut Problem, and TSP) by learning a greedy construction policy that acts as a metaheuristic. \citet{Li2021} proposed a deep RL method for the multi-objective TSP. The multi-objective optimization problem is decomposed into several single-objective subproblems, where each is formulated as an RL decision process under a given weight vector. The trained model functions as a black-box metaheuristic that generates solutions without per-instance reoptimization. Transformer-based methods further improve scalability and sample efficiency \citep[e.g.,][]{Kool2018}. More recently, \citet{Jin2023} proposed an end-to-end solution method named Pointerformer, which can solve instances with up to several hundred nodes. Beyond TSP, \citet{Barrett2020} developed a deep Q-learning approach for graph-structured COPs that integrates GNN encoders with local-search-inspired mechanisms (e.g., vertex flipping and reward shaping). Applied to the Max-Cut problem, their method outperforms prior RL-based baselines and can also be combined with other search heuristics.

\subsubsection{Literature analysis of solving other COPs}
\label{sec_RLsolverOCOP}

Beyond classic COPs, a growing body of research applies RL to other combinatorial problems in domain-specific contexts. Table~\ref{tab_RLOptimizationApp} classifies related studies along four dimensions:
(i) \emph{Research domain};
(ii) \emph{Addressed problem};
(iii) \emph{RL paradigm and architecture}; and
(iv) \emph{Baselines}, i.e., benchmark algorithms used for comparison.

\begin{table}[h]
\centering
\caption{Classification of papers studying RL as an end-to-end and independent solution method for solving domain-specific COPs.}
\label{tab_RLOptimizationApp}
\resizebox{\textwidth}{!}{%
\begin{tabular}{@{} r|l l l l l @{}}
\toprule
Publications & Domain & Problem & RL paradigm / architecture & Baselines  \\
\midrule
\textcolor{blue}{\citet{Wei2017}} & \textcolor{black}{Information retrieval\textsuperscript{*}} & \textcolor{black}{Learning to rank} & \textcolor{black}{Policy gradient / Linear softmax policy} & \textcolor{black}{Learning-to-rank algorithms}\\
\textcolor{black}{\citet{Monaci2024}} & \textcolor{black}{Manufacturing} & \textcolor{black}{Job-shop scheduling} & \textcolor{black}{Actor-critic (PPO) / Double LSTM} & \textcolor{black}{Heuristics, CPLEX, and DRL}\\
\citet{Ma2025} & Supply chain & Product design change & Actor-critic / Bi-level SAC   &  DRL and heuristics\\
\citet{Su2025}& Public transport &  Route-frequency design & PPO / Action masking & Heuristics and Q-learning\\
\citet{Ding2025} & Electric robots & Charging-robot scheduling & Policy-gradient / Transformer & GUROBI, heuristics, and DRL\\
\citet{Tian2025} & Logistics & Service scheduling & Actor-critic / Transformer & Heuristics and RL baselines\\
\citet{Teusch2025} & Urban mobility planning & Facility location & Value-based / DDQN & Simulation and PPO\\
\bottomrule
\end{tabular}%
}
\\[2pt]
\begin{minipage}{\textwidth}
\footnotesize\raggedright
\textit{Notes:} [a] SAC = Soft Actor-Critic, DRL = Deep Reinforcement Learning, DDQN = Double DQN, PPO = Proximal Policy Optimization, LSTM = Long Short-Term Memory.
\par
[b] * marks problems or application domains outside transportation, logistics, and closely related operations problems; the corresponding studies are included for methodological illustration.
\end{minipage}
\end{table}

In supply chains, \citet{Ma2025} modeled product-design change as a bilevel joint optimization problem and solve it end-to-end with a bilevel DRL method, benchmarking against conventional DRL and heuristic baselines. In manufacturing, \citet{Monaci2024} formulated the job-shop scheduling problem as a Markov decision process in which a PPO-based actor-critic agent with recurrent networks selects the next operation to schedule. The trained policy acts as a greedy construction heuristic and is compared with dispatching rules, an adaptive heuristic, CPLEX, and another DRL method. \citet{Wei2017} similarly construct a document ranking for each query position by position and train the ranking policy with REINFORCE. RL-based algorithmic frameworks are also applied to areas related to the public transport network design \citep[e.g.,][]{Su2025}, charging and scheduling for electric robots or vehicles \citep[e.g.,][]{Ding2025}, logistics \citep[e.g.,][]{Tian2025}, and urban mobility planning \citep[e.g.,][]{Teusch2025}. Studies that use RL to learn operational policies for dynamic variants of such problems, for example, job-shop scheduling with new job arrivals \citep{Yu2026}, are reviewed in Section~\ref{Sec_MDP}.

In summary, across various domains, policy-gradient, actor-critic, and PPO variants are commonly used to construct schedules and designs step by step, as in job-shop scheduling, charging robot scheduling, and route-frequency design. Value-based methods such as DDQN are applied to discrete location decisions. Benchmark methods vary considerably across studies, including heuristics, simulators, commercial solvers such as GUROBI and CPLEX, and RL baselines, highlighting the need for standardized evaluation protocols. Several studies benchmark primarily against other RL methods; incorporating OR baselines, such as OR-Tools, CP-SAT, or domain-specific heuristics, would strengthen the evidence for the effectiveness of the proposed RL-based methods.

\subsubsection{Discussion \& future research directions}
\label{sec_RLsolverFuture}

When RL is used as an end-to-end, independent solution method for classic and domain-specific COPs, three commonalities emerge. First, routing problems typically use policy-gradient decoders with Pointer or Transformer architectures, while general graph problems such as the Minimum Vertex Cover Problem and Max-Cut Problem often pair value-based methods with GNN encoders. Domain-specific scheduling and design problems frequently adopt actor-critic or PPO. Second, in routing studies, larger instance sizes are usually addressed with Transformer architectures rather than Pointer networks. Third, benchmark solution methods remain heterogeneous across studies, mixing classical heuristics, simulators, commercial solvers, and RL baselines, which complicates fair comparison.

Despite rapid progress, several challenges remain. First, RL-based standalone solvers often struggle with scalability and generalization, particularly for large-scale instances (e.g., thousands of nodes) and typically require retraining when scale, topology, costs, or constraints change. Second, training times grow quickly with instance size. Third, most studies focus on simplified problems, with real-world constraints and uncertainties often overlooked. Fourth, inconsistent benchmarking and the lack of standardized evaluation protocols hinder comparability. The limited scalability of the RL methods arises mainly from two issues. First, typical encoders rely on $O(n^2)$ attention or high-order message passing, and decoders use step-by-step construction, both of which scale poorly and lead to long training times. Second, feasibility is rarely enforced as a hard constraint, and classical heuristics are not integrated into the generation step, making it difficult for the RL model to learn stable and near-optimal construction rules. Additionally, training is usually performed on a fixed input distribution, preventing transfer across scales or objective coefficients and limiting the reuse of learned policies or value functions.

To address these limitations and challenges, future research could (i) develop more efficient encoding and decoding schemes to reduce training time; (ii) combine RL with heuristic or exact algorithms to improve computational efficiency; (iii) address problems that better reflect real-world complexity; (iv) enhance generalization by employing dataset randomization and systematically evaluating transfer performance across different distributions; and (v) establish unified evaluation benchmarks that incorporate robust OR tools and problem-specific heuristics, following the criteria in Section \ref{sec_benchmarking}.

\subsection{RL-enhanced heuristic algorithms}
\label{sec_RLheuristic}

In this subsection, we analyze algorithms that integrate RL with heuristics. We then outline limitations, challenges, and potential directions for future research.

\subsubsection{Literature analysis}

In this stream of research, RL plays five main roles in algorithms that combine RL with heuristics. Table~\ref{tab_RLHeuristic} summarizes the literature where RL is incorporated into heuristic methods for solving COPs, outlining the role of RL, the RL technique applied, the heuristic component, and the OR problems considered. The studies can be grouped into five categories according to the role of RL: (i) RL-driven cooperative search; (ii) RL constructs initial solutions; (iii) RL selects the most promising heuristic; (iv) 
RL selects the most promising operators or neighborhood; and (v) RL helps to tune intensities/quantities online.

\begin{table}[h]
\centering
\caption{Classification of papers studying the combinations of RL with heuristic algorithms for solving COPs. }
\label{tab_RLHeuristic}
\resizebox{\textwidth}{!}{%
\begin{tabular}{@{} r|l l l l l @{}}
\toprule
Publications & Detailed roles of RL & RL method & \begin{tabular}[c]{@{}c@{}} OR \\ algorithm \end{tabular} & OR problem  \\
\midrule
\citet{Martin2016}& \begin{tabular}[c]{@{}l@{}} Coordinate multi-agent \\ cooperative search \end{tabular}
  & RL-driven cooperation
 & Heuristics  & PFSP, CVRP\\
\citet{Benlic2017} & Help selecting the perturbation &MAB & BLS &  VSP\textsuperscript{*}\\
\citet{Deudon2018}& Help generating initial solutions  & Pointer Network & LS & TSP\\
\citet{Mosadegh2020} & RL selects  heuristics & Tabular Q-learning & SA
& SMMALSP\\
\citet{Lamghari2020} & Selects low-level heuristics & Bandit (choice-function) & 27 heuristics & SMPSP\\
\citet{Ma2021}& Help selecting the neighborhood  & DQN & LS & DPDP\\
\citet{Meijer2021} & \begin{tabular}[c]{@{}l@{}} Constructs initial solutions \\ (SDP-guided rounding) \end{tabular}  & Tabular Q-learning & \begin{tabular}[c]{@{}l@{}} SDP-based \\ rounding \end{tabular}& QCCP\\
\citet{Alicastro2021} & Help selecting the neighborhood & Tabular Q-learning & LS & MSP\\
\citet{Brammer2022} & Help generating initial solutions & PPO & LS & PFSP\\
\citet{Zhang2022} &\begin{tabular}[c]{@{}l@{}}  RL selects parameterized \\ low-level heuristics \end{tabular}  & DDQN
& \begin{tabular}[c]{@{}l@{}}Heuristics \\ with rules\end{tabular}  & COPs 
 \\
\citet{Li2022} & Help selecting the neighborhood & Probability learning & TS & AGAP \\
\citet{Kallestad2023} & Help selecting operators & PPO & ALNS & COPs\\
\citet{Zhang2023} &Help selecting operators & DQN & ALNS & SFTRP \\
\citet{Karimi2023} & \begin{tabular}[c]{@{}l@{}} Help selecting perturbation \\ operators and strength \end{tabular}
 & Tabular Q-learning & IG & PFSP\\
\citet{Li2024} & \begin{tabular}[c]{@{}l@{}}  Guides construction of initial \\ solutions and steer subsequent search \end{tabular} & Tabular Q-learning & LS & CDAP\\
\citet{Li2024COR} &Help selecting operators & Q-learning & GA & FAP \\
\citet{Wu2024} & \begin{tabular}[c]{@{}l@{}}Guide TS to focus on \\ promising add/drop moves \end{tabular} & \begin{tabular}[c]{@{}l@{}}  Learning-automata \\ (probability-matrix) \end{tabular}& TS & \begin{tabular}[c]{@{}l@{}} Clustered Orienteering \\ Problem \end{tabular}\\
\citet{Lu2024} & Help selecting the neighborhood & Probability learning & LS
 & IUCP\textsuperscript{*}\\
\citet{Cui2024} &Help selecting operators & Tabular Q-learning & GA  & PSP  \\
\citet{Zou2024}  &Help selecting the neighborhood & Tabular Q-learning & VND  & LLRP  \\
\citet{Zhang2025}& \begin{tabular}[c]{@{}l@{}} Hyper-heuristic controller that adjusts \\  the number of non-dominated solutions\end{tabular} & Tabular Q-learning & GA & \begin{tabular}[c]{@{}l@{}}  Formation and \\ Scheduling Optimization\end{tabular}\\
\citet{Rolim2025} &Help selecting the neighborhood & Q-learning & \textcolor{black}{SA-based SLS, ALNS} & PBSP \\
\citet{Zhao2025} & Help selecting the neighborhood & Probability learning & VNS & k-CMBCP\textsuperscript{*} \\ 
\citet{An2026} &Help selecting the neighborhood  & MAB&VNS &CBSP \\
\bottomrule
\end{tabular}%
}
\\[2pt]
\begin{minipage}{\textwidth}
\footnotesize\raggedright
\textit{Notes:} [a] PFSP = Permutation
Flowshop Scheduling Problem, CVRP = Capacitated Vehicle Routing Problem, MAB = Multi-Armed Bandit, BLS = Breakout Local Search, VSP = Vertex Separator Problem, LS = Local Search, TSP = Traveling Salesman Problem, SA = Simulated Annealing, SMMALSP = Stochastic Mixed-Model Assembly Line Sequencing Problem, SMPSP = Stochastic open-pit Mine Production Scheduling Problem, DQN = Deep Q-Network, DPDP = Dynamic Pickup and Delivery Problem, SDP = Semidefinite Programming, QCCP = Quadratic Cycle Cover Problem, MSP = Machine Scheduling Problem, PPO = Proximal Policy Optimization, DDQN = Double Deep Q Network, TS = Tabu Search, AGAP = Airport Gate Assignment Problem, ALNS = Adaptive Large Neighborhood Search, SFTRP = Synchromodal Freight Transport Re-planning Problem, IG = Iterated Greedy, CDAP = Cross-dock Door Assignment Problem, \textcolor{black}{FAP = Fleet Allocation Problem,} IUCP = Maximum Independent Union of Cliques Problem, GA = Genetic Algorithm, PSP = Production Scheduling Problem, VND = Variable Neighborhood Descent, LLRP = Latency Location Routing Problem, SLS = Stochastic Local Search, PBSP = Parallel Batch Scheduling Problem, VNS = Variable Neighborhood Search, k-CMBCP = k-clustering Minimum Biclique Completion Problem, CBSP = Customized Bus Scheduling Problem.
\par
[b] * marks problems or application domains outside transportation, logistics, and closely related operations problems; the corresponding studies are included for methodological illustration.
\end{minipage}
\end{table}

\textbf{(i) RL-driven cooperative search.} In this algorithmic framework, multiple agents run distinct metaheuristics/local searches and cooperate via RL-guided information sharing. RL adapts agents’ behaviors (e.g., which peer/solution pattern to trust or exchange) so that the overall search escapes local minima and balances diversification/intensification. For example, \citet{Martin2016} propose a general multi-agent cooperative search framework where each agent runs a different metaheuristic/local-search configuration. Agents communicate asynchronously and adapt via a cooperation protocol based on RL and pattern matching: they identify good patterns and share them, while RL adjusts how agents exploit shared information.

\textbf{(ii) RL constructs initial solutions.} In this research line,  RL constructs an initial solution that is then refined by heuristics \citep[e.g.,][]{Deudon2018, Brammer2022, Li2024}. RL thus provides high-quality starting points and reduces the workload of the post-improvement phase. For instance, \citet{Deudon2018} proposed an algorithm that integrates a policy gradient-based RL with a local search heuristic. A neural network (i.e., Pointer Network with Attention) is used to generate many initial solutions and select the best one for the TSP. Subsequently, a 2-opt local search heuristic is applied as a post-processing step to further refine the solution and reduce the tour length. To solve the Cross-dock Door Assignment Problem, \citet{Li2024} proposed a hybrid algorithm combining RL and a local search algorithm: a Q-learning agent constructs an initial solution, which is then refined through local search in an iterative loop of \enquote{Q-learning-based construction} to \enquote{improvement by local search} to \enquote{reward update}.

\textbf{(iii) RL selects the most promising heuristic.} RL can be used to select the most promising heuristics \citep[see,][]{Mosadegh2020, Lamghari2020}. \citet{Mosadegh2020} proposed a Hyper Simulated Annealing framework, where a tabular Q-learning algorithm is embedded to select the most suitable heuristic online during the search process. Here, RL acts as a controller that dynamically chooses which heuristic to apply at each iteration based on learned Q-values.

\textbf{(iv) RL selects the most promising operators or neighborhood. } In this type of algorithm, RL helps a heuristic select the most promising operator or neighborhood in each iteration \citep[e.g.,][]{Benlic2017, Ma2021, Alicastro2021, Li2022, Kallestad2023, Zhang2023, Li2024COR, Wu2024, Cui2024, Zou2024, Lu2024, Rolim2025, Zhao2025, An2026}. Specifically, to solve the Dynamic Pickup and Delivery Problem (DPDP), \citet{Ma2021} employed RL as an intelligent scheduler for a local search heuristic, using the Deep Q-Network to dynamically select the most promising local perturbation at each decision point. This enhances global search capability and improves convergence quality. For COPs in general, \citet{Kallestad2023} developed a selection hyperheuristic framework that integrates DRL into the Adaptive Large
Neighborhood Search (ALNS) algorithm, where DRL replaces ALNS’s adaptive weights
and is used to choose which destroy/repair operator, or an extra deterministic heuristic, to apply at each iteration. For the Clustered Orienteering
Problem, \citet{Wu2024} designed an algorithm combining RL and tabu search (TS). TS serves as the main search engine, while RL updates a probability matrix that guides neighborhood evaluation, filtering out unpromising neighbors and reinforcing the search.

\textbf{(v) RL helps tuning intensities/quantities online.}
In \textcolor{blue}{this} line of work, RL tunes heuristic intensities or quantities online (e.g., perturbation strength, number of tasks removed, or parameters of low-level heuristics), allowing the search to modulate diversification and intensification on the fly \citep[see,][]{Zhang2022, Karimi2023, Zhang2025}. For instance, \citet{Zhang2022} proposed a DRL-based hyper-heuristic framework for COPs under uncertainty. The approach augments traditional hyper-heuristics with a data-driven heuristic selection module, where DRL is used to select parameter-controlled low-level heuristics to improve performance under uncertainty across diverse problem domains. For the Permutation Flowshop Scheduling Problem, \citet{Karimi2023} embedded Q-learning into an Iterated Greedy (IG) algorithm, where the perturbation operators and their intensities (i.e., the number of jobs removed and reinserted per iteration) are selected online during the search process. In this framework, RL does not directly construct solutions but guides the perturbation phase of IG. The reward function combines local improvements with progress toward the global best solution. \citet{Zhang2025} proposed a multi-objective cooperative co-evolution algorithm combining RL and a genetic algorithm to tackle planning in a hybrid seru system (i.e., a production mode combining seru cells with a flow line). A hypervolume-based Q-learning agent serves as a hyper-heuristic controller that adaptively adjusts the number of cooperative non-dominated solutions participating in coevolution at each period.

\subsubsection{Discussion \& future research directions}

Across the five research streams, the \emph{RL selects the most promising operators or neighborhood} category currently dominates. These methods are relatively easy to integrate with existing metaheuristics (e.g., VNS/ALNS/GA/LS/TS) and often yield noticeable improvements in computational efficiency and/or solution quality. By contrast, the \emph{RL-driven cooperative search}, \emph{RL constructs initial solutions}, and \emph{RL selects the most promising heuristic} streams are less common in the literature, as they require deeper integration and are more sensitive to reward shaping and feasibility safeguards. The \emph{RL helps tuning intensities/quantities online} stream is gaining traction, yet it remains underexplored beyond small and narrowly defined parameter sets.

Specifically, only a limited number of studies have explored the RL-driven cooperative search method, in which each agent runs a metaheuristic. Its effectiveness has only been demonstrated on the Permutation Flowshop Scheduling Problem (PFSP) and the Capacitated Vehicle Routing Problem (CVRP). Its limitations are largely due to the strategy of aggregating locally promising solution fragments from each agent. This approach is less suitable for problems that involve long decision sequences or complex route configurations. Moreover, using only a few agents results in marginal performance gains over single-agent methods, while employing many agents considerably increases communication and computational costs. To address these challenges, future research could explore integrating a broader range of heuristics and designing pattern-discovery mechanisms that coordinate more effectively across different algorithms. It would also be valuable to investigate strategies for extracting and filtering higher-level structural features from candidate solutions, rather than relying solely on local fragments. Finally, developing more efficient concurrent implementations and communication strategies, along with systematic evaluations of how solution quality scales with the number of agents, would contribute to improving scalability and overall performance.

In the domain where RL constructs initial solutions that are then refined by heuristics, current research is largely limited to relatively simple problems, such as single-vehicle orienteering \citep{Li2024TS}. Future research could extend these methods to larger scales and more complex variants, including multi-vehicle orienteering with capacity and time-window constraints. Within the limited literature on RL for heuristic selection and for online tuning of intensities or quantities, existing approaches are mostly problem-specific, and their generalizability to broader families of COPs remains uncertain. Future research should establish cross-domain benchmarks, evaluate transferability across instance sizes and input distributions, and compare RL-augmented methods against non-RL heuristics under consistent computational budgets. 

In the field of using RL to select the most promising operators or neighborhoods, several recurring limitations have been identified. One commonly cited problem is the lack of quantifiable optimality gaps for the solutions obtained. However, this is a limitation of heuristic methods rather than of RL itself. Another frequently noted concern is that RL training times can be relatively long and may even exceed the time required to solve the problem using heuristics alone. However, the reported gains depend on the baseline and the computational budget, and improved solution quality does not necessarily imply shorter solution time or an advantage over well-tuned heuristics. Many existing studies also highlight that RL-based approaches for operator or neighborhood selection often lack generalization capability. When applied to new problem classes, these methods usually require retraining, and their performance cannot be guaranteed. This challenge stems from the nature of heuristics themselves, as effective operators and neighborhood structures are generally problem-specific and do not transfer reliably across different COP families.

In summary, the reviewed studies report gains from RL-enhanced heuristics under specific test conditions, but these findings do not establish a general advantage over well-tuned heuristic methods. Optimality guarantees remain elusive, generalization to unseen instances or new problem classes is often limited, and retraining is usually required when problem characteristics such as scale, distribution, or constraint structure change. Importantly, these shortcomings are largely attributable to the heuristic backbone rather than to RL itself. An interesting future research direction is to investigate whether RL can help mitigate these long-standing limitations, particularly by improving the generalization behavior of heuristic methods and enabling the creation of more transferable operator- or neighborhood-selection strategies.

\subsection{RL-enhanced exact algorithms}
\label{sec_RLexact}

We summarize and analyze the literature that integrates RL with various types of exact algorithms and then discuss key challenges and future research directions.

\subsubsection{Literature analysis}
\begin{table}[h]
\centering
\caption{Classification of papers studying the combinations of RL with exact algorithms for solving COPs.}
\label{tab_RLExact}
\resizebox{\textwidth}{!}{%
\begin{tabular}{@{} r|l   l l @{}}
\toprule
Publications & Roles of RL  & OR algorithm & OR problem  \\
\midrule
\citet{Cappart2019} &  Learn DD ordering & DD bounding
    & Max-Cut\textsuperscript{*} and MIS\textsuperscript{*}
\\
\citet{Tang2020}& Select Gomory cuts
  & \begin{tabular}[c]{@{}l@{}} Cutting-plane, \\ Branch-and-Cut\end{tabular}
 & \begin{tabular}[c]{@{}l@{}} Packing; binary packing; \\ production planning; Max-Cut\textsuperscript{*} \end{tabular}
\\
\citet{Ichnowski2021}&  Tune ADMM parameters  & ADMM & \begin{tabular}[c]{@{}l@{}} QP: control; Huber fitting\textsuperscript{*}; \\ SVM\textsuperscript{*}; Lasso regression\textsuperscript{*}; \\ portfolio optimization\textsuperscript{*}; \\ equality-constrained / random QPs\textsuperscript{*} \end{tabular}\\
\citet{Cappart2021} & Learn CP branching & CP & \begin{tabular}[c]{@{}l@{}} TSP with time windows; 0--1 knapsack;\\ 4-moments portfolio optimization\textsuperscript{*} \end{tabular}\\
\citet{Cappart2022}& Learn DD ordering &  DD bounding and B\&B
    & MIS\textsuperscript{*}
\\
\citet{Tassel2023} & Learn dispatching strategy &CP & Job-shop scheduling \\
\citet{Li2024TS} &VFA look-ahead  &Branch and Cut & DOP-rd\\
\citet{Liu2025} & Generate solver hints & GUROBI & Fleet relocation (0--1 IP)\\
\citet{Harsha2025} & MILP-optimized actor & CPLEX (B\&B) & Inventory replenishment\\
\citet{Chohlas2025} & Bandit + LP policy &Open-source solver & \begin{tabular}[c]{@{}l@{}} Transportation assistance allocation \end{tabular}\\
\bottomrule
\end{tabular}%
}
\\[2pt]
\begin{minipage}{\textwidth}
\footnotesize\raggedright
\textit{Notes:} [a] DDs = Decision Diagrams, MIS = Maximum Independent Set, IP = Integer Programming, ADMM = Alternating Direction Method of Multipliers, QP = Quadratic Programming, CP = Constraint Programming, VFA = Value Function
Approximation, DOP-rd = Orienteering Problem with Stochastic and Dynamic Release Dates, B\&B = Branch and Bound, LP = Linear Programming. TSP = Traveling Salesman Problem; SVM = Support Vector Machine.
\par
[b] * marks problems or application domains outside transportation, logistics, and closely related operations problems; the corresponding studies are included for methodological illustration.
\end{minipage}
\end{table}

A wide range of exact algorithms have been developed in the OR field, including the Alternating Direction Method of Multipliers (ADMM), Branch and Bound (B\&B), Branch and Cut, Branch and Price, Benders Decomposition, Constrained (or Relaxed) Decision Diagrams (DDs), and Constraint Programming (CP) to tackle static combinatorial optimization problems. Commercial solvers such as GUROBI and CPLEX include implementations of many of these methods. Recent studies have explored integrating RL with exact algorithms or commercial solvers to improve computational efficiency and solution quality. Unlike traditional RL approaches that directly generate solutions for COPs, this research stream focuses on using RL to steer, tune, or augment exact algorithmic frameworks. We classify the relevant literature based on the type of exact algorithm integrated with RL and summarize representative studies in Table~\ref{tab_RLExact}.

\textbf{(i) RL integrated with commercial solvers.} Several studies combine RL with commercial solvers, where the solver handles subproblems exactly while RL focuses on high-level policy decisions \citep[e.g.,][]{Liu2025, Harsha2025, Chohlas2025}. These studies show how exact solvers can be embedded in RL policies for dynamic problems. \citet{Liu2025} developed a multi-agent RL framework for the fleet relocation problem of shared autonomous electric vehicles, where agents in the RL algorithm generate relocation preferences fed into a 0-1 integer programming model solved by GUROBI. This approach tightly integrates RL-driven decision guidance with exact optimization. \citet{Harsha2025} introduced a solution framework that combines deep policy iteration with mathematical programming, using neural networks to approximate the value function while exact optimization determines the actions.

\textbf{(ii) RL with Decision Diagram-based exact algorithms.} DDs, such as relaxed or restricted DDs, are graph-based structures that compactly encode feasible regions or objective bounds. The quality of DD-based bounds is highly sensitive to the variable ordering. \citet{Cappart2019} formulated DD construction as a sequential decision process and used Q-learning to learn variable orderings that produce tighter bounds for Max-Cut and Maximum Independent Set problems. Here, the sequential decision process describes the algorithmic order in which a DD is constructed for a COP with fixed inputs, not the evolution of an operational system over time. The agent in the Q-learning algorithm selects the next variable based on the current DD state, with rewards linked to bound improvement. Building on this, \citet{Cappart2022} applied deep RL to learn effective variable orderings that improve both primal and dual bounds. This framework was further integrated into a full-fledged branch-and-bound algorithm, and the results showed that optimization bounds can be significantly enhanced through the use of deep RL algorithm.

\textbf{(iii) RL with the Alternating Direction Method of Multipliers.} ADMM is widely used for convex optimization problems such as quadratic programming (QP). \citet{Ichnowski2021} introduced RLQP, which uses deep RL to compute a policy that adapts the internal parameters of a QP solver to speed up the computation. Their framework significantly accelerates convergence compared to traditional fixed-parameter baselines. Here, RL acts as a meta-controller while the solver handles exact computations.

\textbf{(iv) RL with Branch and Cut and cutting-plane algorithms.} Branch and Cut methods integrate branch-and-bound search with cutting-plane generation, where selecting effective cuts (e.g., Gomory cuts) is key to performance. RL has been used to select cuts within this framework.  \citet{Tang2020} proposed a deep RL policy that selects Gomory cuts during cutting-plane iterations. Embedding this policy into Branch and Cut resulted in fewer nodes and faster convergence across several integer programming instances. The RL component is trained using an attention-based representation and evolution strategies. \citet{Li2024TS} combine RL with Branch and Cut differently: a value function approximation learned by RL is used with a look-ahead model solved by Branch and Cut for the orienteering problem with stochastic and dynamic release dates.

\textbf{(v) RL with Constraint Programming algorithms.} In this research stream, RL is typically used to learn the branching strategy in CP \citep[see,][]{Cappart2021, Tassel2023}. For instance, \citet{Cappart2021} formulated COPs in a unified dynamic programming-based framework and trained policies with Deep Q-Network or Proximal Policy Optimization methods to guide branching decisions, which were then embedded into a CP solver.

\subsubsection{Discussion \& future research directions}

The literature summarized in Table~\ref{tab_RLExact} indicates that RL–exact algorithm integration is particularly effective for large-scale problems where decisions such as branching, variable ordering, or cut selection strongly influence convergence but are difficult to manage through fixed, analytically derived rules. The reported improvements concern different measures, including bound quality and search efficiency. Three settings are especially well suited to RL integration: (i) when search, cutting, or sequencing decisions play a central role in the computational efficiency of the exact algorithms, such as branching in B\&B, cut selection in Branch and Cut, or variable ordering in DDs; (ii) when parameter tuning in exact algorithms is challenging, for example, the update rules in ADMM that affect convergence speed, RL can be used to learn adaptive tuning strategies that enhance algorithmic performance; (iii) when hierarchical frameworks are adopted, where the high-level RL component generates decisions, guidance, or priorities, and the low-level optimization module such as MILP, CP, or QP performs the final optimization. Such structures are common in transportation planning, scheduling, and fleet management.

Within these integrated frameworks, RL typically plays one of the following four roles. First, RL generates high-level decisions (such as preferences or priorities), that are then passed to commercial solvers for refinement. Second, RL provides value functions or policies to guide look-ahead decisions, which are subsequently optimized by exact algorithms such as Branch-and-Cut, CP, or QP. Third, RL learns algorithmic strategies such as branching, pruning, or cut selection. For example, guiding CP branching or Gomory cut selection to improve convergence behavior. Fourth, RL tunes parameters of exact algorithms, such as dynamically adjusting ADMM parameters to reduce the number of iterations and computation time.

Despite promising results, several limitations remain. First, learned branching, cutting, or sequencing strategies often generalize poorly. Their effectiveness tends to diminish or even become counterproductive when applied to problems with different sizes, structures, or distributions. Second, commercial solvers already include sophisticated heuristic and presolve routines for branching, ordering, and cut selection strategies. Achieving meaningful improvements using RL typically requires substantial training effort. Third, explainability remains a considerable challenge, as the rationale behind learned strategies is difficult to interpret and validate. Fourth, the problems tested so far are relatively simple. More realistic and operationally relevant scenarios, such as multi-hub systems, multiple vehicles, batch arrivals, partial charging, or uncertainty in service and charging times, have not yet been thoroughly addressed. Future research should therefore aim to enhance both the generalizability and explainability of RL-enhanced exact algorithms, particularly in frameworks that leverage the complementary strengths of learning and optimization. In parallel, more complex and practical problem settings should be incorporated to test the robustness and scalability of these approaches, such as the applications in the multi-hub location problem, partial charging problem, and integrated optimization problem of routing and charging. Such tests should follow the evaluation criteria in Section \ref{sec_benchmarking}, including the reporting of bounds and of instances solved to proven optimality.

\subsection{Benchmarking and reproducibility}
\label{sec_benchmarking}

In this section, we present common evaluation criteria and specify additional measures for exact methods, including bounds and the number of instances solved to proven optimality. Table \ref{tab_rl_or_comparison} summarizes where the reviewed studies report improvements over OR methods and where the evidence does not establish an advantage. We selected studies that compare RL-based methods with OR methods under stated test conditions. We distinguish solution quality, solution time after training, and improvements in intermediate measures such as bounds. Because the studies differ in their objectives, baselines, and computational budgets, we do not report a pooled percentage improvement. The reported advantages depend on the problem, the outcome measured, the baseline, and the computational budget. Three limits recur. First, the strength of the baseline matters. An advantage over a general-purpose tool such as OR-Tools does not establish an advantage over the strongest problem-specific heuristics, and several comparisons use simple rules or the original algorithm without RL rather than a well-tuned solver.  Second, improvements in intermediate measures, such as tighter bounds or fewer cuts, do not by themselves establish a shorter time to proven optimality. Third, conclusions depend on the computational budget.  In \cite{Kallestad2023}, the advantage over the original ALNS operator selection narrows as the number of iterations increases, and the reported solution times in the routing examples exclude training. These limits motivate the criteria on baselines, performance measures, and computational cost below.

\begin{table}[!htbp]
\centering
\caption{Reported advantages and limits of RL relative to OR methods
in representative studies.}
\label{tab_rl_or_comparison}
\small
\setlength{\tabcolsep}{6pt}
\renewcommand{\arraystretch}{1.15}
\renewcommand{\cellgape}{\Gape[4pt][4pt]}

\resizebox{\textwidth}{!}{%
\begin{tabular}{@{} c|c c @{}}
\hline

\makecell[cr]{\textbf{Use of RL}}
&
\makecell[cl]{\textbf{Reported advantage and comparison}}
&
\makecell[cl]{\textbf{Limit or case without an advantage}}
\\
\hline

\makecell[cc]{
Independent RL\\
solver (Section \ref{sec_RLsolver})
}
&
\makecell[cl]{
\citet{Nazari2018} report better solution\\
quality than classical heuristics and\\
OR-Tools on medium-sized capacitated\\
vehicle routing instances, with comparable\\
solution time after training.
}
&
\makecell[cl]{
An advantage over OR-Tools does not\\
imply an advantage over specialist\\
heuristics on solution quality. On 100-node\\
CVRP instances, the greedy method of\\
\citet{Kool2018} has a smaller gap than\\
OR-Tools (7.34\% vs.\ 9.67\%) but remains\\
7.34\% above LKH3. On 100-node TSP\\
instances, its greedy gap is 4.53\%, versus\\
0\% for Concorde. Sampling reduces the\\
gap to 2.26\%, at a higher computational\\ cost than greedy decoding.
}
\\
\hline

\makecell[cc]{
RL with heuristic\\
algorithms (Section \ref{sec_RLheuristic})
}
&
\makecell[cl]{
\citet{Kallestad2023} report better\\
objective values than the original ALNS\\
operator selection on most tested variants,\\
particularly larger instances.\\
\citet{An2026} report a 4.88\% average\\
gain in computational efficiency from\\
adding RL to VNS.
}
&
\makecell[cl]{
The advantage reported by\\
\citet{Kallestad2023} narrows as the\\
number of search iterations increases, and\\
elapsed time is not uniformly lower. The\\
separate 6.51\% improvement over other\\
heuristics reported by \citet{An2026}\\
is for their complete method, so it should\\
not be attributed solely to RL.
}
\\
\hline

\makecell[cc]{
RL with exact\\
algorithms (Section \ref{sec_RLexact})
}
&
\makecell[cl]{
\citet{Cappart2022} report improved\\
decision-diagram bounds from learned\\
variable ordering on maximum independent\\
set instances. \citet{Tang2020} report\\
fewer cuts and greater integrality-gap\\
closure than standard cut-selection rules.
}
&
\makecell[cl]{
Better bounds are an intermediate result. A\\
comparison of complete methods should\\
also consider solution time, the cost of\\
learned ordering, and instances solved to\\
proven optimality. \citet{Tang2020}\\
compare with simple selection rules in a\\
simplified branch-and-cut, not with the\\
cut management of a commercial solver.
}
\\
\hline

\end{tabular}%
}

\par\smallskip
{\footnotesize\raggedright
\textit{Note:} These percentages measure different outcomes under
different test conditions and should not be interpreted as a typical
improvement from RL.
\par}

\end{table}

Building on this comparison, we outline evaluation criteria and resources for reproducible benchmarking. Public tools and benchmark instances provide a practical starting point for reproducible evaluation. \href{https://github.com/hubbs5/or-gym}{OR-Gym} \citep{Hubbs2020ORGymSoftware} provides environments for inventory management, supply chains, and vehicle routing. \href{https://github.com/ai4co/rl4co}{RL4CO} \citep{Berto2025RL4CO} provides combinatorial optimization environments, RL implementations, and experiment configurations, including TSP and CVRP examples. For standard routing instances, \href{https://comopt.ifi.uni-heidelberg.de/software/TSPLIB95/}{TSPLIB} \citep{Reinelt1991TSPLIB} provides TSP instances and reference solutions. \href{https://galgos.inf.puc-rio.br/cvrplib/index.php/en/}{CVRPLIB} includes CVRP instance sets, such as those of \citet{Uchoa2017CVRPLIB}, together with best known solutions and their optimality status. Strong OR baselines are also available as open-source code. For the CVRP, for example, the hybrid genetic search of \citet{Vidal2022} is a state-of-the-art heuristic against which learning-based methods can be compared. \citet{Accorsi2022} give guidelines for the computational testing of machine learning approaches to vehicle routing problems.

Using a common library alone does not ensure comparable results. For the methods discussed in Sections~\ref{sec_RLsolver}, \ref{sec_RLheuristic}, and~\ref{sec_RLexact}, we suggest the following evaluation criteria:
\begin{enumerate}[leftmargin=*, itemsep=0.4ex, topsep=0.6ex]
    \item \textbf{Test instances.} Use the same test instances, objectives, and constraints for all methods. Keep training, validation, and test sets separate. State the instance sizes and generation distributions, and test generalization to sizes and distributions not used during training.
    \item \textbf{Baselines.} Compare with strong OR and RL baselines and document their parameter settings and tuning budgets. For methods that combine RL with a heuristic or exact algorithm, include the original algorithm without RL to assess the contribution of the learned component.
    \item \textbf{Performance measures.} Report objective values, feasibility rates, and gaps, stating whether each reference value is a proven optimum, a bound, or a best known solution. For exact methods, also report bounds and the number of instances solved to proven optimality.
    \item \textbf{Computational cost.} Report training and parameter tuning costs separately from solution time. Include decoding, local search, and solver calls in solution time. Use common limits on solution time and state the CPU, GPU, and parallel processing settings.
    \item \textbf{Repeated experiments.} Repeat training and stochastic solution procedures with multiple independent random seeds. Report the number of runs and uncertainty, such as confidence intervals, alongside average performance.
\end{enumerate}
To support reproducibility, authors should release code, trained parameters, instance lists or generators, configuration files, software versions, random seeds, and evaluation scripts. A release or commit identifier should be recorded so that the reported experiments can be reproduced after the repository changes.

\section{\texorpdfstring{{Digital Twins for RL Deployment and Evaluation in OR}}{Digital Twins for RL Deployment and Evaluation in OR}}
\label{Sec_digitalTwins}

Digital twins combine digital models of physical systems with operating data to support monitoring, simulation, and operational decisions. Within our taxonomy, DTs provide an optional deployment and evaluation layer that can support either methodological role of RL in OR. They can provide environments for policy training and testing, supply information for operational decisions, or introduce decisions about the resources and data needed to maintain the DT. This layer requires a separate discussion because model calibration, data updates, and the transfer of policies to physical systems raise issues that are not captured by the choice of RL algorithm alone. Section~\ref{sec_analysisDT} compares these uses, and Section~\ref{sec_challengesDT} discusses their conditions of use, limitations, and research priorities.

\subsection{Literature analysis}
\label{sec_analysisDT}

Table~\ref{tab_digitalTwins} groups ten studies by the main role of DTs and RL. The three streams are DTs as training environments for RL, DTs as accelerators for RL, and RL for operational control with DTs. These uses can overlap. In particular, a DT used to accelerate learning also provides a training environment. The distinction concerns whether the study primarily uses the DT to train an operational policy or specifically develops a mechanism to improve learning efficiency.

\begin{table}[h]
\centering
\caption{Classification of papers studying the integration of RL and DTs by their main role.}
\label{tab_digitalTwins}
\resizebox{\textwidth}{!}{%
\begin{tabular}{r|lll}
\hline
Publications & Domain & Roles of DT and RL & RL method \\
\hline
\citet{Liu2023}
& \begin{tabular}[c]{@{}l@{}} IoT networks: virtual network \\ function migration\textsuperscript{*} \end{tabular}    &   DTs as training environments for RL   &  DPPO\\
\citet{Kamal2024}
& Traffic signal control   &   DTs as training environments for RL   &  MADDPG\\
\citet{Schlappa2024}
& Control of waste incineration plants\textsuperscript{*}   &   DTs as training environments for RL   &  Deep Q-network\\
\citet{Tang2023}
& Multi-UAV task assignment   &   DTs as accelerators for RL   &  Deep Q-learning\\
\citet{Wu2021}
&  Autonomous driving  &   DTs as accelerators for RL   &  Actor-critic\\
\citet{Zhang2024}
&  \begin{tabular}[c]{@{}l@{}} Network slicing: \\ resource management\textsuperscript{*} \end{tabular}  &   DTs as accelerators for RL   & Deep
Q-network \\
\citet{Yan2022}  & Scheduling in manufacturing systems       &   RL for operational control with DTs   &  Double-layer Q-learning\\
\citet{Park2022}
& Production control    &   RL for operational control with DTs    & Dueling Q-network \\
\citet{Wang2024}
&  \begin{tabular}[c]{@{}l@{}} Vehicle platoon DT networks: \\ resource allocation \end{tabular}  &   RL for operational control with DTs    &  MADDPG\\
\citet{Xu2024}
&  \begin{tabular}[c]{@{}l@{}} Internet of Energy: \\ DT data mapping\textsuperscript{*} \end{tabular}  &   RL for operational control with DTs    &  PPO \\
\hline
\end{tabular}%
}
\\[2pt]
\begin{minipage}{\textwidth}
\footnotesize\raggedright
\textit{Notes:} [a] MADDPG = Multi-agent Deep Deterministic Policy Gradient, PPO = Proximal Policy Optimization, DPPO = Distributed Proximal Policy Optimization, UAV = Unmanned Aerial Vehicle, IoT = Internet of Things. Roles can overlap.
\par
[b] * marks problems or application domains outside transportation, logistics, and closely related operations problems; the corresponding studies are included for methodological illustration.
\end{minipage}
\end{table}

\textbf{(i) DTs as training environments for RL.} In this stream, RL learns operational policies through interaction with a digital representation of the system. \citet{Kamal2024} use a DT to simulate traffic scenarios and evaluate traffic signal control with MADDPG, using synthetic data and traffic data from Amman. \citet{Schlappa2024} construct a neural network model from waste incineration plant data and use it as an environment for training a deep Q-network. In network management, \citet{Liu2023} construct a virtual network for training migration policies with distributed PPO. Resource demand predictions provide inputs for these decisions, so the contribution of the prediction component should be distinguished from that of the RL policy. These applications illustrate how a DT can support policy development when experiments on the physical system are costly or difficult. However, using operating data to construct or test a DT does not by itself establish that the resulting policy performs well in the physical system.

\textbf{(ii) DTs as accelerators for RL.} These studies use the DT to improve learning efficiency through prediction or pretraining. For autonomous driving, \citet{Wu2021} use a DT to predict transition dynamics and generate experience for RL, with evaluation in a simulated lateral control task. For network slicing, \citet{Zhang2024} construct a DT from historical data and use predicted states and rewards to train a second RL agent. Parameters learned by this agent are used to update the agent interacting with the network environment. The framework also includes knowledge distillation and an extension to offline RL. Their numerical experiments examine learning efficiency and resource allocation performance. For task reassignment in a system with multiple unmanned aerial vehicles, \citet{Tang2023} use DT pretraining to provide Q values for deep Q-learning and reduce the required training episodes. Thus, acceleration can result from predictive models or from transferring information obtained during pretraining. Fewer interactions or training episodes should be assessed alongside the cost of constructing, calibrating, and running the DT.

\textbf{(iii) RL for operational control with DTs.} In this stream, the object of the RL decision differs across studies. \citet{Yan2022} use a DT to detect disturbances and trigger rescheduling in manufacturing, with double-layer Q-learning selecting machines and operations. \citet{Park2022} integrate DTs and RL for production control in a re-entrant job shop. Both address production decisions supported by a digital representation of the system. In contrast, \citet{Wang2024} allocate bandwidth and computational resources needed to maintain platoon DTs through MADDPG. \citet{Xu2024} use deep RL to select mapping data according to system states and user requirements, reducing the amount of mapped data and resource consumption. These latter studies concern the operation of the DT itself. Their results should therefore be assessed using resource and mapping measures, together with the effect on the operational decisions that the DT supports.

The comparison reveals two important distinctions. First, training and acceleration studies address how a policy is learned, whereas operational control studies address decisions made with or about the DT. These functions can be combined within one system, but they should be evaluated separately to identify the contribution of each component. Second, improved performance in a digital environment and improved performance in a physical system are different forms of evidence. For example, the evaluations in \citet{Wu2021}, \citet{Zhang2024}, and \citet{Kamal2024} use simulation. Such results support conclusions within the tested environments. Reliable performance in the physical system requires further validation, including tests under operating conditions not used to construct or train the system. The studies also identify specific limitations. \citet{Liu2023} acknowledge that their theoretical DT network does not account for differences from the physical network. In \citet{Schlappa2024}, RL performance is evaluated in the DT while human operator performance comes from the physical plant, which limits the comparability of the reported rewards. These limitations show why the evaluation setting is central to interpreting the contribution of a DT.

\Needspace{4\baselineskip}
\subsection{Challenges and future research directions}
\label{sec_challengesDT}

The studies in Table~\ref{tab_digitalTwins} demonstrate several uses of DTs, but their objectives and evaluation settings differ. We therefore discuss when the integration is useful before examining the costs and limitations associated with each use.

\Needspace{4\baselineskip}
\subsubsection{When to use RL with DTs}

DTs are particularly relevant for RL training when repeated experimentation on the physical system is risky, costly, or slow, as in plant control and autonomous driving. They are also useful when a policy must be tested under demand changes, equipment failures, or other disturbances that are difficult to observe frequently in operating data. These benefits depend on whether the DT reproduces the system responses that matter for the decision. A visually detailed representation alone is insufficient if its transition dynamics, constraints, or operational costs are inaccurate.

For operational control, a DT is useful when updated information changes the decisions that should be made. The disturbance detection and rescheduling mechanism in \citet{Yan2022} illustrates this connection. RL can also address repeated decisions about limited communication and computational resources, as in \citet{Wang2024}, or data selection, as in \citet{Xu2024}. In these cases, the formulation should identify the decisions available to the agent and explain how resource use or mapping quality affects the operational objective.

The additional complexity of RL with DTs may be difficult to justify when established OR methods already solve the decision problem reliably within the available time, when operating conditions change little, or when the DT cannot be calibrated with sufficient accuracy. Maintaining the DT also requires data collection, model updates, and computational resources. These costs should be compared with the expected improvement in decision quality. Conventional simulation optimization, model predictive control, and suitable heuristics provide relevant alternatives, depending on the problem.

\Needspace{4\baselineskip}
\subsubsection{Technical challenges}

The first challenge is the relationship between model accuracy and policy performance. A DT may predict average system behavior well while misrepresenting rare events or states visited by an RL policy. Increasing model detail can improve some predictions but also increase data requirements, calibration effort, and simulation time. The relevant question is whether greater accuracy improves feasibility and operational performance enough to justify these costs. The plant model of \citet{Schlappa2024} and the predictive DTs of \citet{Wu2021} and \citet{Zhang2024} make this dependence on model quality particularly important.

A second challenge concerns the cost of learning and maintaining the DT. Reducing interactions with a physical system or reducing training episodes does not necessarily reduce total computational cost. Model construction, calibration, simulated experience, and policy training all contribute to that cost. During operation, frequent data updates can provide more current observations but require communication and computation. The resource allocation problem studied by \citet{Wang2024} and the data selection problem studied by \citet{Xu2024} show why these requirements cannot be treated as free inputs. Delayed or missing observations can also affect the decisions made by the policy.

A third challenge is validation beyond the conditions used to construct and train the system. Changes in demand, equipment, network topology, or user behavior can affect both the DT and the policy. Simulation studies can isolate such effects, but evaluation within the same model used for training can conceal errors shared by the model and policy. Tests using operating data from separate periods, alternative simulation assumptions, and, where available, physical system observations are therefore needed. Results should distinguish model prediction accuracy, policy performance in simulation, and performance observed during operation.

\Needspace{4\baselineskip}
\subsubsection{Future research directions}

Future work should relate DT calibration to the quality of the decisions made by RL. Alongside prediction error, studies should report operational costs or rewards, service levels, and constraint violations. Data used to construct and calibrate the DT should be separated from data used for evaluation, for example by time or scenario. Additional tests should cover changes in demand, failures, and delayed observations. Following the evaluation criteria in Section~\ref{sec_benchmarking}, comparisons should include strong OR methods and the same RL method using a DT that is not updated with new operating data. This would help identify whether maintaining an updated DT provides benefits beyond using a fixed simulation model.

Another direction is to study how model detail, data update frequency, and resource allocation affect operational performance together. Existing resource allocation and data mapping studies provide starting points, but joint control of these decisions should be treated as a research direction rather than an established capability of all DT applications. Comparisons should include fixed update schedules and simple resource allocation rules, and report the resulting delay, data volume, computational cost, and decision quality.

Finally, reproducible validation requires reporting the origin and coverage of operating data, the DT construction and calibration procedures, update rules, and the costs of model maintenance and policy training. Shared environments should include operational constraints and documented disturbances. Where suitable historical data are available, offline policy evaluation can help assess a policy before physical deployment, with its assumptions and uncertainty stated explicitly. Reward definitions and operational performance measures should be reported separately to support meaningful comparisons. Such evidence would make it easier to determine when the integration of RL and DTs improves on established OR approaches.

\section{Synthesis and conclusion}
\label{Sec_conclusion}

This section brings together the findings of the review. Section~\ref{sec:guidance} consolidates the guidance on when each RL approach is suitable, and Section~\ref{sec:concluding} summarizes the main conclusions, the remaining challenges, and the priorities for future research.

\subsection{Guidance for selecting RL approaches}
\label{sec:guidance}

Table~\ref{tab:guidance} consolidates the guidance developed in Sections~\ref{Sec_MDP} to~\ref{Sec_digitalTwins}. It follows the taxonomy in Figure~\ref{fig_Taxonomy}: the problem class determines the first choice, and the characteristics listed in the table inform the second. Within each problem class, the rows are not mutually exclusive. For example, model-based components can be combined with actor-critic methods, and RL-generated initial solutions can be refined by a heuristic. Observability and the time horizon should also be considered in each application, but the reviewed studies do not establish clear rules for selecting an RL algorithm family from either characteristic alone.

The table summarizes considerations drawn from the reviewed studies rather than prescriptive selection rules. As Table~\ref{tab_rl_or_comparison} shows, reported advantages depend on the problem, the baseline, and the computational budget. Any proposed RL approach should therefore be compared with a strong OR baseline using the evaluation criteria in Section~\ref{sec_benchmarking}, including the costs of training and solving instances. Where an established OR method already meets the operational requirements, the additional cost of learning and maintaining an RL component needs a clear justification.

\begin{table}[htbp]
\caption{Guidance for selecting RL approaches by problem characteristics.}
\label{tab:guidance}
\footnotesize
\setlength{\tabcolsep}{4pt}
\renewcommand{\arraystretch}{1.25}
\begin{tabular}{>{\raggedright\arraybackslash}p{0.14\linewidth}
                >{\raggedright\arraybackslash}p{0.35\linewidth}
                >{\raggedright\arraybackslash}p{0.34\linewidth}
                >{\raggedright\arraybackslash}p{0.08\linewidth}}
\hline
RL approach & Favorable conditions & Main limitations & Section \\
\hline
\multicolumn{4}{@{}l}{\textit{Dynamic sequential decision problems}} \\
Value-based RL
  & Discrete actions; interaction data or a simulator available
  & Training may require substantial data; standard value-based methods are less direct for continuous or mixed actions
  & \ref{sec_value} \\
Policy-based RL
  & Continuous or mixed actions; online or rolling decisions
  & Training may require many interactions, usually in a simulator; operational constraints need explicit treatment
  & \ref{sec_policy} \\
Actor-critic RL
  & Continuous or high-dimensional actions; rolling control or multi-agent coordination
  & Training can be sensitive to design choices; multi-agent methods may add coordination and communication costs
  & \ref{sec_actor} \\
Model-based RL
  & Real-system interaction is costly or risky; dynamics are known or can be learned; sample efficiency matters
  & Model errors can compound over long horizons; model development and, where used, online planning add costs
  & \ref{sec_model} \\[2pt]
\multicolumn{4}{@{}l}{\textit{Static COPs}} \\
Independent RL solver
  & Many instances of similar size and distribution must be solved; short solution time after training is valuable
  & Solution quality can trail the strongest heuristics and exact solvers; feasibility and transfer to other sizes and distributions require testing
  & \ref{sec_RLsolver} \\
RL with heuristic algorithms
  & An effective heuristic exists and operator, neighborhood, or parameter choices affect its performance
  & No general optimality guarantee; training costs may outweigh gains; learned choices may be problem-specific
  & \ref{sec_RLheuristic} \\
RL with exact algorithms
  & Branching, ordering, cut selection, or tuning affects solver performance, and RL can guide these choices while preserving validity
  & Guarantees hold only if validity is preserved; gains are mostly reported against simple rules rather than tuned solvers; learned strategies may not transfer across sizes and structures
  & \ref{sec_RLexact} \\[2pt]
\multicolumn{4}{@{}l}{\textit{Deployment and evaluation}} \\
Digital twin layer
  & Physical experimentation is costly or risky; a sufficiently representative model can support training or testing
  & Calibration and maintenance require resources; performance in simulation does not establish physical performance
  & \ref{Sec_digitalTwins} \\
\hline
\end{tabular}
\end{table}

\subsection{Concluding remarks}
\label{sec:concluding}
This paper has presented a taxonomy of how RL contributes to OR and reviewed how RL methods are integrated with OR methods. We distinguish two methodological roles. RL learns policies for repeated decisions in dynamic systems, and it constructs solutions or supports heuristic and exact algorithms for static COPs. Digital twins form an optional deployment and evaluation layer that can support both roles through policy training, testing, validation, and operational control. This distinction connects the review of RL methods with the conditions required for their use in physical systems.

Section~\ref{sec:guidance} summarizes conditions under which the reviewed approaches may be considered (Table~\ref{tab:guidance}). For dynamic sequential decision problems, the reviewed applications include high-dimensional or uncertain environments and settings requiring frequent or real-time decisions. The action space and the availability of data or a simulator are among the considerations in choosing an algorithm family. For static COPs, reported gains can arise when RL improves the algorithmic solution process by learning branching, cutting, variable-ordering, operator-selection, neighborhood-selection, or parameter-tuning strategies, although these gains depend on the baseline and the computational budget (Table~\ref{tab_rl_or_comparison}). When digital twins are available, they can further facilitate safe experimentation, scenario generation, and pre-deployment evaluation, although effective deployment depends critically on rigorous calibration to reduce gaps between simulated and physical systems.

Overall, RL methods and established OR methods are complementary. Their integration provides a promising pathway to more efficient and higher-quality decision-making in complex systems and opens broad opportunities for rigorous research and impactful applications. Nevertheless, several important challenges remain. RL's low sample efficiency and limited training stability restrict its broad applicability. Generalization and transfer across problem sizes, topologies, and distributions remain weak. Explanations and theoretical guarantees for learned strategies are often lacking. Benchmarks and evaluation protocols remain fragmented. Digital twins require reliable data and careful model validation, while multi-agent implementations introduce latency and bandwidth constraints.

To address these challenges and advance the field, we outline several priorities. First, establishing standardized benchmarks and unified evaluation protocols is essential for fair comparison and reproducibility. These should build on the criteria in Section~\ref{sec_benchmarking} and include comparisons with strong OR baselines and, where available, bounds on the optimal value. Second, future research should develop hybrid frameworks that combine model-free and model-based RL algorithms to solve dynamic sequential decision problems with high-quality solutions and strong computational efficiency. Third, for static COPs, efforts should focus on strengthening the generalizability and explainability of RL-enhanced heuristic and exact algorithms, with the long-term goal of building unified frameworks that coordinate learning and optimization. Fourth, for deployment with digital twins, reporting model calibration quality, comparing the full costs of training and maintenance, and testing policies under new scenarios and disturbances would strengthen the evidence for practical use.

\bibliographystyle{cas-model2-names}

\bibliography{liter.bib}

@article{Mazyavkina2021,
  author  = {Mazyavkina, Nina and Sviridov, Sergey and Ivanov, Sergei and Burnaev, Evgeny},
  title   = {Reinforcement learning for combinatorial optimization: A survey},
  journal = {Computers \& Operations Research},
  volume  = {134},
  pages   = {105400},
  year    = {2021},
  doi = {10.1016/j.cor.2021.105400}
}

@article{Scavuzzo2026,
  author  = {Scavuzzo, Lara and Aardal, Karen and Lodi, Andrea and Yorke-Smith, Neil},
  title   = {Machine learning augmented branch and bound for mixed integer linear programming},
  journal = {Mathematical Programming},
  volume  = {217},
  pages   = {123--166},
  year    = {2026},
  doi = {10.1007/s10107-024-02130-y}
}

@article{Boute2022,
  author  = {Boute, Robert N. and Gijsbrechts, Joren and {Van Jaarsveld}, Willem and Vanvuchelen, Nathalie},
  title   = {Deep reinforcement learning for inventory control: A roadmap},
  journal = {European Journal of Operational Research},
  volume  = {298},
  number  = {2},
  pages   = {401--412},
  year    = {2022},
  doi = {10.1016/j.ejor.2021.07.016}
}

@article{Hildebrandt2023,
  author  = {Hildebrandt, Florentin D. and Thomas, Barrett W. and Ulmer, Marlin W.},
  title   = {Opportunities for reinforcement learning in stochastic dynamic vehicle routing},
  journal = {Computers \& Operations Research},
  volume  = {150},
  pages   = {106071},
  year    = {2023},
  doi = {10.1016/j.cor.2022.106071}
}

@article{Yan2022review,
  author  = {Yan, Yimo and Chow, Andy H. F. and Ho, Chin Pang and Kuo, Yong-Hong and Wu, Qihao and Ying, Chengshuo},
  title   = {Reinforcement learning for logistics and supply chain management: Methodologies, state of the art, and future opportunities},
  journal = {Transportation Research Part E: Logistics and Transportation Review},
  volume  = {162},
  pages   = {102712},
  year    = {2022},
  doi = {10.1016/j.tre.2022.102712}
}

@article{Farazi2021,
  author  = {Farazi, Nahid Parvez and Zou, Bo and Ahamed, Tanvir and Barua, Limon},
  title   = {Deep reinforcement learning in transportation research: A review},
  journal = {Transportation Research Interdisciplinary Perspectives},
  volume  = {11},
  pages   = {100425},
  year    = {2021},
  doi     = {10.1016/j.trip.2021.100425}
}

@ARTICLE{Kamal2024,
  author={Kamal, Hani and Yánez, Wendy and Hassan, Sara and Sobhy, Dalia},
  journal={IEEE Internet of Things Journal}, 
  title={Digital-twin-based deep reinforcement learning approach for adaptive traffic signal control}, 
  year={2024},
  volume={11},
  number={12},
  pages={21946-21953},
  doi={10.1109/JIOT.2024.3377600}
}

@article{Yan2022,
title = {Digital twin-enabled dynamic scheduling with preventive maintenance using a double-layer Q-learning algorithm},
journal = {Computers \& Operations Research},
volume = {144},
pages = {105823},
year = {2022},
issn = {0305-0548},
doi = {10.1016/j.cor.2022.105823},
author = {Qi Yan and Hongfeng Wang and Fang Wu}
}

@article{Pan2021,
  author  = {Pan, Tianlu and Guo, Renzhong and Lam, William H. K. and Zhong, Renxin and Wang, Weixi and He, Biao},
  title   = {Integrated optimal control strategies for freeway traffic mixed with connected automated vehicles: A model-based reinforcement learning approach},
  journal = {Transportation Research Part C: Emerging Technologies},
  volume  = {123},
  pages   = {102987},
  year    = {2021},
  doi     = {10.1016/j.trc.2021.102987}
}

@article{Barrett2020, title={Exploratory combinatorial optimization with reinforcement learning}, volume={34},DOI={10.1609/aaai.v34i04.5723}, number={04}, journal={Proceedings of the AAAI Conference on Artificial Intelligence}, author={Barrett, Thomas and Clements, William and Foerster, Jakob and Lvovsky, Alex}, year={2020}, pages={3243-3250} }

@article{Jin2023, title={Pointerformer: Deep reinforced multi-pointer transformer for the traveling salesman problem}, volume={37}, DOI={10.1609/aaai.v37i7.25982},
number={7}, journal={Proceedings of the AAAI Conference on Artificial Intelligence}, author={Jin, Yan and Ding, Yuandong and Pan, Xuanhao and He, Kun and Zhao, Li and Qin, Tao and Song, Lei and Bian, Jiang}, year={2023},  pages={8132-8140} }

@inproceedings{Nazari2018,
  author    = {Nazari, MohammadReza and Oroojlooy, Afshin and Snyder, Lawrence V. and Tak{\'a}{\v{c}}, Martin},
  title     = {Reinforcement learning for solving the vehicle routing problem},
  booktitle = {Advances in Neural Information Processing Systems},
  editor    = {Bengio, S. and Wallach, H. and Larochelle, H. and Grauman, K. and Cesa-Bianchi, N. and Garnett, R.},
  volume    = {31},
  year      = {2018},
  publisher = {Curran Associates, Inc.},
  url       = {https://proceedings.neurips.cc/paper_files/paper/2018/file/9fb4651c05b2ed70fba5afe0b039a550-Paper.pdf}
}

@inproceedings{Dai2017,
  author    = {Dai, Hanjun and Khalil, Elias B. and Zhang, Yuyu and Dilkina, Bistra and Song, Le},
  title     = {Learning combinatorial optimization algorithms over graphs},
  booktitle = {Advances in Neural Information Processing Systems},
  editor    = {Guyon, I. and von Luxburg, U. and Bengio, S. and Wallach, H. and Fergus, R. and Vishwanathan, S. and Garnett, R.},
  volume    = {30},
  year      = {2017},
  publisher = {Curran Associates, Inc.},
  url       = {https://proceedings.neurips.cc/paper_files/paper/2017/file/d9896106ca98d3d05b8cbdf4fd8b13a1-Paper.pdf}
}

@inproceedings{Ma2021,
  author    = {Ma, Yi and Hao, Xiaotian and Hao, Jianye and Lu, Jiawen and Liu, Xing and Tong, Xialiang and Yuan, Mingxuan and Li, Zhigang and Tang, Jie and Meng, Zhaopeng},
  title     = {A hierarchical reinforcement learning based optimization framework for large-scale dynamic pickup and delivery problems},
  booktitle = {Advances in Neural Information Processing Systems},
  editor    = {Ranzato, M. and Beygelzimer, A. and Dauphin, Y. and Liang, P. and {Wortman Vaughan}, J.},
  volume    = {34},
  pages     = {23609--23620},
  year      = {2021},
  publisher = {Curran Associates, Inc.},
  url       = {https://proceedings.neurips.cc/paper_files/paper/2021/file/c6a01432c8138d46ba39957a8250e027-Paper.pdf}
}

@inproceedings{Ichnowski2021,
 author = {Ichnowski, Jeffrey and Jain, Paras and Stellato, Bartolomeo and Banjac, Goran and Luo, Michael and Borrelli, Francesco and Gonzalez, Joseph E and Stoica, Ion and Goldberg, Ken},
 booktitle = {Advances in Neural Information Processing Systems},
 editor = {M. Ranzato and A. Beygelzimer and Y. Dauphin and P.S. Liang and J. Wortman Vaughan},
 pages = {21043--21055},
 publisher = {Curran Associates, Inc.},
 title = {Accelerating quadratic optimization with reinforcement learning},
 url = {https://proceedings.neurips.cc/paper_files/paper/2021/file/afdec7005cc9f14302cd0474fd0f3c96-Paper.pdf},
 volume = {34},
 year = {2021}
}

@article{Meng2025,
title = {Dynamic volunteer assignment: Integrating skill diversity, task variability and volunteer preferences},
journal = {Transportation Research Part E: Logistics and Transportation Review},
volume = {197},
pages = {104068},
year = {2025},
issn = {1366-5545},
doi = {10.1016/j.tre.2025.104068},
author = {Qingchun Meng and Bo Feng and Guodong Yu},
}

@article{Tian2025,
title = {PTB: A deep reinforcement learning method for flexible logistics service combination problem with spatial-temporal constraint},
journal = {Transportation Research Part E: Logistics and Transportation Review},
volume = {195},
pages = {103978},
year = {2025},
issn = {1366-5545},
doi = {10.1016/j.tre.2025.103978},
author = {Ran Tian and Longlong Chang and Zhihui Sun and Guanglu Zhao and Xin Lu}
}

@article{Wang2025,
title = {Distributed multi-agent reinforcement learning approach for energy-saving optimization under disturbance conditions},
journal = {Transportation Research Part E: Logistics and Transportation Review},
volume = {200},
pages = {104180},
year = {2025},
issn = {1366-5545},
doi = {10.1016/j.tre.2025.104180},
author = {Dahan Wang and Jianjun Wu and Ximing Chang and Haodong Yin}
}

@article{Cappart2019, title={Improving optimization bounds using machine learning: Decision diagrams meet deep reinforcement learning}, volume={33}, DOI={10.1609/aaai.v33i01.33011443}, number={01}, journal={Proceedings of the AAAI Conference on Artificial Intelligence}, author={Cappart, Quentin and Goutierre, Emmanuel and Bergman, David and Rousseau, Louis-Martin}, year={2019}, pages={1443-1451} }

@article{Teusch2025,
  author  = {Teusch, Julian and {Neumann Saavedra}, Bruno and Scherr, Yannick Oskar and M{\"u}ller, J{\"o}rg P.},
  title   = {Strategic planning of geo-fenced micro-mobility facilities using reinforcement learning},
  journal = {Transportation Research Part E: Logistics and Transportation Review},
  volume  = {194},
  pages   = {103872},
  year    = {2025},
  doi     = {10.1016/j.tre.2024.103872}
}

@article{Dehaybe2024,
title = {Deep Reinforcement Learning for Inventory Optimization with Non-stationary Uncertain Demand},
journal = {European Journal of Operational Research},
volume = {314},
number = {2},
pages = {433-445},
year = {2024},
issn = {0377-2217},
doi = {10.1016/j.ejor.2023.10.007},
author = {Henri Dehaybe and Daniele Catanzaro and Philippe Chevalier}
}

@article{Teck2025,
title = {Deep reinforcement learning for the real-time inventory rack storage assignment and replenishment problem},
journal = {European Journal of Operational Research},
volume = {327},
number = {2},
pages = {606-622},
year = {2025},
issn = {0377-2217},
doi = {10.1016/j.ejor.2025.05.008},
author = {Teck, Sander and Ph{\d{a}}m, T{\'u} San and Rousseau, Louis-Martin and Vansteenwegen, Pieter}
}

@article{Zhang2025,
title = {Multi-objective cooperative co-evolution algorithm with hypervolume-based Q-learning for hybrid seru system},
journal = {European Journal of Operational Research},
volume = {324},
number = {3},
pages = {839-854},
year = {2025},
issn = {0377-2217},
doi = {10.1016/j.ejor.2025.02.025},
author = {Zhecong Zhang and Yang Yu and Xuqiang Qi and Yangguang Lu and Xiaolong Li and Ikou Kaku}
}

@article{Jin2024,
title = {Container port truck dispatching optimization using Real2Sim based deep reinforcement learning},
journal = {European Journal of Operational Research},
volume = {315},
number = {1},
pages = {161-175},
year = {2024},
issn = {0377-2217},
doi = {10.1016/j.ejor.2023.11.038},
author = {Jiahuan Jin and Tianxiang Cui and Ruibin Bai and Rong Qu}
}

@article{Vanvuchelen2024,
title = {Cluster-based lateral transshipments for the Zambian health supply chain},
journal = {European Journal of Operational Research},
volume = {313},
number = {1},
pages = {373-386},
year = {2024},
issn = {0377-2217},
doi = {10.1016/j.ejor.2023.08.005},
author = {Nathalie Vanvuchelen and Kim {De Boeck} and Robert N. Boute}
}

@article{Wu2024,
title = {An effective hybrid evolutionary algorithm for the clustered orienteering problem},
journal = {European Journal of Operational Research},
volume = {313},
number = {2},
pages = {418-434},
year = {2024},
issn = {0377-2217},
doi = {10.1016/j.ejor.2023.08.006},
author = {Qinghua Wu and Mu He and Jin-Kao Hao and Yongliang Lu}
}

@article{Li2024,
title = {A flow based formulation and a reinforcement learning based strategic oscillation for cross-dock door assignment},
journal = {European Journal of Operational Research},
volume = {312},
number = {2},
pages = {473-492},
year = {2024},
issn = {0377-2217},
doi = {10.1016/j.ejor.2023.07.014},
author = {Mingjie Li and Jin-Kao Hao and Qinghua Wu}
}

@article{Kallestad2023,
title = {A general deep reinforcement learning hyperheuristic framework for solving combinatorial optimization problems},
journal = {European Journal of Operational Research},
volume = {309},
number = {1},
pages = {446-468},
year = {2023},
issn = {0377-2217},
doi = {10.1016/j.ejor.2023.01.017},
author = {Jakob Kallestad and Ramin Hasibi and Ahmad Hemmati and Kenneth S{\"o}rensen}
}

@article{Yan2023,
title = {An online reinforcement learning approach to charging and order-dispatching optimization for an e-hailing electric vehicle fleet},
journal = {European Journal of Operational Research},
volume = {310},
number = {3},
pages = {1218-1233},
year = {2023},
issn = {0377-2217},
doi = {10.1016/j.ejor.2023.03.039},
author = {Pengyu Yan and Kaize Yu and Xiuli Chao and Zhibin Chen}
}

@article{Karimi2023,
title = {Learning to select operators in meta-heuristics: An integration of Q-learning into the iterated greedy algorithm for the permutation flowshop scheduling problem},
journal = {European Journal of Operational Research},
volume = {304},
number = {3},
pages = {1296-1330},
year = {2023},
issn = {0377-2217},
doi = {10.1016/j.ejor.2022.03.054},
author = {Maryam Karimi-Mamaghan and Mehrdad Mohammadi and Bastien Pasdeloup and Patrick Meyer}
}

@article{Chen2022,
title = {Deep Q-learning for same-day delivery with vehicles and drones},
journal = {European Journal of Operational Research},
volume = {298},
number = {3},
pages = {939-952},
year = {2022},
issn = {0377-2217},
doi = {10.1016/j.ejor.2021.06.021},
author = {Xinwei Chen and Marlin W. Ulmer and Barrett W. Thomas}
}

@article{Lee2021,
title = {Multi-agent reinforcement learning algorithm to solve a partially-observable multi-agent problem in disaster response},
journal = {European Journal of Operational Research},
volume = {291},
number = {1},
pages = {296-308},
year = {2021},
issn = {0377-2217},
doi = {10.1016/j.ejor.2020.09.018},
author = {Hyun-Rok Lee and Taesik Lee}
}

@article{Mosadegh2020,
  author  = {Mosadegh, H. and {Fatemi Ghomi}, S. M. T. and S{\"u}er, G. A.},
  title   = {Stochastic mixed-model assembly line sequencing problem: Mathematical modeling and {Q}-learning based simulated annealing hyper-heuristics},
  journal = {European Journal of Operational Research},
  volume  = {282},
  number  = {2},
  pages   = {530--544},
  year    = {2020},
  doi     = {10.1016/j.ejor.2019.09.021}
}

@article{Liu2020,
title = {Dynamic selective maintenance optimization for multi-state systems over a finite horizon: A deep reinforcement learning approach},
journal = {European Journal of Operational Research},
volume = {283},
number = {1},
pages = {166-181},
year = {2020},
issn = {0377-2217},
doi = {10.1016/j.ejor.2019.10.049},
author = {Yu Liu and Yiming Chen and Tao Jiang}
}

@article{Martin2016,
title = {A multi-agent based cooperative approach to scheduling and routing},
journal = {European Journal of Operational Research},
volume = {254},
number = {1},
pages = {169-178},
year = {2016},
issn = {0377-2217},
doi = {10.1016/j.ejor.2016.02.045},
author = {Simon Martin and Djamila Ouelhadj and Patrick Beullens and Ender Ozcan and Angel A. Juan and Edmund K. Burke}
}

@article{Brammer2022,
title = {Permutation flow shop scheduling with multiple lines and demand plans using reinforcement learning},
journal = {European Journal of Operational Research},
volume = {299},
number = {1},
pages = {75-86},
year = {2022},
issn = {0377-2217},
doi = {10.1016/j.ejor.2021.08.007},
author = {Janis Brammer and Bernhard Lutz and Dirk Neumann}
}

@article{Hasselt2016,
  author  = {{Van Hasselt}, Hado and Guez, Arthur and Silver, David},
  title   = {Deep reinforcement learning with double {Q}-learning},
  journal = {Proceedings of the AAAI Conference on Artificial Intelligence},
  volume  = {30},
  number  = {1},
  year    = {2016},
  doi     = {10.1609/aaai.v30i1.10295}
}

@InProceedings{Wang2016,
  title = {Dueling Network Architectures for Deep Reinforcement Learning},
  author = {Wang, Ziyu and Schaul, Tom and Hessel, Matteo and
            {Van Hasselt}, Hado and Lanctot, Marc and {De Freitas}, Nando},
  booktitle = {Proceedings of the 33rd International Conference on Machine Learning},
  pages = {1995--2003},
  year = {2016},
  editor = {Balcan, Maria Florina and Weinberger, Kilian Q.},
  volume = {48},
  series = {Proceedings of Machine Learning Research},
  publisher = {PMLR},
  address = {New York, New York, USA},
  url = {https://proceedings.mlr.press/v48/wangf16.html}
}

@article{Watkins1992,
  author  = {Watkins, Christopher J. C. H. and Dayan, Peter},
  title   = {{Q}-learning},
  journal = {Machine Learning},
  volume  = {8},
  pages   = {279--292},
  year    = {1992},
  doi     = {10.1007/BF00992698}
}

@ARTICLE{Arulkumaran2017,
  author={Arulkumaran, Kai and Deisenroth, Marc Peter and Brundage, Miles and Bharath, Anil Anthony},
  journal={IEEE Signal Processing Magazine}, 
  title={Deep reinforcement learning: A brief survey}, 
  year={2017},
  volume={34},
  number={6},
  pages={26-38},
  doi={10.1109/MSP.2017.2743240}
}

@article{Kaelbling1996,
  title={Reinforcement learning: A survey},
  author={Kaelbling, Leslie Pack and Littman, Michael L and Moore, Andrew W},
  journal={Journal of artificial intelligence research},
  volume={4},
  pages={237--285},
  year={1996},
  doi = {10.1613/jair.301}
}

@article{Shakya2023,
title = {Reinforcement learning algorithms: A brief survey},
journal = {Expert Systems with Applications},
volume = {231},
pages = {120495},
year = {2023},
issn = {0957-4174},
doi = {10.1016/j.eswa.2023.120495},
author = {Ashish Kumar Shakya and Gopinatha Pillai and Sohom Chakrabarty}
}

@article{Nian2020,
title = {A review On reinforcement learning: Introduction and applications in industrial process control},
journal = {Computers \& Chemical Engineering},
volume = {139},
pages = {106886},
year = {2020},
issn = {0098-1354},
doi = {10.1016/j.compchemeng.2020.106886},
author = {Rui Nian and Jinfeng Liu and Biao Huang}
}

@article{Ladosz2022,
title = {Exploration in deep reinforcement learning: A survey},
journal = {Information Fusion},
volume = {85},
pages = {1-22},
year = {2022},
issn = {1566-2535},
doi = {10.1016/j.inffus.2022.03.003},
author = {Pawel Ladosz and Lilian Weng and Minwoo Kim and Hyondong Oh}
}

@book{sutton2018,
  title     = {Reinforcement learning: An introduction},
  author    = {Sutton, Richard S. and Barto, Andrew G.},
  year      = {2018},
  edition   = {2nd},
  publisher = {MIT Press},
  address   = {Cambridge, MA},
  isbn      = {978-0262039246}
}

@InProceedings{Mousavi2016,
author={Mousavi, Seyed Sajad
and Schukat, Michael
and Howley, Enda},
editor={Bi, Yaxin
and Kapoor, Supriya
and Bhatia, Rahul},
title={Deep reinforcement learning: An overview},
booktitle={Proceedings of SAI Intelligent Systems Conference (IntelliSys) 2016},
year={2018},
publisher={Springer International Publishing},
address={Cham},
pages={426--440},
doi = {10.1007/978-3-319-56991-8_32}
}

@ARTICLE{Gu2024,
  author={Gu, Shangding and Yang, Long and Du, Yali and Chen, Guang and Walter, Florian and Wang, Jun and Knoll, Alois},
  journal={IEEE Transactions on Pattern Analysis and Machine Intelligence}, 
  title={A review of safe reinforcement learning: Methods, theories, and applications}, 
  year={2024},
  volume={46},
  number={12},
  pages={11216-11235},
  doi={10.1109/TPAMI.2024.3457538}
}

@article{Recht2019,
   author = {Recht, Benjamin},
   title = {A tour of reinforcement learning: The view from continuous control}, 
   journal= {Annual Review of Control, Robotics, and Autonomous Systems},
   year = {2019},
   volume = {2},
   number = {Volume 2, 2019},
   pages = {253-279},
   doi = {10.1146/annurev-control-053018-023825}
  }

@article{Luo2024,
  title={A survey on model-based reinforcement learning},
  author={Luo, Fan-Ming and Xu, Tian and Lai, Hang and Chen, Xiong-Hui and Zhang, Weinan and Yu, Yang},
  journal={Science China Information Sciences},
  volume={67},
  number={2},
  pages={121101},
  year={2024},
  publisher={Springer},
  doi = {10.1007/s11432-022-3696-5}
}

@article{Yu2023,
  author    = {Yu, Chao and Liu, Jiming and Nemati, Shamim and Yin, Guosheng},
  title     = {Reinforcement learning in healthcare: A survey},
  journal   = {ACM Computing Surveys},
  volume    = {55},
  number    = {1},
  pages     = {1--36},
  articleno = {5},
  year      = {2021},
  doi       = {10.1145/3477600}
}

@article{Garnier2021,
title = {A review on deep reinforcement learning for fluid mechanics},
journal = {Computers \& Fluids},
volume = {225},
pages = {104973},
year = {2021},
issn = {0045-7930},
doi = {10.1016/j.compfluid.2021.104973},
author = {Paul Garnier and Jonathan Viquerat and Jean Rabault and Aurélien Larcher and Alexander Kuhnle and Elie Hachem}
}

@inproceedings{Mnih2016,
  title = {Asynchronous methods for deep reinforcement learning},
  author = {Mnih, Volodymyr and Badia, Adria Puigdomenech and Mirza, Mehdi and Graves, Alex and Lillicrap, Timothy and Harley, Tim and Silver, David and Kavukcuoglu, Koray},
  booktitle = {Proceedings of the 33rd International Conference on Machine Learning},
  pages = {1928--1937},
  year = {2016},
  volume = {48},
  series = {Proceedings of Machine Learning Research},
  publisher = {PMLR},
  url = {https://proceedings.mlr.press/v48/mniha16.html}
}

@article{Wang2020,
  title={Deep reinforcement learning: a survey},
  author={Wang, Haonan and Liu, Ning and Zhang, Yiyun and Feng, Dawei and Huang, Feng and Li, Dongsheng and Zhang, Yiming},
  journal={Frontiers of Information Technology \& Electronic Engineering},
  volume={21},
  number={12},
  pages={1726--1744},
  year={2020},
  doi={10.1631/FITEE.1900533}
}

@article{Wang2020control,
title = {Reinforcement learning for building controls: The opportunities and challenges},
journal = {Applied Energy},
volume = {269},
pages = {115036},
year = {2020},
issn = {0306-2619},
doi = {10.1016/j.apenergy.2020.115036},
author = {Zhe Wang and Tianzhen Hong}
}

@incollection{Busoniu2010,
  author    = {Bu{\c{s}}oniu, Lucian and Babu{\v{s}}ka, Robert and {De Schutter}, Bart},
  title     = {Multi-agent reinforcement learning: An overview},
  booktitle = {Innovations in Multi-Agent Systems and Applications - 1},
  editor    = {Srinivasan, Dipti and Jain, Lakhmi C.},
  series    = {Studies in Computational Intelligence},
  volume    = {310},
  pages     = {183--221},
  year      = {2010},
  publisher = {Springer},
  address   = {Berlin, Heidelberg},
  doi       = {10.1007/978-3-642-14435-6_7}
}

@inproceedings{Vijay1999,
  author    = {Konda, Vijay R. and Tsitsiklis, John N.},
  title     = {Actor-Critic Algorithms},
  booktitle = {Advances in Neural Information Processing Systems},
  editor    = {Solla, S. and Leen, T. and M{\"u}ller, K.},
  volume    = {12},
  year      = {1999},
  publisher = {MIT Press},
  url       = {https://proceedings.neurips.cc/paper_files/paper/1999/file/6449f44a102fde848669bdd9eb6b76fa-Paper.pdf}
}

@article{Singh2022,
  title={Reinforcement learning in robotic applications: a comprehensive survey},
  author={Singh, Bharat and Kumar, Rajesh and Singh, Vinay Pratap},
  journal={Artificial Intelligence Review},
  volume={55},
  number={2},
  pages={945--990},
  year={2022},
  doi={10.1007/s10462-021-09997-9}
}

@article{Hu2024,
  title={A review of research on reinforcement learning algorithms for multi-agents},
  author={Hu, Kai and Li, Mingyang and Song, Zhiqiang and Xu, Keer and Xia, Qingfeng and Sun, Ning and Zhou, Peng and Xia, Min},
  journal={Neurocomputing},
  volume={599},
  pages={128068},
  year={2024},
  doi={10.1016/j.neucom.2024.128068}
}

@ARTICLE{Sureview2025,
  author={Su, Tong and Wu, Tong and Zhao, Junbo and Scaglione, Anna and Xie, Le},
  journal={Proceedings of the IEEE}, 
  title={A review of safe reinforcement learning methods for modern power systems}, 
  year={2025},
  volume={113},
  number={3},
  pages={213-255},
  doi={10.1109/JPROC.2025.3584656}
}

@article{Karimi2022,
title = {Machine learning at the service of meta-heuristics for solving combinatorial optimization problems: A state-of-the-art},
journal = {European Journal of Operational Research},
volume = {296},
number = {2},
pages = {393-422},
year = {2022},
issn = {0377-2217},
doi = {10.1016/j.ejor.2021.04.032},
author = {Maryam Karimi-Mamaghan and Mehrdad Mohammadi and Patrick Meyer and Amir Mohammad Karimi-Mamaghan and El-Ghazali Talbi}
}

@article{Song2019,
  title={A review on the self and dual interactions between machine learning and optimisation},
  author={Song, Heda and Triguero, Isaac and {\"O}zcan, Ender},
  journal={Progress in Artificial Intelligence},
  volume={8},
  number={2},
  pages={143--165},
  year={2019},
  doi={10.1007/s13748-019-00185-z}
}

@article{Talbi2016,
  title={Combining metaheuristics with mathematical programming, constraint programming and machine learning},
  author={Talbi, El-Ghazali},
  journal={Annals of Operations Research},
  volume={240},
  number={1},
  pages={171--215},
  year={2016},
  doi={10.1007/s10479-015-2034-y}
}

@article{Bengio2021,
title = {Machine learning for combinatorial optimization: A methodological tour d’horizon},
journal = {European Journal of Operational Research},
volume = {290},
number = {2},
pages = {405-421},
year = {2021},
issn = {0377-2217},
doi = {10.1016/j.ejor.2020.07.063},
author = {Yoshua Bengio and Andrea Lodi and Antoine Prouvost}
}

@inproceedings{Vinyals2015,
 author = {Vinyals, Oriol and Fortunato, Meire and Jaitly, Navdeep},
 booktitle = {Advances in Neural Information Processing Systems},
 editor = {C. Cortes and N. Lawrence and D. Lee and M. Sugiyama and R. Garnett},
 publisher = {Curran Associates, Inc.},
 title = {Pointer Networks},
 url = {https://proceedings.neurips.cc/paper_files/paper/2015/file/29921001f2f04bd3baee84a12e98098f-Paper.pdf},
 volume = {28},
 year = {2015}
}

@InProceedings{Vaswani2018,
  title = 	 {Image Transformer},
  author =       {Parmar, Niki and Vaswani, Ashish and Uszkoreit, Jakob and Kaiser, Lukasz and Shazeer, Noam and Ku, Alexander and Tran, Dustin},
  booktitle = 	 {Proceedings of the 35th International Conference on Machine Learning},
  pages = 	 {4055--4064},
  year = 	 {2018},
  editor = 	 {Dy, Jennifer and Krause, Andreas},
  volume = 	 {80},
  series = 	 {Proceedings of Machine Learning Research},
  month = 	 {10--15 Jul},
  publisher =    {PMLR},
  url = 	 {https://proceedings.mlr.press/v80/parmar18a.html}
}

@book{powell2007approximate,
  title={Approximate dynamic programming: Solving the curses of dimensionality},
  author={Powell, Warren B},
  volume={703},
  year={2007},
  publisher={John Wiley \& Sons}
}

@inproceedings{Lillicrap2015,
  title = {Continuous control with deep reinforcement learning},
  author = {Lillicrap, Timothy P. and Hunt, Jonathan J. and Pritzel, Alexander and Heess, Nicolas and Erez, Tom and Tassa, Yuval and Silver, David and Wierstra, Daan},
  booktitle = {Proceedings of the 4th International Conference on Learning Representations},
  year = {2016}
}

@inproceedings{Haarnoja2018,
  title = {Soft actor-critic: Off-policy maximum entropy deep reinforcement learning with a stochastic actor},
  author = {Haarnoja, Tuomas and Zhou, Aurick and Abbeel, Pieter and Levine, Sergey},
  booktitle = {Proceedings of the 35th International Conference on Machine Learning},
  pages = {1861--1870},
  year = {2018},
  volume = {80},
  series = {Proceedings of Machine Learning Research},
  publisher = {PMLR},
  url = {https://proceedings.mlr.press/v80/haarnoja18b.html}
}

@inproceedings{Fujimoto2018,
  title = {Addressing function approximation error in actor-critic methods},
  author = {Fujimoto, Scott and van Hoof, Herke and Meger, David},
  booktitle = {Proceedings of the 35th International Conference on Machine Learning},
  pages = {1587--1596},
  year = {2018},
  volume = {80},
  series = {Proceedings of Machine Learning Research},
  publisher = {PMLR},
  url = {https://proceedings.mlr.press/v80/fujimoto18a.html}
}

@article{Sutton1991,
author = {Sutton, Richard S.},
title = {Dyna, an integrated architecture for learning, planning, and reacting},
year = {1991},
issue_date = {Aug. 1991},
publisher = {Association for Computing Machinery},
address = {New York, NY, USA},
volume = {2},
number = {4},
issn = {0163-5719},
doi = {10.1145/122344.122377},
journal = {SIGART Bull.},
pages = {160–163},
numpages = {4}
}

@inproceedings{Buckman2019,
  title = {Sample-efficient reinforcement learning with stochastic ensemble value expansion},
  author = {Buckman, Jacob and Hafner, Danijar and Tucker, George and Brevdo, Eugene and Lee, Honglak},
  booktitle = {Advances in Neural Information Processing Systems},
  volume = {31},
  year = {2018},
  publisher = {Curran Associates, Inc.},
  url = {https://proceedings.neurips.cc/paper_files/paper/2018/hash/f02208a057804ee16ac72ff4d3cec53b-Abstract.html}
}

@inproceedings{Janner2019,
 author = {Janner, Michael and Fu, Justin and Zhang, Marvin and Levine, Sergey},
 booktitle = {Advances in Neural Information Processing Systems},
 editor = {H. Wallach and H. Larochelle and A. Beygelzimer and F. d\textquotesingle Alch\'{e}-Buc and E. Fox and R. Garnett},
 publisher = {Curran Associates, Inc.},
 title = {When to trust your model: Model-based policy optimization},
 url = {https://proceedings.neurips.cc/paper_files/paper/2019/file/5faf461eff3099671ad63c6f3f094f7f-Paper.pdf},
 volume = {32},
 year = {2019}
}

@article{WangSongLiuZhao2026,
  title   = {Reinforcement learning algorithm for reusable resource allocation with unknown rental time distribution},
  author  = {Wang, Ziwei and Song, Jie and Liu, Yixuan and Zhao, Jingtong},
  journal = {European Journal of Operational Research},
  volume  = {331},
  number  = {1},
  pages   = {186--199},
  year    = {2026},
  doi     = {10.1016/j.ejor.2025.09.012}
}

@inproceedings{Schulman2017TRPO,
  title = {Trust region policy optimization},
  author = {Schulman, John and Levine, Sergey and Abbeel, Pieter and Jordan, Michael and Moritz, Philipp},
  booktitle = {Proceedings of the 32nd International Conference on Machine Learning},
  pages = {1889--1897},
  year = {2015},
  volume = {37},
  series = {Proceedings of Machine Learning Research},
  publisher = {PMLR},
  url = {https://proceedings.mlr.press/v37/schulman15.html}
}

@article{Panda2024,
title = {Dynamic resource matching in manufacturing using deep reinforcement learning},
journal = {European Journal of Operational Research},
volume = {318},
number = {2},
pages = {408-423},
year = {2024},
issn = {0377-2217},
doi = {10.1016/j.ejor.2024.05.027},
author = {Saunak Kumar Panda and Yisha Xiang and Ruiqi Liu}
}

@article{Alfonso2024,
title = {Optimizing credit limit adjustments under adversarial goals using reinforcement learning},
journal = {European Journal of Operational Research},
volume = {315},
number = {2},
pages = {802-817},
year = {2024},
issn = {0377-2217},
doi = {10.1016/j.ejor.2023.12.025},
author = {Sherly Alfonso-Sánchez and Jesús Solano and Alejandro Correa-Bahnsen and Kristina P. Sendova and Cristián Bravo}
}

@article{Verleijsdonk2024,
title = {Scalable policies for the dynamic traveling multi-maintainer problem with alerts},
journal = {European Journal of Operational Research},
volume = {319},
number = {1},
pages = {121-134},
year = {2024},
issn = {0377-2217},
doi = {10.1016/j.ejor.2024.05.049},
author = {Peter Verleijsdonk and Willem {Van Jaarsveld} and Stella Kapodistria}
}

@article{Monaci2024,
title = {An actor-critic algorithm with policy gradients to solve the job shop scheduling problem using deep double recurrent agents},
journal = {European Journal of Operational Research},
volume = {312},
number = {3},
pages = {910-926},
year = {2024},
issn = {0377-2217},
doi = {10.1016/j.ejor.2023.07.037},
author = {Marta Monaci and Valerio Agasucci and Giorgio Grani}
}

@article{Chen2023,
title = {Same-day delivery with fair customer service},
journal = {European Journal of Operational Research},
volume = {308},
number = {2},
pages = {738-751},
year = {2023},
issn = {0377-2217},
doi = {10.1016/j.ejor.2022.12.009},
author = {Xinwei Chen and Tong Wang and Barrett W. Thomas and Marlin W. Ulmer},
}

@article{Benlic2017,
title = {A hybrid breakout local search and reinforcement learning approach to the vertex separator problem},
journal = {European Journal of Operational Research},
volume = {261},
number = {3},
pages = {803-818},
year = {2017},
issn = {0377-2217},
doi = {10.1016/j.ejor.2017.01.023},
author = {Una Benlic and Michael G. Epitropakis and Edmund K. Burke}
}

@article{Rolim2025,
title = {On the integration of reinforcement learning and simulated annealing for the parallel batch scheduling problem with setups},
journal = {European Journal of Operational Research},
volume = {326},
number = {2},
pages = {220-233},
year = {2025},
issn = {0377-2217},
doi = {10.1016/j.ejor.2025.04.042},
author = {Gustavo Alencar Rolim and Caio Paziani Tomazella and Marcelo Seido Nagano}
}

@article{Li2022,
title = {Learning-driven feasible and infeasible tabu search for airport gate assignment},
journal = {European Journal of Operational Research},
volume = {302},
number = {1},
pages = {172-186},
year = {2022},
issn = {0377-2217},
doi = {10.1016/j.ejor.2021.12.019},
author = {Mingjie Li and Jin-Kao Hao and Qinghua Wu}
}

@article{Zhang2023,
title = {Synchromodal freight transport re-planning under service time uncertainty: An online model-assisted reinforcement learning},
journal = {Transportation Research Part C: Emerging Technologies},
volume = {156},
pages = {104355},
year = {2023},
issn = {0968-090X},
doi = {10.1016/j.trc.2023.104355},
author = {Yimeng Zhang and Rudy R. Negenborn and Bilge Atasoy}
}

@article{Semrov2016,
title = {Reinforcement learning approach for train rescheduling on a single-track railway},
journal = {Transportation Research Part B: Methodological},
volume = {86},
pages = {250-267},
year = {2016},
issn = {0191-2615},
doi = {10.1016/j.trb.2016.01.004},
author = {D. Semrov and R. Marsetič and M. Žura and L. Todorovski and A. Srdic}
}

@techreport{Wu2025,
  title = {Deep Reinforcement Learning for Combinatorial Optimization: A Tutorial},
  author = {Wu, Yaoxin and Bukhsh, Zaharah and Zhang, Yingqian},
  year = {2025},
  institution = {Eindhoven University of Technology},
  type = {Technical report},
  note = {Unpublished tutorial, retained as a background reference},
  url = {https://research.tue.nl/en/publications/deep-reinforcement-learning-for-combinatorial-optimization-a-tuto/}
}

@article{Ahamed2021,
title = {Deep reinforcement learning for crowdsourced urban delivery},
journal = {Transportation Research Part B: Methodological},
volume = {152},
pages = {227-257},
year = {2021},
issn = {0191-2615},
doi = {10.1016/j.trb.2021.08.015},
author = {Tanvir Ahamed and Bo Zou and Nahid Parvez Farazi and Theja Tulabandhula}
}

@article{Litrain2022,
title = {Train timetabling with the general learning environment and multi-agent deep reinforcement learning},
journal = {Transportation Research Part B: Methodological},
volume = {157},
pages = {230-251},
year = {2022},
doi = {10.1016/j.trb.2022.02.006},
author = {Wenqing Li and Shaoquan Ni}
}

@article{Yu2025,
title = {Interpretable state-space model of urban dynamics for human-machine collaborative transportation planning},
journal = {Transportation Research Part B: Methodological},
volume = {192},
pages = {103134},
year = {2025},
issn = {0191-2615},
doi = {10.1016/j.trb.2024.103134},
author = {Jiangbo Yu and Michael F. Hyland}
}

@article{Ying2024,
  author  = {Ying, Chengshuo and Chow, Andy H. F. and Yan, Yimo and Kuo, Yong-Hong and Wang, Shouyang},
  title   = {Adaptive rescheduling of rail transit services with short-turnings under disruptions via a multi-agent deep reinforcement learning approach},
  journal = {Transportation Research Part B: Methodological},
  volume  = {188},
  pages   = {103067},
  year    = {2024},
  doi     = {10.1016/j.trb.2024.103067}
}

@article{Mao2025,
author = {Mao, Weichao and Zhang, Kaiqing and Zhu, Ruihao and Simchi-Levi, David and Ba\c{s}ar, Tamer},
title = {Model-free nonstationary reinforcement learning: Near-optimal regret and applications in multiagent reinforcement learning and inventory control},
journal = {Management Science},
volume = {71},
number = {2},
pages = {1564-1580},
year = {2025},
doi = {10.1287/mnsc.2022.02533}
}

@article{Wang2023,
author = {Wang, Wen and Li, Beibei and Luo, Xueming and Wang, Xiaoyi},
title = {Deep reinforcement learning for sequential targeting},
journal = {Management Science},
volume = {69},
number = {9},
pages = {5439-5460},
year = {2023},
doi = {10.1287/mnsc.2022.4621}
}

@article{Ying2020,
  author  = {Ying, Chengshuo and Chow, Andy H. F. and Chin, Kwaisang},
  title   = {An actor-critic deep reinforcement learning approach for metro train scheduling with rolling stock circulation under stochastic demand},
  journal = {Transportation Research Part B: Methodological},
  volume  = {140},
  pages   = {210--235},
  year    = {2020},
  doi     = {10.1016/j.trb.2020.08.005}
}

@article{Guo2022,
  author    = {Guo, W. and Atasoy, B. and Negenborn, R. R.},
  title     = {Global synchromodal shipment matching problem with dynamic and stochastic travel times: a reinforcement learning approach},
  journal   = {Annals of Operations Research},
  volume    = {350},
  number    = {1},
  pages     = {63--94},
  year      = {2025},
  publisher = {Springer},
  doi       = {10.1007/s10479-021-04489-z}
}

@InProceedings{Enders2023,
  title = 	 {Hybrid multi-agent deep reinforcement learning for autonomous mobility on demand systems},
  author =       {Enders, Tobias and Harrison, James and Pavone, Marco and Schiffer, Maximilian},
  booktitle = 	 {Proceedings of The 5th Annual Learning for Dynamics and Control Conference},
  pages = 	 {1284--1296},
  year = 	 {2023},
  editor = 	 {Matni, Nikolai and Morari, Manfred and Pappas, George J.},
  volume = 	 {211},
  series = 	 {Proceedings of Machine Learning Research},
  month = 	 {15--16 Jun},
  publisher =    {PMLR},
  url = 	 {https://proceedings.mlr.press/v211/enders23a.html}
}

@inproceedings{Rigter2022,
  author    = {Rigter, Marc and Lacerda, Bruno and Hawes, Nick},
  title     = {{RAMBO-RL}: Robust Adversarial Model-Based Offline Reinforcement Learning},
  booktitle = {Advances in Neural Information Processing Systems},
  editor    = {S. Koyejo and S. Mohamed and A. Agarwal and D. Belgrave and K. Cho and A. Oh},
  volume    = {35},
  pages     = {16082--16097},
  publisher = {Curran Associates, Inc.},
  year      = {2022},
  doi       = {10.52202/068431-1170},
  url       = {https://proceedings.neurips.cc/paper_files/paper/2022/file/6691c5e4a199b72dffd9c90acb63bcd6-Paper-Conference.pdf}
}

@article{Ma2024,
  title   = {Efficient and scalable reinforcement learning for large-scale network control},
  author  = {Chengdong Ma and Aming Li and Yali Du and Hao Dong and Yaodong Yang},
  journal = {Nature Machine Intelligence},
  year    = {2024},
  volume  = {6},
  pages   = {1006--1020},
  doi     = {10.1038/s42256-024-00879-7}
}

@inproceedings{Jayant2022,
 author = {Jayant, Ashish K and Bhatnagar, Shalabh},
 booktitle = {Advances in Neural Information Processing Systems},
 editor = {S. Koyejo and S. Mohamed and A. Agarwal and D. Belgrave and K. Cho and A. Oh},
 pages = {24432--24445},
 publisher = {Curran Associates, Inc.},
 title = {Model-based safe deep reinforcement learning via a constrained proximal policy optimization algorithm},
 url = {https://proceedings.neurips.cc/paper_files/paper/2022/file/9a8eb202c060b7d81f5889631cbcd47e-Paper-Conference.pdf},
 volume = {35},
 year = {2022}
}

@inproceedings{Wei2017,
author = {Wei, Zeng and Xu, Jun and Lan, Yanyan and Guo, Jiafeng and Cheng, Xueqi},
title = {Reinforcement learning to rank with Markov Decision Process},
year = {2017},
isbn = {9781450350228},
publisher = {Association for Computing Machinery},
address = {New York, NY, USA},
doi = {10.1145/3077136.3080685},
booktitle = {Proceedings of the 40th International ACM SIGIR Conference on Research and Development in Information Retrieval},
pages = {945–948},
numpages = {4},
location = {Shinjuku, Tokyo, Japan},
series = {SIGIR '17}
}

@inproceedings{Park2024,
  title = {Is value learning really the main bottleneck in offline {RL}?},
  author = {Park, Seohong and Frans, Kevin and Levine, Sergey and Kumar, Aviral},
  booktitle = {Advances in Neural Information Processing Systems},
  volume = {37},
  pages = {79029--79056},
  year = {2024},
  doi = {10.52202/079017-2509},
  url = {https://proceedings.neurips.cc/paper_files/paper/2024/hash/8ffb4e3118280a66b192b6f06e0e2596-Abstract-Conference.html}
}

@article{Zhu2021,
title = {A mean-field Markov decision process model for spatial-temporal subsidies in ride-sourcing markets},
journal = {Transportation Research Part B: Methodological},
volume = {150},
pages = {540-565},
year = {2021},
issn = {0191-2615},
doi = {10.1016/j.trb.2021.06.014},
author = {Zheng Zhu and Jintao Ke and Hai Wang}
}

@article{Ying2022,
  author  = {Ying, Chengshuo and Chow, Andy H. F. and Nguyen, Hoa T. M. and Chin, KwaiSang},
  title   = {Multi-agent deep reinforcement learning for adaptive coordinated metro service operations with flexible train composition},
  journal = {Transportation Research Part B: Methodological},
  volume  = {161},
  pages   = {36--59},
  year    = {2022},
  doi     = {10.1016/j.trb.2022.05.001}
}

@article{Oroojlooyjadid2022,
author = {Oroojlooyjadid, Afshin and Nazari, MohammadReza and Snyder, Lawrence V. and Tak\'{a}\v{c}, Martin},
title = {A deep Q-network for the beer game: Deep reinforcement learning for inventory optimization},
journal = {Manufacturing \& Service Operations Management},
volume = {24},
number = {1},
pages = {285-304},
year = {2022},
doi = {10.1287/msom.2020.0939}
}

@article{Meijer2021,
author = {{De Meijer}, Frank and Sotirov, Renata},
title = {SDP-based bounds for the quadratic cycle cover problem via cutting-plane augmented lagrangian methods and reinforcement learning},
journal = {INFORMS Journal on Computing},
volume = {33},
number = {4},
pages = {1262-1276},
year = {2021},
doi = {10.1287/ijoc.2021.1075}
}

@article{Cappart2022,
author = {Cappart, Quentin and Bergman, David and Rousseau, Louis-Martin and Pr\'{e}mont-Schwarz, Isabeau and Parjadis, Augustin},
title = {Improving variable orderings of approximate decision diagrams using reinforcement learning},
journal = {INFORMS Journal on Computing},
volume = {34},
number = {5},
pages = {2552-2570},
year = {2022},
doi = {10.1287/ijoc.2022.1194}
}

@article{Li2024TS,
author = {Li, Yuanyuan and Archetti, Claudia and Ljubi\'{c}, Ivana},
title = {Reinforcement learning approaches for the orienteering problem with stochastic and dynamic release dates},
journal = {Transportation Science},
volume = {58},
number = {5},
pages = {1143-1165},
year = {2024},
doi = {10.1287/trsc.2022.0366}
}

@article{Harsha2025,
author = {Harsha, Pavithra and Jagmohan, Ashish and Kalagnanam, Jayant and Quanz, Brian and Singhvi, Divya},
title = {Deep policy iteration with integer programming for inventory management},
journal = {Manufacturing \& Service Operations Management},
volume = {27},
number = {2},
pages = {369-388},
year = {2025},
doi = {10.1287/msom.2022.0617}
}

@article{Chohlas2025,
author = {Chohlas-Wood, Alex and Coots, Madison and Zhu, Henry and Brunskill, Emma and Goel, Sharad},
title = {Learning to Be Fair: A Consequentialist Approach to Equitable Decision Making},
journal = {Management Science},
volume = {72},
number = {1},
pages = {456-473},
year = {2026},
doi = {10.1287/mnsc.2022.00345}
}

@article{Cappart2021, title={Combining reinforcement learning and constraint programming for combinatorial optimization}, volume={35}, DOI={10.1609/aaai.v35i5.16484}, number={5}, journal={Proceedings of the AAAI Conference on Artificial Intelligence}, author={Cappart, Quentin and Moisan, Thierry and Rousseau, Louis-Martin and Prémont-Schwarz, Isabeau and Cire, Andre A.}, year={2021}, month={May}, pages={3677-3687} }

@article{Tassel2023, title={An end-to-end reinforcement learning approach for job-shop scheduling problems based on constraint programming}, volume={33}, DOI={10.1609/icaps.v33i1.27243}, number={1}, journal={Proceedings of the International Conference on Automated Planning and Scheduling}, author={Tassel, Pierre and Gebser, Martin and Schekotihin, Konstantin}, year={2023}, month={Jul.}, pages={614-622} }

@article{Lamghari2020,
title = {Hyper-heuristic approaches for strategic mine planning under uncertainty},
journal = {Computers \& Operations Research},
volume = {115},
pages = {104590},
year = {2020},
issn = {0305-0548},
doi = {10.1016/j.cor.2018.11.010},
author = {Amina Lamghari and Roussos Dimitrakopoulos}
}

@article{Alicastro2021,
title = {A reinforcement learning iterated local search for makespan minimization in additive manufacturing machine scheduling problems},
journal = {Computers \& Operations Research},
volume = {131},
pages = {105272},
year = {2021},
issn = {0305-0548},
doi = {10.1016/j.cor.2021.105272},
author = {Mirko Alicastro and Daniele Ferone and Paola Festa and Serena Fugaro and Tommaso Pastore}
}

@article{Zou2024,
title = {A reinforcement learning guided hybrid evolutionary algorithm for the latency location routing problem},
journal = {Computers \& Operations Research},
volume = {170},
pages = {106758},
year = {2024},
issn = {0305-0548},
doi = {10.1016/j.cor.2024.106758},
author = {Yuji Zou and Jin-Kao Hao and Qinghua Wu}
}

@article{Cui2024,
title = {A hybrid genetic algorithm based on reinforcement learning for the energy-aware production scheduling in the photovoltaic glass industry},
journal = {Computers \& Operations Research},
volume = {163},
pages = {106521},
year = {2024},
issn = {0305-0548},
doi = {10.1016/j.cor.2023.106521},
author = {Weiwei Cui and Biao Yuan}
}

@article{Lu2024,
title = {Learning driven three-phase search for the maximum independent union of cliques problem},
journal = {Computers \& Operations Research},
volume = {164},
pages = {106549},
year = {2024},
issn = {0305-0548},
doi = {10.1016/j.cor.2024.106549},
author = {Zhi Lu and Jian Gao and Jin-Kao Hao and Pingle Yang and Lixin Zhou}
}

@article{Zhao2025,
title = {Reinforcement learning-enhanced variable neighborhood search strategies for the k-clustering minimum biclique completion problem},
journal = {Computers \& Operations Research},
volume = {178},
pages = {107008},
year = {2025},
issn = {0305-0548},
doi = {10.1016/j.cor.2025.107008},
author = {Juntao Zhao and Mhand Hifi}
}

@article{An2026,
title = {Flexible scheduling of customized bus for green mega-events: A distributionally robust optimization approach},
journal = {Computers \& Operations Research},
volume = {185},
pages = {107249},
year = {2026},
issn = {0305-0548},
doi = {10.1016/j.cor.2025.107249},
author = {Xiaojie An and Xiang Li and Bowen Zhang}
}

@article{Li2024COR,
title = {Minimizing passenger waiting time in the multi-route bus fleet allocation problem through distributionally robust optimization and reinforcement learning},
journal = {Computers \& Operations Research},
volume = {164},
pages = {106568},
year = {2024},
issn = {0305-0548},
doi = {10.1016/j.cor.2024.106568},
author = {Xiang Li and Xiaojie An and Bowen Zhang}
}

@article{Zhang2022,
title = {A deep reinforcement learning based hyper-heuristic for combinatorial optimisation with uncertainties},
journal = {European Journal of Operational Research},
volume = {300},
number = {2},
pages = {418-427},
year = {2022},
issn = {0377-2217},
doi = {10.1016/j.ejor.2021.10.032},
author = {Yuchang Zhang and Ruibin Bai and Rong Qu and Chaofan Tu and Jiahuan Jin}
}

@InProceedings{Deudon2018,
  author    = {Deudon, Michel and Cournut, Pierre and
               Lacoste, Alexandre and Adulyasak, Yossiri and
               Rousseau, Louis-Martin},
  editor    = {{Van Hoeve}, Willem-Jan},
  title     = {Learning Heuristics for the {TSP} by Policy Gradient},
  booktitle = {Integration of Constraint Programming, Artificial Intelligence, and Operations Research},
  series    = {Lecture Notes in Computer Science},
  volume    = {10848},
  year      = {2018},
  publisher = {Springer International Publishing},
  address   = {Cham},
  pages     = {170--181},
  doi       = {10.1007/978-3-319-93031-2_12}
}

@InProceedings{Tang2020,
  title = 	 {Reinforcement learning for integer programming: learning to cut},
  author =       {Tang, Yunhao and Agrawal, Shipra and Faenza, Yuri},
  booktitle = 	 {Proceedings of the 37th International Conference on Machine Learning},
  pages = 	 {9367--9376},
  year = 	 {2020},
  editor = 	 {Daum{\'e} III, Hal and Singh, Aarti},
  volume = 	 {119},
  series = 	 {Proceedings of Machine Learning Research},
  publisher =    {PMLR},
  url = 	 {https://proceedings.mlr.press/v119/tang20a.html}
}

@ARTICLE{Li2021,
  author={Li, Kaiwen and Zhang, Tao and Wang, Rui},
  journal={IEEE Transactions on Cybernetics}, 
  title={Deep Reinforcement Learning for Multiobjective Optimization}, 
  year={2021},
  volume={51},
  number={6},
  pages={3103-3114},
  doi={10.1109/TCYB.2020.2977661}}

@article{Ma2025,
  author  = {Ma, Yujie and Xia, Xin and Liu, Peiyao and Zhang, Chen},
  title   = {Bilevel joint optimization for product design changes with a resilient supply chain based on deep reinforcement learning},
  journal = {International Journal of Production Economics},
  volume  = {290},
  pages   = {109791},
  year    = {2025},
  doi     = {10.1016/j.ijpe.2025.109791}
}

@article{Yu2026,
title = {A deep reinforcement learning approach for dynamic job-shop scheduling problem considering time variable and new job arrivals},
journal = {Computers \& Operations Research},
volume = {185},
pages = {107263},
year = {2026},
issn = {0305-0548},
doi = {10.1016/j.cor.2025.107263},
author = {Haoyang Yu and Wenbin Gu and Na Tang and Zhenyang Guo}
}

@article{Li2025,
title = {Robust optimization for adaptive bus service scheduling with adversarial reinforcement learning under demand uncertainties},
journal = {Transportation Research Part C: Emerging Technologies},
volume = {178},
pages = {105222},
year = {2025},
doi = {10.1016/j.trc.2025.105222},
author = {Guangyu Li and Andy H.F. Chow and Chengshuo Ying}
}

@article{Su2025,
title = {Enhancing feeder bus service coverage with Multi-Agent Reinforcement Learning: A case study in Hong Kong},
journal = {Transportation Research Part E: Logistics and Transportation Review},
volume = {196},
pages = {103997},
year = {2025},
issn = {1366-5545},
doi = {10.1016/j.tre.2025.103997},
author = {Yang Su and Hai Yang}
}

@INPROCEEDINGS{Huang2021,
  author={Huang, Xingshuai and Wu, Di and Jenkin, Michael and Boulet, Benoit},
  booktitle={2024 IEEE Canadian Conference on Electrical and Computer Engineering (CCECE)}, 
  title={Sample-Efficient Meta-RL for Traffic Signal Control}, 
  year={2024},
  pages={312-318},
  doi={10.1109/CCECE59415.2024.10667211}}

@ARTICLE{Devailly2022,
  author={Devailly, Fran{\c{c}}ois-Xavier and Larocque, Denis and Charlin, Laurent},
  journal={IEEE Open Journal of Intelligent Transportation Systems}, 
  title={Model-Based Graph Reinforcement Learning for Inductive Traffic Signal Control}, 
  year={2024},
  volume={5},
  pages={238-250},
  doi={10.1109/OJITS.2024.3376583}}

@ARTICLE{Wu2022,
  author={Wu, Jingda and Huang, Zhiyu and Lv, Chen},
  journal={IEEE Transactions on Intelligent Vehicles}, 
  title={Uncertainty-Aware Model-Based Reinforcement Learning: Methodology and Application in Autonomous Driving}, 
  year={2023},
  volume={8},
  number={1},
  pages={194-203},
  doi={10.1109/TIV.2022.3185159}}

@article{Diehl2022,
  title = {Uncertainty-Aware Model-Based Offline Reinforcement Learning for Automated Driving},
  author = {Diehl, Christopher and Sievernich, Timo Sebastian and Kr{\"u}ger, Martin and Hoffmann, Frank and Bertram, Torsten},
  journal = {IEEE Robotics and Automation Letters},
  year = {2023},
  volume = {8},
  number = {2},
  pages = {1167--1174},
  doi = {10.1109/LRA.2023.3236579}
}

@article{Gao2025,
  title = {Dream to Drive With Predictive Individual World Model},
  author = {Gao, Yinfeng and Zhang, Qichao and Ding, Da-Wei and Zhao, Dongbin},
  journal = {IEEE Transactions on Intelligent Vehicles},
  year = {2024},
  volume = {9},
  number = {12},
  pages = {8224--8238},
  doi = {10.1109/TIV.2024.3408830}
}

@ARTICLE{Zhang2024,
  author={Zhang, Zhengming and Huang, Yongming and Zhang, Cheng and Zheng, Qingbi and Yang, Luxi and You, Xiaohu},
  journal={IEEE Transactions on Communications}, 
  title={Digital twin-enhanced deep reinforcement learning for resource management in networks slicing}, 
  year={2024},
  volume={72},
  number={10},
  pages={6209-6224},
  doi={10.1109/TCOMM.2024.3395698}}

@article{Liu2025,
title = {Multi-Agent reinforcement learning framework for addressing demand-supply imbalance of shared autonomous electric vehicle},
journal = {Transportation Research Part E: Logistics and Transportation Review},
volume = {197},
pages = {104062},
year = {2025},
issn = {1366-5545},
doi = {10.1016/j.tre.2025.104062},
author = {Chengqi Liu and Zelin Wang and Zhiyuan Liu and Kai Huang}
}

@article{Ding2025,
title = {Scheduling intelligent charging robots for electric vehicle: A deep reinforcement learning approach},
journal = {Transportation Research Part E: Logistics and Transportation Review},
volume = {200},
pages = {104090},
year = {2025},
issn = {1366-5545},
doi = {10.1016/j.tre.2025.104090},
author = {Yi Ding and Ming Deng and Ginger Y. Ke and Yingjun Shen and Lianmin Zhang}
}

@article{Park2022,
author = {Kyu Tae Park and Seung-Woo Jeon and Sang Do Noh},
title = {Digital twin application with horizontal coordination for reinforcement-learning-based production control in a re-entrant job shop},
journal = {International Journal of Production Research},
volume = {60},
number = {7},
pages = {2151--2167},
year = {2022},
publisher = {Taylor \& Francis},
doi = {10.1080/00207543.2021.1884309}

}

@INPROCEEDINGS{Wu2021,
  author={Wu, Jingda and Huang, Zhiyu and Hang, Peng and Huang, Chao and De Boer, Niels and Lv, Chen},
  booktitle={2021 IEEE 1st International Conference on Digital Twins and Parallel Intelligence (DTPI)}, 
  title={Digital twin-enabled reinforcement learning for end-to-end autonomous driving}, 
  year={2021},
  pages={62-65},
  doi={10.1109/DTPI52967.2021.9540179}}

@ARTICLE{Wang2024,
  author={Wang, Lei and Liang, Hongbin and Mao, Guotao and Zhao, Dongmei and Liu, Qian and Yao, Yiting and Zhang, Han},
  journal={IEEE Transactions on Vehicular Technology}, 
  title={Resource allocation for dynamic platoon digital twin networks: A multi-agent deep reinforcement learning method}, 
  year={2024},
  volume={73},
  number={10},
  pages={15609-15620},
  doi={10.1109/TVT.2024.3414447}}

@ARTICLE{Xu2024,
  author={Xu, Siyu and Guan, Xin and Peng, Yu and Liu, Yang and Cui, Chen and Chen, Hongyang and Ohtsuki, Tomoaki and Han, Zhu},
  journal={IEEE Transactions on Network Science and Engineering}, 
  title={Deep reinforcement learning based data-driven mapping mechanism of digital twin for internet of energy}, 
  year={2024},
  volume={11},
  number={4},
  pages={3876-3890},
  doi={10.1109/TNSE.2024.3390797}
}

@ARTICLE{Liu2023,
  author={Liu, Qinghai and Tang, Lun and Wu, Ting and Chen, Qianbin},
  journal={IEEE Internet of Things Journal}, 
  title={Deep reinforcement learning for resource demand prediction and virtual function network migration in digital twin network}, 
  year={2023},
  volume={10},
  number={21},
  pages={19102-19116},
  doi={10.1109/JIOT.2023.3281678}
}

@ARTICLE{Schlappa2024,
  author={Schlappa, Martin and Hegemann, Jonas and Spinler, Stefan},
  journal={IEEE Transactions on Engineering Management}, 
  title={Optimizing control of waste incineration plants using reinforcement learning and digital twins}, 
  year={2024},
  volume={71},
  number={},
  pages={3076-3087},
  doi={10.1109/TEM.2022.3201434}
}

@ARTICLE{Tang2023,
  author={Tang, Xin and Li, Xiaohuan and Yu, Rong and Wu, Yuan and Ye, Jin and Tang, Fengzhu and Chen, Qian},
  journal={IEEE Internet of Things Journal}, 
  title={Digital-twin-assisted task assignment in multi-UAV systems: A deep reinforcement learning approach}, 
  year={2023},
  volume={10},
  number={17},
  pages={15362-15375},
  doi={10.1109/JIOT.2023.3263574}}

@article{Mnih2015,
  title = {Human-level control through deep reinforcement learning},
  author = {Mnih, Volodymyr and Kavukcuoglu, Koray and Silver, David and Rusu, Andrei A. and Veness, Joel and Bellemare, Marc G. and Graves, Alex and Riedmiller, Martin and Fidjeland, Andreas K. and Ostrovski, Georg and Petersen, Stig and Beattie, Charles and Sadik, Amir and Antonoglou, Ioannis and King, Helen and Kumaran, Dharshan and Wierstra, Daan and Legg, Shane and Hassabis, Demis},
  journal = {Nature},
  volume = {518},
  number = {7540},
  pages = {529--533},
  year = {2015},
  doi = {10.1038/nature14236}
}

@inproceedings{Kool2018,
  title = {Attention, Learn to Solve Routing Problems!},
  author = {Kool, Wouter and van Hoof, Herke and Welling, Max},
  booktitle = {International Conference on Learning Representations},
  year = {2019},
  url = {https://openreview.net/forum?id=ByxBFsRqYm}
}

@book{BirgeLouveaux2011,
  author    = {Birge, John R. and Louveaux, Fran{\c{c}}ois},
  title     = {Introduction to Stochastic Programming},
  edition   = {2},
  series    = {Springer Series in Operations Research and Financial Engineering},
  publisher = {Springer},
  address   = {New York},
  year      = {2011},
  doi       = {10.1007/978-1-4614-0237-4}
}

@article{Powell2019,
  author    = {Powell, Warren B.},
  title     = {A unified framework for stochastic optimization},
  journal   = {European Journal of Operational Research},
  volume    = {275},
  number    = {3},
  pages     = {795--821},
  year      = {2019},
  doi       = {10.1016/j.ejor.2018.07.014}
}

@article{LiSeru2025,
  title   = {Dynamic worker allocation in Seru production systems with actor--critic and pointer networks},
  author  = {Li, Dongni and Jin, Hongbo and Zhang, Yaoxin},
  journal = {European Journal of Operational Research},
  volume  = {324},
  number  = {1},
  pages   = {62--74},
  year    = {2025},
  doi     = {10.1016/j.ejor.2025.01.012}
}

@book{Powell2022,
  author    = {Powell, Warren B.},
  title     = {Reinforcement Learning and Stochastic Optimization: A Unified Framework for Sequential Decisions},
  publisher = {John Wiley \& Sons},
  address   = {Hoboken, NJ},
  year      = {2022},
  doi       = {10.1002/9781119815068}
}

@book{Shapiro2021,
  author    = {Shapiro, Alexander and Dentcheva, Darinka and Ruszczy{\'n}ski, Andrzej},
  title     = {Lectures on Stochastic Programming: Modeling and Theory},
  edition   = {3rd},
  publisher = {SIAM},
  address   = {Philadelphia},
  year      = {2021},
    doi       = {10.1137/1.9781611976595}

}

@article{Alcaraz2022,
  title = {Online model-based reinforcement learning for decision-making in long distance routes},
  author = {Alcaraz, Juan J. and Losilla, Fernando and Caballero-Arnaldos, Luis},
  journal = {Transportation Research Part E: Logistics and Transportation Review},
  volume = {164},
  pages = {102790},
  year = {2022},
  doi = {10.1016/j.tre.2022.102790}
}

@article{Leech2025,
  title = {Model-Based Deep Reinforcement Learning for Energy Efficient Routing of a Connected and Automated Vehicle},
  author = {Leech, David R. and Yoon, Hwan-Sik},
  journal = {Sustainability},
  volume = {17},
  number = {13},
  pages = {5727},
  year = {2025},
  doi = {10.3390/su17135727}
}

@inproceedings{Andersen2019,
  author = {Andersen, Per-Arne and Goodwin, Morten and Granmo, Ole-Christoffer},
  title = {Towards model-based reinforcement learning for industry-near environments},
  booktitle = {Artificial Intelligence XXXVI},
  editor = {Bramer, Max and Petridis, Miltos},
  series = {Lecture Notes in Computer Science},
  volume = {11927},
  pages = {36--49},
  year = {2019},
  publisher = {Springer},
  doi = {10.1007/978-3-030-34885-4_3}
}

@article{Demizu2023,
  title = {Inventory management of new products in retailers using model-based deep reinforcement learning},
  author = {Demizu, Tsukasa and Fukazawa, Yusuke and Morita, Hiroshi},
  journal = {Expert Systems with Applications},
  volume = {229},
  pages = {120256},
  year = {2023},
  doi = {10.1016/j.eswa.2023.120256}
}

@article{Zhu2025VCTS,
  title = {{IRS}-guided policy search in model based reinforcement learning for virtually coupled train sets control},
  author = {Zhu, Li and Miao, Jia and Chen, Cheng and Yan, Baicheng and Gong, Taiyuan and Yu, F. Richard and Tang, Tao},
  journal = {Computers \& Industrial Engineering},
  volume = {208},
  pages = {111372},
  year = {2025},
  doi = {10.1016/j.cie.2025.111372}
}

@article{Wang2023TrafficInference,
  title = {Model-Based Deep Reinforcement Learning with Traffic Inference for Traffic Signal Control},
  author = {Hao Wang and Jinan Zhu and Bao Gu},
  year = {2023},
  journal = {Applied Sciences},
  volume = {13},
  number = {6},
  pages = {4010},
  doi = {10.3390/app13064010}
}

@inproceedings{Li2025Think2Drive,
  title = {{Think2Drive}: Efficient Reinforcement Learning by Thinking with Latent World Model for Autonomous Driving (in {CARLA-V2})},
  author = {Qifeng Li and Xiaosong Jia and Shaobo Wang and Junchi Yan},
  year = {2025},
  booktitle = {Computer Vision -- ECCV 2024},
  series = {Lecture Notes in Computer Science},
  volume = {15103},
  pages = {142--158},
  publisher = {Springer Nature Switzerland},
  doi = {10.1007/978-3-031-72995-9_9}
}

@article{Zhu2025SLFMLight,
  author = {Zhu, Yifan and Lv, Yisheng and Lin, Shu and Xu, Jungang},
  title = {A Stochastic Traffic Flow Model-Based Reinforcement Learning Framework for Advanced Traffic Signal Control},
  journal = {IEEE Transactions on Intelligent Transportation Systems},
  year = {2025},
  volume = {26},
  number = {1},
  pages = {714--723},
  doi = {10.1109/TITS.2024.3494251}
}

@misc{Hubbs2020ORGymSoftware,
  author = {Hubbs, Christian D. and Perez, Hector D. and Sarwar, Owais and Sahinidis, Nikolaos V. and Grossmann, Ignacio E. and Wassick, John M.},
  title = {{OR-Gym}: Environments for {OR} and {RL} Research},
  year = {2020},
  howpublished = {Software repository},
  url = {https://github.com/hubbs5/or-gym},
  note = {Accessed 30 August 2026}
}

@inproceedings{Berto2025RL4CO,
  author = {Berto, Federico and Hua, Chuanbo and Park, Junyoung and Luttmann, Laurin and Ma, Yining and Bu, Fanchen and Wang, Jiarui and Ye, Haoran and Kim, Minsu and Choi, Sanghyeok and Zepeda, Nayeli Gast and Hottung, Andr{\'e} and Zhou, Jianan and Bi, Jieyi and Hu, Yu and Liu, Fei and Kim, Hyeonah and Son, Jiwoo and Kim, Haeyeon and Angioni, Davide and Kool, Wouter and Cao, Zhiguang and Zhang, Qingfu and Kim, Joungho and Zhang, Jie and Shin, Kijung and Wu, Cathy and Ahn, Sungsoo and Song, Guojie and Kwon, Changhyun and Tierney, Kevin and Xie, Lin and Park, Jinkyoo},
  title = {{RL4CO}: An Extensive Reinforcement Learning for Combinatorial Optimization Benchmark},
  booktitle = {Proceedings of the 31st ACM SIGKDD Conference on Knowledge Discovery and Data Mining V.2},
  year = {2025},
  pages = {5278--5289},
  publisher = {Association for Computing Machinery},
  doi = {10.1145/3711896.3737433},
  url = {https://github.com/ai4co/rl4co},
  note = {Software repository accessed 30 August 2026}
}

@article{Reinelt1991TSPLIB,
  author = {Reinelt, Gerhard},
  title = {{TSPLIB}---A Traveling Salesman Problem Library},
  journal = {ORSA Journal on Computing},
  year = {1991},
  volume = {3},
  number = {4},
  pages = {376--384},
  doi = {10.1287/ijoc.3.4.376},
  url = {https://comopt.ifi.uni-heidelberg.de/software/TSPLIB95/},
  note = {Benchmark website accessed 30 August 2026}
}

@article{Uchoa2017CVRPLIB,
  author = {Uchoa, Eduardo and Pecin, Diego and Pessoa, Artur and Poggi, Marcus and Vidal, Thibaut and Subramanian, Anand},
  title = {New benchmark instances for the Capacitated Vehicle Routing Problem},
  journal = {European Journal of Operational Research},
  year = {2017},
  volume = {257},
  number = {3},
  pages = {845--858},
  doi = {10.1016/j.ejor.2016.08.012},
  url = {https://galgos.inf.puc-rio.br/cvrplib/index.php/en/},
  note = {Benchmark website accessed 30 August 2026}
}

@book{Puterman1994,
  author = {Puterman, Martin L.},
  title = {Markov Decision Processes: Discrete Stochastic Dynamic Programming},
  publisher = {John Wiley \& Sons},
  year = {1994},
  doi = {10.1002/9780470316887}
}

@book{Powell2011,
  author = {Powell, Warren B.},
  title = {Approximate Dynamic Programming: Solving the Curses of Dimensionality},
  edition = {2nd},
  publisher = {John Wiley \& Sons},
  year = {2011},
  doi = {10.1002/9781118029176}
}

@book{BertsekasTsitsiklis1996,
  author = {Bertsekas, Dimitri P. and Tsitsiklis, John N.},
  title = {Neuro-Dynamic Programming},
  publisher = {Athena Scientific},
  address = {Belmont, MA},
  year = {1996}
}

@article{deFariasVanRoy2003,
  author = {De Farias, Daniela P. and Van Roy, Benjamin},
  title = {The Linear Programming Approach to Approximate Dynamic Programming},
  journal = {Operations Research},
  year = {2003},
  volume = {51},
  number = {6},
  pages = {850--865},
  doi = {10.1287/opre.51.6.850.24925}
}

@article{Adelman2007,
  author = {Adelman, Daniel},
  title = {Dynamic Bid Prices in Revenue Management},
  journal = {Operations Research},
  year = {2007},
  volume = {55},
  number = {4},
  pages = {647--661},
  doi = {10.1287/opre.1060.0368}
}

@article{Sen1999,
  author = {Sen, Suvrajeet and Higle, Julia L.},
  title = {An Introductory Tutorial on Stochastic Linear Programming Models},
  journal = {Interfaces},
  year = {1999},
  volume = {29},
  number = {2},
  pages = {33--61},
  doi = {10.1287/inte.29.2.33}
}

@article{Bayraksan2006,
  author = {Bayraksan, G{\"u}zin and Morton, David P.},
  title = {Assessing solution quality in stochastic programs},
  journal = {Mathematical Programming},
  year = {2006},
  volume = {108},
  pages = {495--514},
  doi = {10.1007/s10107-006-0720-x}
}

@article{Chang2005,
  author = {Chang, Hyeong Soo and Fu, Michael C. and Hu, Jiaqiao and Marcus, Steven I.},
  title = {An Adaptive Sampling Algorithm for Solving {Markov} Decision Processes},
  journal = {Operations Research},
  year = {2005},
  volume = {53},
  number = {1},
  pages = {126--139},
  doi = {10.1287/opre.1040.0145}
}

@article{Pereira1991,
  author = {Pereira, M. V. F. and Pinto, L. M. V. G.},
  title = {Multi-stage stochastic optimization applied to energy planning},
  journal = {Mathematical Programming},
  year = {1991},
  volume = {52},
  pages = {359--375},
  doi = {10.1007/BF01582895}
}

@article{Brown2010,
  author = {Brown, David B. and Smith, James E. and Sun, Peng},
  title = {Information Relaxations and Duality in Stochastic Dynamic Programs},
  journal = {Operations Research},
  year = {2010},
  volume = {58},
  number = {4, Part 1},
  pages = {785--801},
  doi = {10.1287/opre.1090.0796}
}

@article{Iyengar2005,
  author = {Iyengar, Garud N.},
  title = {Robust Dynamic Programming},
  journal = {Mathematics of Operations Research},
  year = {2005},
  volume = {30},
  number = {2},
  pages = {257--280},
  doi = {10.1287/moor.1040.0129}
}

@article{XuMannor2012,
  author = {Xu, Huan and Mannor, Shie},
  title = {Distributionally Robust {Markov} Decision Processes},
  journal = {Mathematics of Operations Research},
  year = {2012},
  volume = {37},
  number = {2},
  pages = {288--300},
  doi = {10.1287/moor.1120.0540}
}

@book{Bertsekas2012,
  author = {Bertsekas, Dimitri P.},
  title = {Dynamic Programming and Optimal Control},
  volume = {II},
  edition = {4th},
  year = {2012},
  publisher = {Athena Scientific},
  isbn = {9781886529441}
}

@article{Vidal2022,
  author  = {Vidal, T.},
  title   = {Hybrid genetic search for the {CVRP}: Open-source implementation and {SWAP}* neighborhood},
  journal = {Computers \& Operations Research},
  volume  = {140},
  pages   = {105643},
  year    = {2022},
  doi     = {10.1016/j.cor.2021.105643}
}

@article{Accorsi2022,
  author  = {Accorsi, L. and Lodi, A. and Vigo, D.},
  title   = {Guidelines for the computational testing of machine learning approaches to vehicle routing problems},
  journal = {Operations Research Letters},
  volume  = {50},
  number  = {2},
  pages   = {229--234},
  year    = {2022},
  doi ={10.1016/j.orl.2022.01.018}
}

\end{document}